 \input amstexl


\catcode`\@=11
\ifx\amstexloaded@\relax\else
 \errmessage{AmS-TeX must be loaded before LamS-TeX}\fi
\ifx\laxread@\undefined\else\catcode`\@=\active \fi
\def\err@#1{\errmessage{LamS-TeX error: #1}}
\def^^L{\par}
\let\+\tabalign
\def\newcount{\alloc@0\count\countdef\insc@unt}
\def\newdimen{\alloc@1\dimen\dimendef\insc@unt}
\def\newskip{\alloc@2\skip\skipdef\insc@unt}
\def\newmuskip{\alloc@3\muskip\muskipdef\@cclvi}
\def\newbox{\alloc@4\box\chardef\insc@unt}
\let\newtoks\relax
\def\newhelp#1#2{\newtoks#1#1\expandafter{\csname#2\endcsname}}
\def\newtoks{\alloc@5\toks\toksdef\@cclvi}
\def\newread{\alloc@6\read\chardef\sixt@@n}
\def\newwrite{\alloc@7\write\chardef\sixt@@n}
\def\newfam{\alloc@8\fam\chardef\sixt@@n}
\def\newlanguage{\alloc@9\language\chardef\@cclvi}
\def\newinsert#1{\global\advance\insc@unt by\m@ne
  \ch@ck0\insc@unt\count
  \ch@ck1\insc@unt\dimen
  \ch@ck2\insc@unt\skip
  \ch@ck4\insc@unt\box
  \allocationnumber=\insc@unt
  \global\chardef#1=\allocationnumber
  \wlog{\string#1=\string\insert\the\allocationnumber}}
\def\newif#1{\count@\escapechar \escapechar\m@ne
  \expandafter\expandafter\expandafter
   \edef\@if#1{true}{\let\noexpand#1=\noexpand\iftrue}%
  \expandafter\expandafter\expandafter
   \edef\@if#1{false}{\let\noexpand#1=\noexpand\iffalse}%
  \@if#1{false}\escapechar\count@}

\def\Err@#1{\errhelp\defaulthelp@\err@{#1}}
{\catcode`\@=\active
 \edef\next{\gdef\noexpand@{\futurelet\noexpand\next
  \csname at\string@\endcsname}}
 \next
}
\def\at@{\ifcat\noexpand\next a\let\next@\at@@\else
 \ifcat\noexpand\next0\let\next@\at@@\else
 \ifcat\noexpand\next\relax\let\next@\at@@\else
 \let\next@\at@@@\fi\fi\fi\next@}
\def\at@@@{\errhelp\athelp@\err@{Invalid use of @}}
\def\at@@#1{\expandafter
 \ifx\csname\string#1@at\endcsname\relax\let\next@\at@@@\else
 \DN@{\csname\string#1@at\endcsname}\fi\next@}
\def\atdef@#1{\expandafter\def\csname\string#1@at\endcsname}
\newif\iftest@
\def\tagin@#1{\tagin@false
 \DN@##1\tag##2##3\next@{\test@true\ifx\tagin@##2\test@false\fi}%
 \next@#1\tag\tagin@\next@\tagin@false\iftest@\tagin@true\fi}
\let\lkerns@\relax
\def\nolinebreak{\RIfM@\mathmodeerr@\nolinebreak\else
 \ifhmode\saveskip@\lastskip\unskip
 \nobreak\ifdim\saveskip@>\z@\hskip\saveskip@\fi\lkerns@
 \else\vmodeerr@\nolinebreak\fi\fi}
\def\allowlinebreak{\RIfM@\mathmodeerr@\allowlinebreak\else
 \ifhmode\saveskip@\lastskip\unskip
 \allowbreak\ifdim\saveskip@>\z@\hskip\saveskip@\fi\lkerns@
 \else\vmodeerr@\allowlinebreak\fi\fi}
\def\linebreak{\RIfM@\mathmodeerr@\linebreak\else
 \ifhmode\unskip\unkern\break\lkerns@
 \else\vmodeerr@\linebreak\fi\fi}
\let\nkerns@\relax
\def\newline{\RIfM@\mathmodeerr@\newline\else
 \ifhmode\unskip\unkern\null\hfill\break\nkerns@
 \else\vmodeerr@\newline\fi\fi}%
\def\newbox@{\alloc@@4\box\chardef\insc@unt}
\def\newcount@{\alloc@@0\count\countdef\insc@unt}
\def\accentedsymbol#1#2{\expandafter\newbox@\csname\exstring@#1@box\endcsname
 \setbox\csname\exstring@#1@box\endcsname\hbox{$\m@th#2$}%
 \define#1{\copy\csname\exstring@#1@box\endcsname{}}}
\def\rightadd@#1\to#2{\toks@{\\#1}\toks@@\expandafter{#2}\xdef#2{\the\toks@@
 \the\toks@}\toks@{}\toks@@{}}
\def\fontlist@{\\\tenrm\\\sevenrm\\\fiverm\\\teni\\\seveni\\\fivei
 \\\tensy\\\sevensy\\\fivesy\\\tenex\\\tenbf\\\sevenbf\\\fivebf
 \\\tensl\\\tenit}
\def\font@#1=#2 {\rightadd@#1\to\fontlist@\font#1=#2 }
\def\ismember@#1#2{\global\let\Next@ F\let\next@= #2%
 {\def\\##1{\let\nextii@##1\ifx\nextii@\next@\global\let\Next@ T\fi}#1}%
 \test@false\ifx\Next@ T\test@true\fi\let\next@\relax}
\def\FNSS@#1{\let\FNSS@@#1\FN@\FNSS@@@}
\def\FNSS@@@{\ifx\next\space@\def\FNSS@@@@. {\FN@\FNSS@@@}\else
 \def\FNSS@@@@.{\FNSS@@}\fi\FNSS@@@@.}
\atdef@"{\unskip
 \DN@{\ifx\next`\DN@`{\FN@\nextii@}%
  \else\ifx\next\lq\DN@\lq{\FN@\nextii@}%
  \else\DN@####1{\FN@\nextiii@}\fi\fi
  \next@}%
 \DNii@{\ifx\next`\DN@`{\sldl@``}%
  \else\ifx\next\lq\DN@\lq{\sldl@``}%
  \else\DN@{\dlsl@`}\fi\fi\next@}%
 \def\nextiii@{\ifx\next'\DN@'{\srdr@''}%
  \else\ifx\next\rq\DN@\rq{\srdr@''}%
  \else\DN@{\drsr@'}\fi\fi\next@}%
 \FNSS@\next@}
\def\root{%
  \DN@{\ifx\next\uproot\let\next@\nextii@\else
   \ifx\next\leftroot\let\next@\nextiii@\else
   \let\next@\plainroot@\fi\fi\next@}%
  \DNii@\uproot##1{\uproot@##1\relax\FNSS@\nextiv@}%
  \def\nextiv@{\ifx\next\leftroot\let\next@\nextv@\else
   \let\next@\plainroot@\fi\next@}%
  \def\nextv@\leftroot##1{\leftroot@##1\relax\plainroot@}%
  \def\nextiii@\leftroot##1{\leftroot@##1\relax\FNSS@\nextvi@}%
  \def\nextvi@{\ifx\next\uproot\let\next@\nextvii@\else
   \let\next@\plainroot@\fi\next@}%
  \def\nextvii@\uproot##1{\uproot@##1\relax\plainroot@}%
  \bgroup\uproot@\z@\leftroot@\z@
 \FNSS@\next@}
\def\loop#1\repeat{\def\iterate{#1\relax\expandafter\iterate\fi}%
 \iterate\let\iterate\relax}
\def\gloop@#1\repeat{\gdef\iterate@{#1\relax\expandafter\iterate@\fi}%
 \iterate@\global\let\iterate@\relax}
\def\printoptions{\W@{Do you want S(yntax check),
  G(alleys) or P(ages)?^^JType S, G or P, follow by <return>: }\loop
 \read\m@ne to\ans@
 \edef\next@{\def\noexpand\Ans@{\ans@}}\uppercase\expandafter{\next@}%
 \ifx\Ans@\S@\test@true\syntax\else
 \ifx\Ans@\G@\test@true\galleys\else
 \ifx\Ans@\P@\test@true\else
 \test@false\fi\fi\fi
 \iftest@\else\W@{Type S, G or P, follow by <return>: }%
 \repeat}
\expandafter\let\csname A@;\endcsname;
\expandafter\let\csname A@:\endcsname:
\expandafter\let\csname A@?\endcsname?
\expandafter\let\csname A@!\endcsname!
\def\APdef#1{\def\next@{\expandafter\let\csname A@\string#1\endcsname#1}%
 \afterassignment\next@\def#1}
\let\fextra@\,
\def\tdots@{\unskip
 \DN@{$\m@th\mathinner{\ldotp\ldotp\ldotp}\,
   \ifx\next,\,$\else\ifx\next.\,$\else
   \ifx\next;\,$\else
   \expandafter\ifx\csname A@\string;\endcsname\next\fextra@$\else
   \ifx\next:\,$\else
   \expandafter\ifx\csname A@\string:\endcsname\next\fextra@$\else
   \ifx\next?\,$\else
   \expandafter\ifx\csname A@\string?\endcsname\next\fextra@$\else
   \ifx\next!\,$\else
   \expandafter\ifx\csname A@\string!\endcsname\next\fextra@$\else
   $ \fi\fi\fi\fi\fi\fi\fi\fi\fi\fi}%
 \ \FN@\next@}
\def\extrap@#1{%
 \ifx\next,\DN@{#1\,}\else
 \ifx\next;\DN@{#1\,}\else
 \expandafter\ifx\csname A@\string;\endcsname\next\DN@{#1\fextra@}\else
 \ifx\next.\DN@{#1\,}\else\extra@
 \ifextra@\DN@{#1\,}\else
 \let\next@#1\fi\fi\fi\fi\fi\next@}
\def\dotsc{\DN@{\ifx\next;\plainldots@\,\else
 \expandafter\ifx\csname A@\string;\endcsname\next\plainldots@\fextra@\else
 \ifx\next.\plainldots@\,\else\extra@\plainldots@
 \ifextra@\,\fi\fi\fi\fi}%
 \FN@\next@}
\def\keybin@{\keybin@true
 \ifx\next+\else\ifx\next=\else\ifx\next<\else\ifx\next>\else\ifx\next-\else
 \ifx\next*\else\ifx\next:\else
 \expandafter\ifx\csname A@\string;\endcsname\next\else
 \keybin@false\fi\fi\fi\fi\fi\fi\fi\fi}
\def\boldkey#1{\ifcat\noexpand#1A%
  \ifcmmibloaded@{\fam\cmmibfam#1}\else
   \Err@{First bold symbol font not loaded}\fi
 \else
 \let\next=#1%
 \ifx#1!\mathchar"5\bffam@21 \else
 \expandafter\ifx\csname A@\string!\endcsname\next\mathchar"5\bffam@21 \else
 \ifx#1(\mathchar"4\bffam@28 \else\ifx#1)\mathchar"5\bffam@29 \else
 \ifx#1+\mathchar"2\bffam@2B \else\ifx#1:\mathchar"3\bffam@3A \else
 \expandafter\ifx\csname A@\string:\endcsname\next\mathchar"3\bffam@3A \else
 \ifx#1;\mathchar"6\bffam@3B \else
 \expandafter\ifx\csname A@\string;\endcsname\next\mathchar"6\bffam@3B \else
 \ifx#1=\mathchar"3\bffam@3D \else
 \ifx#1?\mathchar"5\bffam@3F \else
 \expandafter\ifx\csname A@\string?\endcsname\next\mathchar"5\bffam@3F \else
 \ifx#1[\mathchar"4\bffam@5B \else
 \ifx#1]\mathchar"5\bffam@5D \else
 \ifx#1,\mathchari@63B \else
 \ifx#1-\mathcharii@200 \else
 \ifx#1.\mathchari@03A \else
 \ifx#1/\mathchari@03D \else
 \ifx#1<\mathchari@33C \else
 \ifx#1>\mathchari@33E \else
 \ifx#1*\mathcharii@203 \else
 \ifx#1|\mathcharii@06A \else
 \ifx#10\bold0\else\ifx#11\bold1\else\ifx#12\bold2\else\ifx#13\bold3\else
 \ifx#14\bold4\else\ifx#15\bold5\else\ifx#16\bold6\else\ifx#17\bold7\else
 \ifx#18\bold8\else\ifx#19\bold9\else
  \Err@{\noexpand\boldkey can't be used with #1}%
 \fi\fi\fi\fi\fi\fi\fi\fi\fi\fi\fi\fi\fi\fi\fi
 \fi\fi\fi\fi\fi\fi\fi\fi\fi\fi\fi\fi\fi\fi\fi\fi\fi\fi}
\def\arabic#1{#1}
\def\alph#1{\count@#1\relax\advance\count@96 \ifnum\count@>122
 \Err@{\noexpand\alph invalid for numbers > 26}\else\char\count@\fi}
\def\Alph#1{\count@#1\relax\advance\count@64 \ifnum\count@>90
 \Err@{\noexpand\Alph invalid for numbers > 26}\else\char\count@\fi}

\def\Roman#1{\uppercase\expandafter{\romannumeral#1}}
\def\fnsymbol#1{\count@#1\relax
 \count@@\count@
 \advance\count@\m@ne\divide\count@7
 \count@@@\count@\advance\count@@@\@ne
 \multiply\count@7 \advance\count@@-\count@
 \count@\count@@@
 {\loop
  \ifcase\count@@\or*\or\dag\or\ddag\or\P\or\S\or\text{$\|$}\or\#\fi
  \advance\count@\m@ne\ifnum\count@>\z@\repeat}}
\def\cardnine@#1{\ifcase#1\or one\or two\or three\or four\or five\or
 six\or seven\or eight\or nine\fi}
\let\alloc@\alloc@@
\newcount\ten@
\ten@10
\def\cardinal#1{\count@#1\relax
 \ifnum\count@>99 \number\count@
 \else
  \ifnum\count@=\z@ zero%
  \else
   \ifnum\count@<\ten@\cardnine@\count@
   \else
    \ifnum\count@<20
     \advance\count@-\ten@
     \ifcase\count@ ten\or eleven\or twelve\or thirteen\or fourteen\or
      fifteen\or sixteen\or seventeen\or eighteen\or nineteen\fi
    \else
     \count@@\count@\count@@@\count@@
     \divide\count@\ten@\multiply\count@\ten@
     \advance\count@@@-\count@\divide\count@\ten@
     \ifcase\count@\or\or twenty\or thirty\or forty\or fifty\or sixty\or
      seventy\or eighty\or ninety\fi
     \ifnum\count@@@=\z@\else-\cardnine@\count@@@\fi
    \fi
   \fi
  \fi
 \fi}
\def\ordnine@#1{\ifcase#1\or first\or second\or third\or fourth\or fifth\or
 sixth\or seventh\or eighth\or ninth\fi}
\newcount\count@@@@
\def\ordsuffix@{\count@@@@\count@
 \divide\count@\ten@
 \count@@@\count@\count@@\count@
 \divide\count@@\ten@\multiply\count@@\ten@
 \advance\count@@@-\count@@
 \ifnum\count@@@=\@ne th%
 \else
  \count@@@\count@@@@
  \count@@\count@@@@
  \divide\count@@\ten@\multiply\count@@\ten@
  \advance\count@@@-\count@@
  \ifcase\count@@@ th\or st\or nd\or rd\else th\fi
 \fi}
\def\nordinal#1{\count@#1\relax\number\count@\ordsuffix@}
\def\spordinal#1{\count@#1\relax\number\count@$^{\text{\ordsuffix@}}$}
\def\ordinal#1{\count@#1\relax
 \ifnum\count@>99 \number\count@\ordsuffix@
 \else
   \ifnum\count@=\z@ zeroth%
  \else
    \ifnum\count@<\ten@\ordnine@\count@
    \else
     \ifnum\count@<20 \advance\count@-\ten@
      \ifcase\count@ tenth\or eleventh\or twelfth\or thirteenth\or
       fourteenth\or fifteenth\or sixteenth\or seventeenth\or eighteenth\or
       nineteenth\fi
     \else
      \count@@\count@
      \divide\count@\ten@\multiply\count@\ten@
      \count@@@\count@@\advance\count@@@-\count@
      \divide\count@\ten@
      \ifcase\count@\or\or twent\or thirt\or fort\or fift\or sixt\or sevent\or
       eight\or ninet\fi
      \ifnum\count@@@=\z@ ieth\else y-\ordnine@\count@@@\fi
     \fi
    \fi
  \fi
 \fi}
\font@\tensmc=cmcsc10
\textonlyfont@\smc\tensmc
\newtoks\noexpandtoks@
\noexpandtoks@{\let\arabic\relax\let\alph\relax\let\Alph\relax
 \let\Roman\relax\let\fnsymbol\relax\let\rm\relax
 \let\it\relax\let\bf\relax\let\sl\relax\let\smc\relax
 \let\/\relax\let\null\relax}
\def\noexpands@{\the\noexpandtoks@}
\def\Nonexpanding#1{\global\noexpandtoks@
 \expandafter{\the\noexpandtoks@\let#1\relax}}
\def\prevanish@{\saveskip@\z@\ifhmode\saveskip@\lastskip\unskip\fi}
\def\postvanish@{\ifdim\saveskip@>\z@\hskip\saveskip@\fi\FN@\postvanish@@}
\def\postvanish@@{\DN@.{}%
 \ifx\next\space@\ifdim\saveskip@>\z@\DN@. {}\fi\fi\next@.}
\def\invisible#1{\prevanish@\ignorespaces#1\unskip\postvanish@}
\def\vanishlist@{\\\invisible}
\let\noindent@\noindent
\def\noindent{\par\noindent@\FN@\pretendspace@}
\def\pretendspace@{\ismember@\vanishlist@\next
 \iftest@\nobreak\hskip-\p@\hskip\p@\fi}

\newtoks\everypartoks@
\def\noindent@@{\par\everypartoks@\expandafter{\the\everypar}\everypar{}%
 \noindent@\everypar\expandafter{\the\everypartoks@}}
\def\page{\Err@{\noexpand\page has no meaning by itself}}
\let\page@C\pageno
\let\page@P\empty
\let\page@Q\empty
\def\page@S#1{#1\/}
\def\page@F{\rm}
\def\page@N{\arabic}   
\newif\ifindexing@
\def\indexfile{\ifindexing@\else
 \alloc@@7\write\chardef\sixt@@n\ndx@
 \immediate\openout\ndx@=\jobname.ndx
 \global\indexing@true\fi}
\global\advance\insc@unt\m@ne
\ch@ck0\insc@unt\count
\ch@ck1\insc@unt\dimen
\ch@ck2\insc@unt\skip
\ch@ck4\insc@unt\box
\allocationnumber\insc@unt
\global\chardef\margin@\allocationnumber
\dimen\margin@\maxdimen
\count\margin@\z@
\skip\margin@\z@
\newif\ifindexproofing@
\def\indexproofing{\indexproofing@true}
\def\noindexproofing{\indexproofing@false}
\def\unmacro@#1:#2->#3\unmacro@{\def\macpar@{#2}\def\macdef@{#3}}
\def\starparts@#1{\def\stari@{#1}\def\starii@{#1}\let\stariii@\empty
 \test@false
 \DN@##1*##2##3\next@{\ifx\starparts@##2\test@false\else\test@true\fi}%
 \next@#1*\starparts@\next@
 \iftest@\DN@{\starparts@@#1\starparts@@}\else\let\next@\relax\fi\next@}
\def\starparts@@#1*#2\starparts@@{\def\starii@{#1}\def\stariii@{*#2}}
\def\windex@{\ifindexing@
 \expandafter\unmacro@\meaning\stari@\unmacro@
 \edef\macdef@{\string"\macdef@\string"}%
 \edef\next@{\write\ndx@{\macdef@}}\next@
 \write\ndx@{{\number\pageno}{\page@N}{\page@P}{\page@Q}}%
 \fi
 \ifindexproofing@
  \ifx\stariii@\empty\else
   \expandafter\unmacro@\meaning\stariii@\unmacro@\fi
  \insert\margin@{\hbox{\rm\vrule\height9\p@\depth2\p@\width\z@\starii@
  \ifx\stariii@\empty\else\tt\macdef@\fi}}\fi}
\catcode`\"=\active
\def"{\FN@\quote@}
\def\quote@{\ifx\next"\expandafter\quote@@\else\expandafter\quote@@@\fi}
\def\quote@@@#1"{\starparts@{#1}\starii@\windex@}
\def\quote@@"#1"{\prevanish@\starparts@{#1}\windex@\FN@\quote@@@@}
\def\quote@@@@{\ifx\next"\DN@"{\postvanish@}\else
 \let\next@\postvanish@\fi\next@}
\rightadd@"\to\vanishlist@
\def\idefine#1{\DN@{#1}\DNii@{\noexpand#1}%
 \afterassignment\idefine@\def\nextiii@}
\def\idefine@{\ifindexing@
 \expandafter\let\next@\nextiii@
 \expandafter\unmacro@\meaning\nextiii@\unmacro@
 \immediate\write\ndx@{\noexpand\define\nextii@\macpar@{\macdef@}}\fi}
\def\iabbrev*#1#2{\ifindexing@\toks@{#2}%
 \immediate\write\ndx@{\noexpand\abbrev*\noexpand#1{\the\toks@}}\fi}
\newread\laxread@
\newwrite\laxwrite@
\let\fnpages@\empty
\def\Finit@#1#2\Finit@{\let\nextii@#1\def\nextiii@{#2}}
\catcode`\~=11
\def\getparts@ @#1~#2~#3~#4~#5~#6{\def\nextiv@{#1}%
 \def\nextiii@{#2~#3~#4~#5~}\count@#6\relax}
\newif\ifdocument@
\def\document{\ifdocument@\else\global\document@true
 \let\fontlist@\empty
 \immediate\openin\laxread@=\jobname.lax\relax
 {\endlinechar\m@ne\noexpands@\catcode`\@=11 \catcode`\~=11
  \loop\ifeof\laxread@\else
   \read\laxread@ to\next@
   \ifx\next@\empty
   \else
    \expandafter\Finit@\next@\Finit@
    \if\nextii@ F%
     \expandafter\rightadd@\nextiii@\to\fnpages@
    \else
     \expandafter\getparts@\next@
     \edef\next@{\gdef\csname\nextiv@ @L\endcsname{\nextiii@\number\count@}}%
     \next@
    \fi
   \fi
  \repeat}%
 \immediate\closein\laxread@
 \immediate\openout\laxwrite@=\jobname.lax\relax\fi}
\let\thelabel@\relax
\def\thelabels@{\thelabel@ ~\thelabel@@ ~\thelabel@@@ ~\thelabel@@@@ ~}
\def\label#1{\prevanish@
 \ifx\thelabel@\relax
  \Err@{There's nothing here to be labelled}%
 \else
  {\noexpands@
  \expandafter\ifx\csname#1@L\endcsname\relax
   \expandafter\xdef\csname#1@L\endcsname{\thelabels@0}%
   \immediate\write\laxwrite@{@#1~\thelabels@1}%
  \else
   \edef\next@{@~\csname#1@L\endcsname}%
    \expandafter\getparts@\next@
    \ifodd\count@
    \expandafter\xdef\csname#1@L\endcsname{\thelabels@0}%
    \immediate\write\laxwrite@{@#1~\thelabels@1}%
   \else
    \Err@{Label #1 already used}%
   \fi
  \fi
  }%
 \fi
 \postvanish@}
\rightadd@\label\to\vanishlist@
\def\thepages@{\page@N{\number\page@C}~%
 \page@S{\page@P\page@N{\number\page@C}\page@Q}~%
 \number\page@C ~\page@P\page@N{\number\page@C}\page@Q ~}
\def\pagelabel#1{\prevanish@
 \expandafter\ifx\csname#1@L\endcsname\relax
  {\noexpands@
  \expandafter\xdef\csname#1@L\endcsname{\thepages@2}}%
  \write\laxwrite@{@#1~\thepages@3}%
 \else
  {\noexpands@
  \edef\next@{@~\csname#1@L\endcsname}%
  \expandafter\getparts@\next@
  \ifodd\count@
   \ifnum\count@=\@ne
    \expandafter\xdef\csname#1@L\endcsname{\thelabels@2}%
   \fi
   \write\laxwrite@{@#1~\thepages@3}%
  \else
   \Err@{Label #1 already used}%
  \fi
  }%
 \fi
 \postvanish@}
\rightadd@\pagelabel\to\vanishlist@
\newif\ifreferr@
\referr@true
\def\RefErrors{\global\referr@true}
\def\RefWarnings{\global\referr@false}
\setbox\z@\hbox{\global\count@=`^^30}
\ifnum\count@=48 \let\versionthree@\relax\fi
\def\nolabel@#1#2#3{\expandafter\ifx\csname#2@L\endcsname\relax
 \ifreferr@\Err@{No \noexpand\label found for #2}\else
 \W@{Warning: No \noexpand\label found for #2.}%
 \ifx\versionthree@\relax\W@{l.\number\inputlineno\space ... \string#1{#2}}\fi
 \fi#3\else}
\def\csL@#1{{\noexpands@\xdef\Next@{\csname#1@L\endcsname}}}
\def\ref#1{\nolabel@\ref{#1}\relax
 \DNii@##1~##2\nextii@{##1}%
 \csL@{#1}\expandafter\nextii@\Next@\nextii@\fi}
\def\Ref#1{\nolabel@\Ref{#1}\relax
 \DNii@##1~##2~##3\nextii@{##2}%
 \csL@{#1}\expandafter\nextii@\Next@\nextii@\fi}
\def\nref#1{\nolabel@\nref{#1}\relax
 \DNii@##1~##2~##3~##4\nextii@{##3}%
 \csL@{#1}\expandafter\nextii@\Next@\nextii@\fi}
\def\pref#1{\nolabel@\pref{#1}\relax
 \DNii@##1~##2~##3~##4~##5\nextii@{##4}%
 \csL@{#1}\expandafter\nextii@\Next@\nextii@\fi}
\let\pref@\pref
\def\Evaluatenref#1{\nolabel@\Evaluatenref{#1}{\gdef\Nref{-10000 }}%
 \DNii@##1~##2~##3~##4\nextii@{\DNii@{##3}}%
 \csL@{#1}\expandafter\nextii@\Next@\nextii@
 \xdef\Nref{\nextii@}\fi}
\def\Evaluatepref#1{\nolabel@\Evaluatepref{#1}{\global\let\Pref\empty}%
 \DNii@##1~##2~##3~##4~##5\nextii@{\DNii@{##4}}%
 \csL@{#1}\expandafter\nextii@\Next@\nextii@
 \xdef\Pref{\nextii@}\fi}
\def\readlax#1{\immediate\openin\laxread@=#1.lax\relax
 \ifeof\laxread@\W@{}\W@{File #1.lax not found.}\W@{}\fi
 {\endlinechar\m@ne\noexpands@\catcode`\@=11 \catcode`\~=11
  \loop\ifeof\laxread@\else
   \read\laxread@ to\nextv@
   \ifx\nextv@\empty
   \else
    \expandafter\Finit@\nextv@\Finit@
    \ifx\nextii@ F%
    \else
     \expandafter\getparts@\nextv@
     \expandafter\ifx\csname\nextiv@ @L\endcsname\relax
      \edef\next@{\gdef\csname\nextiv@ @L\endcsname
       {\nextiii@\ifnum\count@=\@ne0\else2\fi}}%
      \next@
     \else
      \Err@{Label \nextiv@\space in #1.lax already used}%
     \fi
    \fi
   \fi
  \repeat}%
 \immediate\closein\laxread@}
\catcode`\~=\active

\def\FNSSP@{\FNSS@\pretendspace@}
\everydisplay{\csname displaymath \endcsname}
\expandafter\def\csname displaymath \endcsname#1$${#1$$\FNSSP@}
\def\locallabel@{\let\thelabel@\Thelabel@\let\thelabel@@\Thelabel@@
 \let\thelabel@@@\Thelabel@@@\let\thelabel@@@@\Thelabel@@@@}
\newcount\tag@C
\tag@C\z@
\let\tag@P\empty
\let\tag@Q\empty
\def\tag@S#1{{\rm(}{#1\/}{\rm)}}
\let\tag@N\arabic
\def\tag@F{\rm}
\def\maketag@{\FN@\maketag@@}
\def\maketag@@{\ifx\next\relax\DN@\relax{\FN@\maketag@@}\else
 \ifx\next"\let\next@\maketag@@@\else
 \let\next@\maketag@@@@\fi\fi\next@}
\def\xdefThelabel@#1{\xdef\Thelabel@{#1{\Thelabel@@@}}}
\def\xdefThelabel@@#1{\xdef\Thelabel@@{#1{\Thelabel@@@@}}}
\def\maketag@@@@#1\maketag@{\global\advance\tag@C\@ne
 {\noexpands@
  \xdef\Thelabel@@@{\number\tag@C}%
  \xdefThelabel@\tag@N
  \xdef\Thelabel@@@@{\ifmathtags@$\tag@P\Thelabel@\tag@Q$\else
   \tag@P\Thelabel@\tag@Q\fi}%
  \xdefThelabel@@\tag@S
  }%
 \locallabel@
 \hbox{\tag@F\thelabel@@}%
 #1}
\def\Qlabel@#1{{\noexpands@\xdef\Thelabel@@{#1}%
 \let\style\empty\xdef\Thelabel@@@@{#1}%
 \let\pre\empty\let\post\empty\xdef\Thelabel@{#1}%
 \let\numstyle\empty\xdef\Thelabel@@@{#1}}}
\def\maketag@@@"#1"#2\maketag@{%
 {\let\pre\tag@P\let\post\tag@Q\let\style\tag@S\let\numstyle\tag@N
  \hbox{\tag@F#1}%
  \noexpands@
  \Qlabel@{#1}%
  }%
 \locallabel@
 #2}
\def\align@{\inalign@true\inany@true
 \vspace@\allowdisplaybreak@\displaybreak@\intertext@
 \def\tag{\global\tag@true\ifnum\and@=\z@
  \DN@{&\omit\global\rwidth@\z@&\relax}\else
  \DN@{&\relax}\fi\next@}%
 \iftagsleft@\DN@{\csname align \endcsname}\else
  \DN@{\csname align \space\endcsname}\fi\next@}
\def\noset@{\def\Offset##1##2{\prevanish@\postvanish@}%
 \def\Reset##1##2{\prevanish@\postvanish@}}
\def\measure@#1\endalign{\global\lwidth@\z@\global\rwidth@\z@
 \global\maxlwidth@\z@\global\maxrwidth@\z@
 \global\and@\z@
 \setbox\z@\vbox
  {\noset@\everycr{\noalign{\global\tag@false\global\and@\z@}}\Let@
  \halign{\setboxz@h{$\m@th\displaystyle{\@lign##}$}%
   \global\lwidth@\wdz@
   \ifdim\lwidth@>\maxlwidth@\global\maxlwidth@\lwidth@\fi
   \global\advance\and@\@ne
   &\setboxz@h{$\m@th\displaystyle{{}\@lign##}$}\global\rwidth@\wdz@
   \ifdim\rwidth@>\maxrwidth@\global\maxrwidth@\rwidth@\fi
   \global\advance\and@\@ne
   &\Tag@\eat@{##}\crcr#1\crcr}}%
 \totwidth@\maxlwidth@\advance\totwidth@\maxrwidth@}
\def\prepost@{\global\let\tag@P@\tag@P\global\let\tag@Q@\tag@Q}
\def\reprepost@{\let\tag@P\tag@P@\let\tag@Q\tag@Q@}
\expandafter\def\csname align \space\endcsname#1\endalign
 {\measure@#1\endalign\global\and@\z@
 \ifingather@\everycr{\noalign{\global\and@\z@}}\else\displ@y@\fi
 \Let@\tabskip\centering@
 \halign to\displaywidth
  {\hfil\strut@\setboxz@h{$\m@th\displaystyle{\@lign##\prepost@}$}%
  \boxz@\global\advance\and@\@ne
  \tabskip\z@skip
  &\setboxz@h{$\m@th\displaystyle{{}\@lign##\prepost@}$}%
  \global\rwidth@\wdz@\boxz@\hfil\global\advance\and@\@ne
  \tabskip\centering@
  &\setboxz@h{\@lign\strut@\reprepost@\maketag@##\maketag@}%
  \dimen@\displaywidth\advance\dimen@-\totwidth@
  \divide\dimen@\tw@\advance\dimen@\maxrwidth@\advance\dimen@-\rwidth@
  \ifdim\dimen@<\tw@\wdz@\llap{\vtop{\normalbaselines\null\boxz@}}%
  \else\llap{\boxz@}\fi
  \tabskip\z@skip
  \crcr#1\crcr
  \black@\totwidth@}}
\expandafter\def\csname align \endcsname#1\endalign{\measure@#1\endalign
 \global\and@\z@
 \ifdim\totwidth@>\displaywidth\let\displaywidth@\totwidth@\else
  \let\displaywidth@\displaywidth\fi
 \ifingather@\everycr{\noalign{\global\and@\z@}}\else\displ@y@\fi
 \Let@\tabskip\centering@\halign to\displaywidth
  {\hfil\strut@\setboxz@h{$\m@th\displaystyle{\@lign##\prepost@}$}%
  \global\lwidth@\wdz@\global\lineht@\ht\z@
  \boxz@\global\advance\and@\@ne
  \tabskip\z@skip&\setboxz@h{$\m@th\displaystyle{{}\@lign##\prepost@}$}%
  \ifdim\ht\z@>\lineht@\global\lineht@\ht\z@\fi
  \boxz@\hfil\global\advance\and@\@ne
  \tabskip\centering@&\kern-\displaywidth@
  \setboxz@h{\@lign\strut@\reprepost@\maketag@##\maketag@}%
  \dimen@\displaywidth\advance\dimen@-\totwidth@
  \divide\dimen@\tw@\advance\dimen@\maxlwidth@\advance\dimen@-\lwidth@
  \ifdim\dimen@<\tw@\wdz@
   \rlap{\vbox{\normalbaselines\boxz@\vbox to\lineht@{}}}\else
   \rlap{\boxz@}\fi
  \tabskip\displaywidth@\crcr#1\crcr\black@\totwidth@}}
\def\attag@#1{\let\Maketag@\maketag@\let\TAG@\Tag@
 \let\Prepost@\prepost@\let\Reprepost@\reprepost@
 \let\Tag@\relax\let\maketag@\relax
 \let\prepost@\relax\let\reprepost@\relax
 \ifmeasuring@
  \def\llap@##1{\setboxz@h{##1}\hbox to\tw@\wdz@{}}%
  \def\rlap@##1{\setboxz@h{##1}\hbox to\tw@\wdz@{}}%
 \else\let\llap@\llap\let\rlap@\rlap\fi
 \toks@{\hfil\strut@
  $\m@th\displaystyle{\@lign\the\hashtoks@\prepost@}$%
  \tabskip\z@skip\global\advance\and@\@ne&
  $\m@th\displaystyle{{}\@lign\the\hashtoks@\prepost@}$\hfil
  \ifxat@\tabskip\centering@\fi\global\advance\and@\@ne}%
 \iftagsleft@
  \toks@@{\tabskip\centering@&\Tag@\kern-\displaywidth
   \rlap@{\@lign\reprepost@\maketag@\the\hashtoks@\maketag@}%
   \global\advance\and@\@ne\tabskip\displaywidth}\else
  \toks@@{\tabskip\centering@&\Tag@\llap@{\@lign\reprepost@\maketag@
   \the\hashtoks@\maketag@}\global\advance\and@\@ne\tabskip\z@skip}\fi
 \atcount@#1\relax\advance\atcount@\m@ne
 \loop\ifnum\atcount@>\z@
  \toks@\expandafter{\the\toks@&\hfil$\m@th\displaystyle{\@lign
  \the\hashtoks@\prepost@}$\global\advance\and@\@ne
  \tabskip\z@skip
  &$\m@th\displaystyle{{}\@lign\the\hashtoks@\prepost@}$\hfil\ifxat@
  \tabskip\centering@\fi\global\advance\and@\@ne}\advance\atcount@\m@ne
 \repeat
 \edef\preamble@{\the\toks@\the\toks@@}%
 \edef\preamble@@{\preamble@}%
 \let\maketag@\Maketag@\let\Tag@\TAG@
 \let\prepost@\Prepost@\let\reprepost@\Reprepost@}
\def\unlabel@{\def\label##1{\prevanish@\postvanish@}%
 \def\pagelabel##1{\prevanish@\postvanish@}}
\newcount\tag@CC
\expandafter\def\csname alignat \endcsname#1#2\endalignat
 {\inany@true\xat@false
 \def\tag{\global\tag@true
  \count@#1\relax\multiply\count@\tw@\advance\count@\m@ne
  \gdef\tag@{&}%
  \loop\ifnum\count@>\and@\xdef\tag@{&\omit\tag@}%
  \advance\count@\m@ne\repeat
  \tag@\relax}%
 \vspace@\allowdisplaybreak@\displaybreak@\intertext@
 \displ@y@\measuring@true\tag@CC\tag@C
 \setbox\savealignat@\hbox{\noset@\unlabel@$\m@th\displaystyle\Let@
  \attag@{#1}\vbox{\halign{\span\preamble@@\crcr#2\crcr}}$}%
 \measuring@false
 \Let@\attag@{#1}\tag@C\tag@CC
 \tabskip\centering@\halign to\displaywidth
  {\span\preamble@@\crcr#2\crcr\black@{\wd\savealignat@}}}
\expandafter\def\csname xalignat \endcsname#1#2\endxalignat
 {\inany@true\xat@true
 \def\tag{\global\tag@true
  \count@#1\relax\multiply\count@\tw@\advance\count@\m@ne
  \gdef\tag@{&}%
  \loop\ifnum\count@>\and@\xdef\tag@{&\omit\tag@}%
  \advance\count@\m@ne\repeat
  \tag@\relax}%
 \vspace@\allowdisplaybreak@\displaybreak@\intertext@
 \displ@y@\measuring@true\tag@CC\tag@C
 \setbox\savealignat@\hbox{\noset@\unlabel@$\m@th\displaystyle\Let@
  \attag@{#1}\vbox{\halign{\span\preamble@@\crcr#2\crcr}}$}%
 \measuring@false\Let@\attag@{#1}\tag@C\tag@CC
 \tabskip\centering@\halign to\displaywidth
 {\span\preamble@@\crcr#2\crcr\black@{\wd\savealignat@}}}
\def\gather{\RIfMIfI@\DN@{\onlydmatherr@\gather}\else
 \ingather@true\inany@true\def\tag{&\relax}%
 \vspace@\allowdisplaybreak@\displaybreak@\intertext@
 \displ@y\Let@
 \iftagsleft@\DN@{\csname gather \endcsname}\else
  \DN@{\csname gather \space\endcsname}\fi\fi
 \else\DN@{\onlydmatherr@\gather}\fi\next@}
\def\exstring@{\expandafter\eat@\string}
\def\newcounter#1{\define#1{}%
 \edef\next@{\def\noexpand#1{\futurelet\noexpand\next
  \csname\exstring@#1@Z\endcsname}}\next@
 \edef\next@{\def\csname\exstring@#1@Z\endcsname
  {\global\advance\csname\exstring@#1@C\endcsname\@ne
  {\csname\exstring@#1@F\endcsname\csname\exstring@#1@S\endcsname
   {\csname\exstring@#1@P\endcsname\csname\exstring@#1@N\endcsname
   {\noexpand\number\csname\exstring@#1@C\endcsname}%
   \csname\exstring@#1@Q\endcsname}}%
  \noexpand\ifx\noexpand\next\noexpand\label
   \def\noexpand\next@\noexpand\label########1{{\noexpand\noexpands@
    \xdef\noexpand\Thelabel@{\csname\exstring@#1@N\endcsname
     {\noexpand\number\csname\exstring@#1@C\endcsname}}%
    \xdef\noexpand\Thelabel@@@{\noexpand\number
     \csname\exstring@#1@C\endcsname}%
    \xdef\noexpand\Thelabel@@{\csname\exstring@#1@S\endcsname
     {\csname\exstring@#1@P\endcsname
     \csname\exstring@#1@N\endcsname
     {\noexpand\number\csname\exstring@#1@C\endcsname}%
     \csname\exstring@#1@Q\endcsname}}%
    \xdef\noexpand\Thelabel@@@@{\csname\exstring@#1@P\endcsname
     \csname\exstring@#1@N\endcsname
     {\noexpand\number\csname\exstring@#1@C\endcsname}%
     \csname\exstring@#1@Q\endcsname}}%
    {\noexpand\locallabel@\noexpand\label{########1}}}%
   \noexpand\else\let\noexpand\next@\relax\noexpand\fi\noexpand\next@}}\next@
 \expandafter\newcount@\csname\exstring@#1@C\endcsname
 \expandafter\let\csname\exstring@#1@N\endcsname\arabic
 \expandafter\def\csname\exstring@#1@S\endcsname##1{##1\/}%
 \expandafter\let\csname\exstring@#1@P\endcsname\empty
 \expandafter\let\csname\exstring@#1@Q\endcsname\empty
 \expandafter\def\csname\exstring@#1@F\endcsname{\rm}%
 }
\def\HASH@#1#2{\ifnum#2=\z@\else
 \edef\next@{\toks@{\the\toks@\the\hashtoks@#2}%
 \toks@@{\the\toks@@{\the\hashtoks@#2}}}\next@\expandafter\HASH@\fi}
\def\HASH@@{\toks@{}\toks@@{}\expandafter\HASH@\macpar@00}
\def\usecounter#1#2{\expandafter\ifx\csname\exstring@#1@Z\endcsname
 \relax\Err@{\noexpand#1not created with \string\newcounter}\fi
 \expandafter\let\csname\exstring@#1@@Z\endcsname\relax
 \expandafter\let\csname\exstring@#1@@Z@\endcsname\relax
 \expandafter\let\csname\exstring@#1@@Z@@\endcsname\relax
 \edef\next@{\def\noexpand#2{\futurelet\noexpand\next
  \csname\exstring@#1@@Z\endcsname}}\next@
 \edef\next@{\def\csname\exstring@#1@@Z\endcsname{\noexpand\ifx
  \noexpand\next\noexpand\label\def\noexpand\next@\noexpand\label
   ########1{\csname\exstring@#1@@Z@\endcsname
   {\noexpand#1\noexpand\label{########1}}}%
   \noexpand\else\noexpand\ifx\noexpand\next
   \noexpand"\def\noexpand\next@\noexpand"########1\noexpand"%
   {\csname\exstring@#1@@Z@\endcsname{{\expandafter\noexpand
   \csname\exstring@#1@F\endcsname
   \let\noexpand\pre\expandafter\noexpand\csname\exstring@#1@P\endcsname
   \let\noexpand\post\expandafter\noexpand\csname\exstring@#1@Q\endcsname
   \let\noexpand\style\expandafter\noexpand\csname\exstring@#1@S\endcsname
   \let\noexpand\numstyle\expandafter\noexpand\csname\exstring@#1@N\endcsname
   ########1}}}\noexpand\else
   \def\noexpand\next@{\csname\exstring@#1@@Z@\endcsname{\noexpand#1}}%
   \noexpand\fi\noexpand\fi\noexpand\next@}}\next@
 \def\next@{\expandafter\expandafter\expandafter\unmacro@\expandafter
  \meaning\csname\exstring@#1@@Z@@\endcsname\unmacro@
  \HASH@@
  \edef\next@{\def\csname\exstring@#1@@Z@\endcsname\the\toks@{%
   \expandafter\noexpand\csname\exstring@#1@@Z@@\endcsname\the\toks@@
   \noexpand\FNSSP@}}\next@}%
 \afterassignment\next@
 \expandafter\def\csname\exstring@#1@@Z@@\endcsname}
\def\listbi@{\penalty50 \medskip}
\def\listbii@{\penalty100 \smallskip}
\let\listbiii@\relax
\let\listbiv@\relax
\let\listbv@\relax
\def\listmi@{\advance\leftskip30\p@\relax}
\let\listmii@\listmi@
\let\listmiii@\listmi@
\let\listmiv@\listmi@
\let\listmv@\listmi@
\def\itemi@#1{\noindent@@\llap{#1\hskip5\p@}}
\let\itemii@\itemi@
\let\itemiii@\itemi@
\let\itemiv@\itemi@
\let\itemv@\itemi@
\def\liste@{\penalty-50 \medskip}
\def\listei@{\penalty-100 \smallskip}
\let\listeii@\relax
\let\listeiii@\relax
\let\listeiv@\relax
\expandafter\newcount\csname list@C1\endcsname
\csname list@C1\endcsname\z@
\expandafter\newcount\csname list@C2\endcsname
\csname list@C2\endcsname\z@
\expandafter\newcount\csname list@C3\endcsname
\csname list@C3\endcsname\z@
\expandafter\newcount\csname list@C4\endcsname
\csname list@C4\endcsname\z@
\expandafter\newcount\csname list@C5\endcsname
\csname list@C5\endcsname\z@
\expandafter\let\csname list@P1\endcsname\empty
\expandafter\let\csname list@P2\endcsname\empty
\expandafter\let\csname list@P3\endcsname\empty
\expandafter\let\csname list@P4\endcsname\empty
\expandafter\let\csname list@P5\endcsname\empty
\expandafter\let\csname list@Q1\endcsname\empty
\expandafter\let\csname list@Q2\endcsname\empty
\expandafter\let\csname list@Q3\endcsname\empty
\expandafter\let\csname list@Q4\endcsname\empty
\expandafter\let\csname list@Q5\endcsname\empty
\expandafter\def\csname list@S1\endcsname#1{{\rm(}{#1\/}{\rm)}}
\expandafter\def\csname list@S2\endcsname#1{{\rm(}{#1\/}{\rm)}}
\expandafter\def\csname list@S3\endcsname#1{{\rm(}{#1\/}{\rm)}}
\expandafter\def\csname list@S4\endcsname#1{{\rm(}{#1\/}{\rm)}}
\expandafter\def\csname list@S5\endcsname#1{{\rm(}{#1\/}{\rm)}}
\expandafter\let\csname list@N1\endcsname\arabic
\expandafter\let\csname list@N2\endcsname\arabic
\expandafter\let\csname list@N3\endcsname\arabic
\expandafter\let\csname list@N4\endcsname\arabic
\expandafter\let\csname list@N5\endcsname\arabic
\expandafter\def\csname list@F1\endcsname{\rm}
\expandafter\def\csname list@F2\endcsname{\rm}
\expandafter\def\csname list@F3\endcsname{\rm}
\expandafter\def\csname list@F4\endcsname{\rm}
\expandafter\def\csname list@F5\endcsname{\rm}
\newcount\listlevel@
\listlevel@\z@
\def\list@@C{\csname list@C\number\listlevel@\endcsname}
\def\list@@P{\csname list@P\number\listlevel@\endcsname}
\def\list@@Q{\csname list@Q\number\listlevel@\endcsname}
\def\list@@S{\csname list@S\number\listlevel@\endcsname}
\def\list@@N{\csname list@N\number\listlevel@\endcsname}
\def\list@@F{\csname list@F\number\listlevel@\endcsname}
\newif\iffirstitemi@
\newif\iffirstitemii@
\newif\iffirstitemiii@
\newif\iffirstitemiv@
\newif\iffirstitemv@
\def\Firstitem@true{\csname firstitem\romannumeral\listlevel@
 @true\endcsname}
\def\Firstitem@false{\csname firstitem\romannumeral\listlevel@
 @false\endcsname}
\def\Listm@{\csname listm\romannumeral\listlevel@ @\endcsname}
\def\Item@{\csname item\romannumeral\listlevel@ @\endcsname}
\def\Liste@{\csname liste\romannumeral\listlevel@ @\endcsname}
\newif\iflistcontinue@
\def\keepitem{\listcontinue@true}
\newcount\list@C@
\def\list{%
 \iflistcontinue@\csname list@C1\endcsname\csname list@C@\endcsname\fi
 \global\csname list@C2\endcsname\z@
 \global\csname list@C3\endcsname\z@
 \global\csname list@C4\endcsname\z@
 \global\csname list@C5\endcsname\z@
 \begingroup
 \firstitemi@true
 \listlevel@\@ne
 \def\item{\FN@\item@}%
 \FN@\list@}
\Invalid@\runinitem
\def\list@{\ifx\next\par
 \DN@\par{\FN@\list@}\else
 \ifx\next\runinitem
  \DN@\runinitem{\FN@\runinitem@}\else
  \DN@{\par\dimen@\parskip\parskip\dimen@}\fi\fi\next@}
\newif\ifoutlevel@
\newif\ifrunin@
\def\item@{%
 \ifoutlevel@\Liste@\outlevel@false\fi
 \ifrunin@\runin@false\par
  \dimen@\parskip\parskip\dimen@
  \Listm@\fi
 \iffirstitemi@\listbi@\listmi@\firstitemi@false\else\par\fi
 \iffirstitemii@\listbii@\listmii@\firstitemii@false\else\par\fi
 \iffirstitemiii@\listbiii@\listmiii@\firstitemiii@false\else\par\fi
 \iffirstitemiv@\listbiv@\listmiv@\firstitemiv@false\else\par\fi
 \iffirstitemv@\listbv@\listmv@\firstitemv@false\else\par\fi
 \DN@"##1"{{\let\pre\list@@P\let\post\list@@Q
  \let\style\list@@S\let\numstyle\list@@N
  \vskip-\parskip
  \Item@{\list@@F##1}%
  \noexpands@
  \Qlabel@{##1}}%
  \locallabel@
  \FNSSP@}%
 \DNii@{\global\advance\list@@C\@ne
  {\noexpands@
   \xdef\Thelabel@@@{\number\list@@C}%
   \xdefThelabel@\list@@N
   \xdef\Thelabel@@@@{\list@@P\Thelabel@\list@@Q}%
   \xdefThelabel@@\list@@S
  }%
  \locallabel@
  \vskip-\parskip
  \Item@{\list@@F\thelabel@@}%
  \FN@\pretendspace@}%
 \ifx\next"\expandafter\next@\else\expandafter\nextii@\fi}
\def\runinitem@{%
  \runin@true
  \Firstitem@false
  \DN@"##1"{{\let\pre\list@@P\let\post\list@@Q
   \let\style\list@@S\let\numstyle\list@@N
   \unskip\space{\list@@F##1} %
   \noexpands@
   \Qlabel@{##1}}%
   \locallabel@
   \ignorespaces}%
  \DNii@{\global\advance\list@@C\@ne
   {\noexpands@
    \xdef\Thelabel@@@{\number\list@@C}%
    \xdefThelabel@\list@@N
    \xdef\Thelabel@@@@{\list@@P\Thelabel@\list@@Q}%
    \xdefThelabel@@\list@@S
   }%
   \locallabel@
   \unskip\space{\list@@F\thelabel@@} }%
  \ifx\next"\expandafter\next@\else\expandafter\nextii@\fi}
\def\inlevel{\ifnum\listlevel@=5
 \DN@{\Err@{Already 5 levels down}}\else
 \DN@{\begingroup\advance\listlevel@\@ne
 \Firstitem@true\FN@\inlevel@}\fi\next@}
\def\inlevel@{\ifx\next\par
 \DN@\par{\FN@\inlevel@}\else
 \ifx\next\runinitem
  \DN@\runinitem{\FN@\runinitem@}\else
  \let\next@\relax\fi\fi\next@}
\def\outlevel{\ifnum\listlevel@=\@ne
 \Err@{At top level}\else
 \par\global\list@@C\z@\endgroup\outlevel@true\fi}
\def\endlist{%
 \expandafter\global\csname list@C@\endcsname\csname list@C1\endcsname
 \par
 \global\toks\@ne{}\count@\listlevel@
 {\loop
  \ifnum\count@>\z@\global\toks\@ne\expandafter{\the\toks\@ne\endgroup}%
  \advance\count@\m@ne
  \repeat}%
 \the\toks\@ne
 \liste@
 \listcontinue@false\global\csname list@C1\endcsname\z@
 \vskip-\parskip
 \noindent@@
 \FN@\pretendspace@}
\newif\iffirstdescribe@
\def\describe{\par
 \begingroup\firstdescribe@true
 \def\item##1{%
  \iffirstdescribe@\penalty50 \medskip\vskip-\parskip
  \firstdescribe@false\else\par\fi
  \noindent@@\hangindent2pc\hangafter\@ne
  {\bf##1}\hskip.5em}}

\Invalid@\pullin
\Invalid@\pullinmore
\newif\iffirstpull@
\def\margins{\par\begingroup\firstpull@true
 \def\pullin##1##2{\par
  \iffirstpull@\firstpull@false\else\endgroup\fi
  \begingroup\DN@{##1}%
  \ifx\next@\empty\leftskip\z@\else\ifx\next@\space\leftskip\z@
  \else\leftskip##1\fi\fi
  \DN@{##2}\ifx\next@\empty\rightskip\z@\else\ifx\next@\space
  \rightskip\z@\else\rightskip##2\fi\fi\ignorespaces}%
 \def\pullinmore##1##2{\par
  \xdef\Next@{\leftskip\the\leftskip\relax\rightskip\the\rightskip\relax}%
  \iffirstpull@\firstpull@false\else\endgroup\fi
  \begingroup\Next@
  \DN@{##1}%
  \ifx\next@\empty\else\ifx\next@\space\else\advance\leftskip##1\fi\fi
  \DN@{##2}\ifx\next@\empty\else\ifx\next@\space\else
  \advance\rightskip##2\fi\fi\ignorespaces}}

\newif\ifnopunct@
\newif\ifnospace@
\newif\ifoverlong@
\let\nofrillslist@\empty
\let\overlonglist@\empty
\def\nopunct{\nopunct@true\FN@\nopunct@}
\def\nospace{\nospace@true\FN@\nospace@}
\def\overlong{\overlong@true\FN@\overlong@}
\def\nopunct@{\ifx\next\nospace
 \DN@\nospace{\nospace@true\FN@\nopnos@}\else\ifx\next\overlong
 \DN@\overlong{\overlong@true\FN@\nopol@}\else
 \let\next@\nopunct@@\fi\fi\next@}
\def\nopunct@@#1{\ismember@\nofrillslist@#1%
 \iftest@\let\next@#1\else
 \DN@{\nopunct@false\Err@{\noexpand\nopunct can't be used with
 \string#1}#1}\fi\next@}
\def\nospace@{\ifx\next\nopunct
 \DN@\nopunct{\nopunct@true\FN@\nopnos@}\else\ifx\next\overlong
 \DN@\overlong{\overlong@true\FN@\nosol@}\else
 \let\next@\nospace@@\fi\fi\next@}
\def\nospace@@#1{\ismember@\nofrillslist@#1%
 \iftest@\let\next@#1\else
 \DN@{\nospace@false\Err@{\noexpand\nospace can't be used with
 \string#1}#1}\fi\next@}
\def\overlong@{\ifx\next\nopunct
 \DN@\nopunct{\nopunct@true\FN@\nopol@}\else\ifx\next\nospace
 \DN@\nospace{\nospace@true\FN@\nosol@}\else
 \let\next@\overlong@@\fi\fi\next@}
\def\overlong@@#1{\ismember@\overlonglist@#1%
 \iftest@\let\next@#1\else
 \DN@{\overlong@false\Err@{\noexpand\overlong can't be used with
 \string#1}#1}\fi\next@}
\def\nopnos@{\ifx\next\overlong
 \DN@\overlong{\overlong@true\nopnosol@}\else
 \let\next@\nopnos@@\fi\next@}
\def\nopol@{\ifx\next\nospace
 \DN@\nospace{\nospace@true\nopnosol@}\else
 \let\next@\nopol@@\fi\next@}
\def\nosol@{\ifx\next\nopunct
 \DN@\nopunct{\nopunct@true\nopnosol@}\else
 \let\next@\nosol@@\fi\next@}
\def\nopnos@@#1{\ismember@\nofrillslist@#1%
 \iftest@\let\next@#1\else
 \DN@{\nopunct@false\nospace@false
  \Err@{\noexpand\nopunct\noexpand\nospace
   can't be used with \string#1}#1}\fi\next@}
\def\testii@#1{\ismember@\nofrillslist@#1%
 \iftest@\let\nextiii@ T\else\let\nextiii@ F\fi
 \ismember@\overlonglist@#1%
 \iftest@\let\nextiv@ T\else\let\nextiv@ F\fi
 \test@false\if\nextiii@ T\if\nextiv@ T\test@true\fi\fi}
\def\nopol@@#1{\testii@{#1}%
 \iftest@\let\next@#1%
 \else\DN@{\if\nextiii@ T\else\nopunct@false\fi
  \if\nextiv@ T\else\overlong@false\fi
  \Err@{\if\nextiii@ T\else\noexpand\nopunct\fi
  \if\nextiv@ T\else\noexpand\overlong\fi can't be used
  with \string#1}#1}\fi\next@}
\def\nosol@@#1{\testii@{#1}%
 \iftest@\let\next@#1%
 \else\DN@{\if\nextiii@ T\else\nospace@false\fi
  \if\nextiv@ T\else\overlong@false\fi
  \Err@{\if\nextiii@ T\else\noexpand\nospace\fi
  \if\nextiv@ T\else\noexpand\overlong\fi can't be used
  with \string#1}#1}\fi\next@}
\def\nopnosol@#1{\testii@{#1}%
 \iftest@\let\next@#1%
 \else\DN@{\if\nextiii@ T\else\nopunct@false\nospace@false\fi
  \if\nextiv@ T\else\overlong@false\fi
  \Err@{\if\nextiii@ T\else\noexpand\nopunct\noexpand\nospace\fi
  \if\nextiv@ T\else\noexpand\overlong\fi can't be used
  with \string#1}#1}\fi\next@}
\def\punct@#1{\ifnopunct@\else#1\fi}
\def\addspace@#1{\ifnospace@\else#1\fi}
\def\hss@{\ifoverlong@\z@ plus\@m\p@ minus\@m\p@
 \else \z@ plus\@m\p@\fi}
\rightadd@\demo\to\nofrillslist@
\newif\ifclaim@
\def\exxx@{\expandafter\expandafter\expandafter\eat@\expandafter\string}
\let\colon@:
\def\demo#1{\ifclaim@
 \Err@{Previous \expandafter\noexpand\claimtype@ has
  no matching \string\end\exxx@\claimtype@}%
 \let\next@\relax
 \else
  \par
  \ifdim\lastskip<\smallskipamount\removelastskip\smallskip\fi
  \begingroup
  \noindent@@{\smc\ignorespaces#1\unskip
   \punct@{\null\colon@}\addspace@\enspace}%
  \nopunct@false\nospace@false
  \rm
  \DN@{\FNSSP@}%
 \fi
 \next@}
\def\enddemo{\par\endgroup\nopunct@false\nospace@false\smallskip}
\rightadd@\claim\to\nofrillslist@
\def\claim@F{\smc}
\def\claim@@@F{\csname\exxx@\claimtype@ @F\endcsname}
\def\claimformat@#1#2#3{%
 \medbreak\noindent@@{\smc#1 {\claim@@@F#2} #3%
 \punct@{\null.}\addspace@\enspace}\sl}
\def\claimformat@@#1#2{\claimformat@{\ignorespaces#1\unskip}%
 {\ifx\thelabel@@\empty\unskip\else\thelabel@@\fi}%
 {\ignorespaces#2\unskip}%
 \let\Claimformat@@\claimformat@@\FNSSP@}
\let\Claimformat@@\claimformat@@
\def\claim@@@P{\csname\exxx@\claimtype@ @P\endcsname}
\def\claim@@@Q{\csname\exxx@\claimtype@ @Q\endcsname}
\def\claim@@@S{\csname\exxx@\claimtype@ @S\endcsname}
\def\claim@@@N{\csname\exxx@\claimtype@ @N\endcsname}
\def\claim@@@C{\csname claim@C\claimclass@\endcsname}
\newcount\claim@C
\claim@C\z@
\let\claim@P\empty
\let\claim@Q\empty
\def\claim@S#1{#1\/}
\let\claim@N\arabic
\def\claim{\claim@true\let\claimclass@\empty
 \def\claimtype@{\claim}\FN@\claim@}
\def\claim@{%
 \ifx\next\c
  \let\next@\claim@c
 \else
  \ifx\next"%
   \let\next@\claim@q
  \else
   \begingroup\global\advance\claim@C\@ne
   {\noexpands@
    \xdef\Thelabel@@@{\number\claim@C}%
    \xdefThelabel@\claim@N
    \xdef\Thelabel@@@@{\claim@P\Thelabel@\claim@Q}%
    \xdefThelabel@@\claim@S
   }%
   \locallabel@
   \let\next@\Claimformat@@
  \fi
 \fi
 \next@}
\def\claim@c\c#1{\claim@true\begingroup
 \expandafter
 \ifx\csname claim@C#1\endcsname\relax
  \expandafter\newcount@\csname claim@C#1\endcsname
  \global\csname claim@C#1\endcsname\@ne
 \else
  \global\advance\csname claim@C#1\endcsname\@ne
 \fi
 \def\claimclass@{#1}%
 {\noexpands@
  \xdef\Thelabel@@@{\number\claim@@@C}%
  \xdefThelabel@\claim@@@N
  \xdef\Thelabel@@@@{\claim@@@P\Thelabel@\claim@@@Q}%
  \xdefThelabel@@\claim@@@S
 }%
 \locallabel@
 \FNSS@\claim@c@}
\def\claim@q"#1"{\begingroup
 {\let\pre\claim@@@P\let\post\claim@@@Q
  \let\style\claim@@@S\let\numstyle\claim@@@N
  \noexpands@
  \Qlabel@{#1}}%
 \locallabel@
 \FNSS@\claim@q@}
\def\claim@c@{\ifx\next"%
 \global\advance\claim@@@C\m@ne\let\next@\claim@cq
 \else\let\next@\Claimformat@@\fi\next@}
\def\claim@cq"#1"{{\let\pre\claim@@@P\let\post\claim@@@Q
 \let\style\claim@@@S\let\numstyle\claim@@@N
 \noexpands@
 \Qlabel@{#1}}%
 \locallabel@
 \FNSS@\Claimformat@@}
\def\claim@q@{\ifx\next\c\expandafter\claim@qc
 \else\expandafter\Claimformat@@\fi}
\def\claim@qc\c#1{\expandafter\ifx\csname claim@C#1\endcsname\relax
 \expandafter\newcount@\csname claim@C#1\endcsname
 \global\csname claim@C#1\endcsname\z@\fi
 \FNSS@\Claimformat@@}
\def\endclaim{\endgroup\claim@false\nopunct@false\nospace@false
 \let\Claimformat@@\claimformat@@\medbreak}
\Invalid@\claimclause
\def\newclaim{\FN@\newclaim@}
\def\newclaim@{\ifx\next\claimclause
 \DN@\claimclause##1{\newclaim@@{##1}}\else
 \DN@{\newclaim@@\relax}\fi\next@}
\def\claimlist@{\\\claim}
\newtoks\claim@i
\newtoks\claim@v
\let\noclaimclause@=F
\def\newclaim@@#1#2#3\c#4#5{\define#2{}%
 \rightadd@#2\to\claimlist@\rightadd@#2\to\nofrillslist@%
 \expandafter\def\csname\exstring@#2@P\endcsname{\claim@P}%
 \expandafter\def\csname\exstring@#2@Q\endcsname{\claim@Q}%
 \expandafter\def\csname\exstring@#2@S\endcsname{\claim@S}%
 \expandafter\def\csname\exstring@#2@N\endcsname{\claim@N}%
 \expandafter\def\csname\exstring@#2@F\endcsname{\claim@F}%
 \expandafter\def\csname end\exstring@#2\endcsname{\endclaim}%
 \expandafter\ifx\csname claim@C#4\endcsname\relax
  \expandafter\newcount@\csname claim@C#4\endcsname
  \global\csname claim@C#4\endcsname\z@\fi
 \edef\next@{\let\csname\exstring@#2@C\endcsname
   \csname claim@C#4\endcsname}\next@
 \def#2{\ifx\noclaimclause@ T\else#1\fi
  \global\claim@i{#1}\gdef\claim@iv{#4}\global\claim@v{#5}%
  \def\claimtype@{#2}\def\Claimformat@@{\claimformat@@{#5}}\claim@c\c{#4}}}
\def\shortenclaim#1#2{\define#2{}%
 \ismember@\claimlist@#1%
 \iftest@
  \rightadd@#2\to\nofrillslist@%
  \expandafter\def\csname\exstring@#2@P\endcsname
   {\csname\exstring@#1@P\endcsname}%
  \expandafter\def\csname\exstring@#2@Q\endcsname
   {\csname\exstring@#1@Q\endcsname}%
  \expandafter\def\csname\exstring@#2@S\endcsname
   {\csname\exstring@#1@S\endcsname}%
  \expandafter\def\csname\exstring@#2@N\endcsname
   {\csname\exstring@#1@N\endcsname}%
  \expandafter\def\csname\exstring@#2@F\endcsname
   {\csname\exstring@#1@F\endcsname}%
  \expandafter\def\csname end\exstring@#2\endcsname{\endclaim}%
  \edef\next@{\let\csname\exstring@#2@C\endcsname
    \csname claim\exstring@#1C\endcsname}\next@
  \setbox\z@\vbox{\let\noclaimclause@ T#1""\relax\endgroup}%
  \edef#2{\the\claim@i
   \def\noexpand\claimtype@{\noexpand#2}%
   \def\noexpand\Claimformat@@{\noexpand\claimformat@@{\the\claim@v}\relax}%
   \noexpand\claim@c\noexpand\c{\claim@iv}}%
 \else
  \Err@{\noexpand#1not yet created by \string\newclaim}%
 \fi}
\def\classtest@#1{\DN@{#1}\ifx\next@\claimclass@
 \test@true\else\test@false\fi}
\def\typetest@#1{\DN@{#1}\ifx\next@\claimtype@\test@true\else
  \test@false\fi}
\newif\iftoc@
\def\tocfile{\iftoc@\else\alloc@@7\write\chardef\sixt@@n\toc@
 \immediate\openout\toc@=\jobname.toc
 \alloc@@7\write\chardef\sixt@@n\tic@
 \immediate\openout\tic@=\jobname.tic
 \global\toc@true\fi}
\rightadd@\hl\to\nofrillslist@
\rightadd@\HL\to\overlonglist@
\def\HL@@C{\csname HL@C\HLlevel@\endcsname}
\def\HL@@P{\csname HL@P\HLlevel@\endcsname}
\def\HL@@Q{\csname HL@Q\HLlevel@\endcsname}
\def\HL@@S{\csname HL@S\HLlevel@\endcsname}
\def\HL@@N{\csname HL@N\HLlevel@\endcsname}
\def\HL@@F{\csname HL@F\HLlevel@\endcsname}
\def\HL@@@C{\csname\exxx@\HLtype@ @C\endcsname}
\def\HL@@@P{\csname\exxx@\HLtype@ @P\endcsname}
\def\HL@@@Q{\csname\exxx@\HLtype@ @Q\endcsname}
\def\HL@@@S{\csname\exxx@\HLtype@ @S\endcsname}
\def\HL@@@N{\csname\exxx@\HLtype@ @N\endcsname}
\def\HL#1{\expandafter
 \ifx\csname HL@C#1\endcsname\relax
  \DN@{\Err@{\string\HL#1 not defined in this style}}%
 \else
  \DN@{\gdef\HLlevel@{#1}\def\HLname@{\HL{#1}}\let\HLtype@\relax\FNSS@\HL@}%
 \fi
 \next@}%
\newif\ifquoted@
\let\aftertoc@\relax
\def\HL@{%
 \DN@"##1"##2\endHL{\def\entry@{##2}\quoted@true
  {\noexpands@
  \ifx\HLtype@\relax
   \let\pre\HL@@P\let\post\HL@@Q\let\style\HL@@S\let\numstyle\HL@@N
  \else
   \let\pre\HL@@@P\let\post\HL@@@Q\let\style\HL@@@S\let\numstyle\HL@@@N
  \fi
  \Qlabel@{##1}\let\style\relax\xdef\Qlabel@@@@{##1}%
  \xdef\Thepref@{\Thelabel@@@@}}%
  \csname HL@\HLlevel@\endcsname##2\endHL
  \let\pref\Thepref@
  \csname HL@I\HLlevel@\endcsname
  \csname HL@J\HLlevel@\endcsname
  \let\pref\pref@
  \HLtoc@
  \aftertoc@
  \let\aftertoc@\relax\overlong@false}%
 \DNii@##1\endHL{\def\entry@{##1}\quoted@false
  {\noexpands@
  \ifx\HLtype@\relax
   \global\advance\HL@@C\@ne
   \xdef\Thelabel@@@{\number\HL@@C}%
   \xdefThelabel@{\HL@@N}%
   \xdef\Thelabel@@@@{\HL@@P\Thelabel@\HL@@Q}%
   \xdefThelabel@@{\HL@@S}%
  \else
   \global\advance\HL@@@C\@ne
   \xdef\Thelabel@@@{\number\HL@@@C}%
   \xdefThelabel@{\HL@@@N}%
   \xdef\Thelabel@@@@{\HL@@@P\Thelabel@\HL@@@Q}%
   \xdefThelabel@@{\HL@@@S}%
  \fi
  \xdef\Thepref@{\Thelabel@@@@}}%
  \csname HL@\HLlevel@\endcsname##1\endHL
  \let\pref\Thepref@
  \csname HL@I\HLlevel@\endcsname
  \csname HL@J\HLlevel@\endcsname
  \let\pref\pref@
  \HLtoc@
  \aftertoc@
  \let\aftertoc@\relax\overlong@false}%
 \ifx\next"\expandafter\next@\else\expandafter\nextii@\fi}%
\Invalid@\endHL
\def\hl@@C{\csname hl@C\hllevel@\endcsname}
\def\hl@@P{\csname hl@P\hllevel@\endcsname}
\def\hl@@Q{\csname hl@Q\hllevel@\endcsname}
\def\hl@@S{\csname hl@S\hllevel@\endcsname}
\def\hl@@N{\csname hl@N\hllevel@\endcsname}
\def\hl@@F{\csname hl@F\hllevel@\endcsname}
\def\hl@@@C{\csname\exxx@\hltype@ @C\endcsname}
\def\hl@@@P{\csname\exxx@\hltype@ @P\endcsname}
\def\hl@@@Q{\csname\exxx@\hltype@ @Q\endcsname}
\def\hl@@@S{\csname\exxx@\hltype@ @S\endcsname}
\def\hl@@@N{\csname\exxx@\hltype@ @N\endcsname}
\def\hl#1{\expandafter
 \ifx\csname hl@C#1\endcsname\relax
  \DN@{\Err@{\string\hl#1 not defined in this style}}%
 \else
  \DN@{\gdef\hllevel@{#1}\def\hlname@{\hl{#1}}\let\hltype@\relax\FNSS@\hl@}%
 \fi
 \next@}
\def\hl@{%
 \DN@"##1"##2{\def\entry@{##2}\quoted@true
  {\noexpands@
  \ifx\hltype@\relax
   \let\pre\hl@@P\let\post\hl@@Q\let\style\hl@@S\let\numstyle\hl@@N
  \else
   \let\pre\hl@@@P\let\post\hl@@@Q\let\style\hl@@@S\let\numstyle\hl@@@N
  \fi
  \Qlabel@{##1}\let\style\relax\xdef\Qlabel@@@@{##1}%
  \xdef\Thepref@{\Thelabel@@@@}}%
  \csname hl@\hllevel@\endcsname{##2}%
  \let\pref\Thepref@
  \csname hl@I\hllevel@\endcsname
  \csname hl@J\hllevel@\endcsname
  \let\pref\pref@
  \hltoc@
  \aftertoc@
  \let\aftertoc@\relax\nopunct@false\nospace@false\FNSSP@}%
 \DNii@##1{\def\entry@{##1}\quoted@false
  {\noexpands@
  \ifx\hltype@\relax
   \global\advance\hl@@C\@ne
   \xdef\Thelabel@@@{\number\hl@@C}%
   \xdefThelabel@{\hl@@N}%
   \xdef\Thelabel@@@@{\hl@@P\Thelabel@\hl@@Q}%
   \xdefThelabel@@{\hl@@S}%
  \else
   \global\advance\hl@@@C\@ne
   \xdef\Thelabel@@@{\number\hl@@@C}%
   \xdefThelabel@{\hl@@@N}%
   \xdef\Thelabel@@@@{\hl@@@P\Thelabel@\hl@@@Q}%
   \xdefThelabel@@{\hl@@@S}%
  \fi
  \xdef\Thepref@{\Thelabel@@@@}}%
  \csname hl@\hllevel@\endcsname{##1}%
  \let\pref\Thepref@
  \csname hl@I\hllevel@\endcsname
  \csname hl@J\hllevel@\endcsname
  \let\pref\pref@
  \hltoc@
  \aftertoc@
  \let\aftertoc@\relax\nopunct@false\nospace@false\FNSSP@}%
 \ifx\next"\expandafter\next@\else\expandafter\nextii@\fi}%
\def\six@#1#2 #3 #4 #5 #6 #7 {\DN@{#2}\ifx\next@\empty
 \DN@##1\six@{}\else
 \write#1{ #2 #3 #4 #5 #6 #7}\DN@{\six@#1}\fi
 \next@}
\def\Sixtoc@{\ifx\macdef@\empty\else
 \DN@##1##2\next@{\def\macdef@{##1##2}}%
 \expandafter\next@\macdef@\next@
 \edef\next@
  {\noexpand\six@\toc@\macdef@
  \space\space\space\space\space\space\space\space\space\space\space\space
  \noexpand\six@}%
 \next@\let\macdef@\relax\fi}
\def\QorThelabel@@@@{\ifquoted@
 \noexpand\noexpand\noexpand"\Qlabel@@@@\noexpand\noexpand\noexpand"\else
 \Thelabel@@@@\fi}
\def\HLtoc@{%
 \iftoc@
 \expandafter\expandafter\expandafter\unmacro@
  \expandafter\meaning\csname HL@W\HLlevel@\endcsname\unmacro@
  {\noexpands@\let\style\relax
   \edef\next@{\write\toc@{\noexpand\noexpand\expandafter\noexpand\HLname@
   {\macdef@}{\QorThelabel@@@@}}}%
  \next@}%
  \expandafter\unmacro@\meaning\entry@\unmacro@
  \Sixtoc@
  \write\toc@{\noexpand\Page{\number\pageno}{\page@N}%
   {\page@P}{\page@Q}^^J}%
 \fi}
\def\hltoc@{%
 \iftoc@
 \expandafter\expandafter\expandafter\unmacro@
  \expandafter\meaning\csname hl@W\hllevel@\endcsname\unmacro@
  {\noexpands@\let\style\relax
  \edef\next@{\write\toc@{%
   \ifnopunct@\noexpand\noexpand\noexpand\nopunct\fi
   \ifnospace@\noexpand\noexpand\noexpand\nospace\fi
   \noexpand\noexpand\expandafter\noexpand\hlname@
   {\macdef@}{\QorThelabel@@@@}}}%
  \next@}%
  \expandafter\unmacro@\meaning\entry@\unmacro@
  \Sixtoc@
  \write\toc@{\noexpand\Page{\number\pageno}{\page@N}%
   {\page@P}{\page@Q}^^J}%
 \fi}
\def\mainfile#1{\def\mainfile@{#1}}
\def\checkmainfile@{\ifx\mainfile@\undefined
 \Err@{No \noexpand\mainfile specified}\fi}
\expandafter\newcount@\csname HL@C1\endcsname
\csname HL@C1\endcsname\z@
\expandafter\def\csname HL@S1\endcsname#1{#1\null.}
\expandafter\let\csname HL@N1\endcsname\arabic
\expandafter\let\csname HL@P1\endcsname\empty
\expandafter\let\csname HL@Q1\endcsname\empty
\expandafter\def\csname HL@F1\endcsname{\bf}
\expandafter\let\csname HL@W1\endcsname\empty
\expandafter\newcount@\csname hl@C1\endcsname
\csname hl@C1\endcsname\z@
\expandafter\def\csname hl@S1\endcsname#1{#1\/}
\expandafter\let\csname hl@N1\endcsname\arabic
\expandafter\let\csname hl@P1\endcsname\empty
\expandafter\let\csname hl@Q1\endcsname\empty
\expandafter\def\csname hl@F1\endcsname{\bf}
\expandafter\let\csname hl@W1\endcsname\empty
\expandafter\def\csname HL@1\endcsname#1\endHL{\bigbreak
 {\locallabel@
  \global\setbox\@ne\vbox{\Let@\tabskip\hss@
  \halign to\hsize{\bf\hfil\ignorespaces##\unskip\hfil\cr
  \expandafter\ifx\csname HL@W1\endcsname\empty\else
   \csname HL@W1\endcsname\space\fi
  {\HL@@F\ifx\thelabel@@\empty\else\thelabel@@\space\fi}%
  \ignorespaces#1\crcr}}%
  }%
 \unvbox\@ne\nobreak\medskip}
\expandafter\def\csname hl@1\endcsname#1{\medbreak\noindent@@
 {\locallabel@
 \bf{\hl@@F\ifx\thelabel@@\empty\else\thelabel@@\space\fi}%
 \ignorespaces#1\unskip\punct@{\null.}\addspace@\enspace}}
\expandafter\def\csname HL@I1\endcsname{\Reset\hl1{1}%
 \ifx\pref\empty\newpre\hl1{}\else\newpre\hl1{\pref.}\fi}
\def\NameHL#1#2{\define#2{}%
 \expandafter\ifx\csname HL@R#1\endcsname\relax
 \else
  \def\nextiv@{\let\nextiii@}%
  \expandafter\nextiv@\csname HL@R#1\endcsname
  \expandafter\let\nextiii@\undefined
  \expandafter\let\csname\exxx@\nextiii@ @C\endcsname\relax
  \expandafter\let\csname\exxx@\nextiii@ @P\endcsname\relax
  \expandafter\let\csname\exxx@\nextiii@ @Q\endcsname\relax
  \expandafter\let\csname\exxx@\nextiii@ @S\endcsname\relax
  \expandafter\let\csname\exxx@\nextiii@ @N\endcsname\relax
  \expandafter\let\csname\exxx@\nextiii@ @F\endcsname\relax
  \expandafter\let\csname\exxx@\nextiii@ @W\endcsname\relax
  \expandafter\let\csname end\exxx@\nextiii@\endcsname\undefined
 \fi
 \expandafter\gdef\csname HL@R#1\endcsname{#2}%
 \expandafter\gdef\csname\exstring@#2@R\endcsname{{HL}{#1}}%
 \iftoc@\write\toc@{\noexpand\NameHL#1\noexpand#2^^J}\fi
 \rightadd@#2\to\overlonglist@
 \edef\next@{\let\csname\exstring@#2@C\endcsname\expandafter\noexpand
  \csname HL@C#1\endcsname}\next@
 \edef\next@{\let\csname\exstring@#2@P\endcsname\expandafter\noexpand
  \csname HL@P#1\endcsname}\next@
 \edef\next@{\let\csname\exstring@#2@Q\endcsname\expandafter\noexpand
  \csname HL@Q#1\endcsname}\next@
 \edef\next@{\let\csname\exstring@#2@S\endcsname\expandafter\noexpand
  \csname HL@S#1\endcsname}\next@
 \edef\next@{\let\csname\exstring@#2@N\endcsname\expandafter\noexpand
  \csname HL@N#1\endcsname}\next@
 \edef\next@{\let\csname\exstring@#2@F\endcsname\expandafter\noexpand
  \csname HL@F#1\endcsname}\next@
 \edef\next@{\let\csname\exstring@#2@W\endcsname\expandafter\noexpand
  \csname HL@W#1\endcsname}\next@
 \edef\next@{\def\noexpand#2####1\expandafter\noexpand
  \csname end\exstring@#2\endcsname
  {\def\noexpand\HLtype@{\noexpand#2}%
   \def\noexpand\HLname@{\noexpand#2}%
   \gdef\noexpand\HLlevel@{#1}%
   \noexpand\FNSS@\noexpand\HL@####1\noexpand\endHL}}%
  \next@
 \edef\next@{\noexpand\Invalid@\expandafter\noexpand
  \csname end\exstring@#2\endcsname}%
 \next@}
\def\Namehl#1#2{\define#2{}%
 \expandafter\ifx\csname hl@R#1\endcsname\relax
 \else
  \def\nextiv@{\let\nextiii@}%
  \expandafter\nextiv@\csname hl@R#1\endcsname
  \expandafter\let\nextiii@\undefined
  \expandafter\let\csname\exxx@\nextiii@ @C\endcsname\relax
  \expandafter\let\csname\exxx@\nextiii@ @P\endcsname\relax
  \expandafter\let\csname\exxx@\nextiii@ @Q\endcsname\relax
  \expandafter\let\csname\exxx@\nextiii@ @S\endcsname\relax
  \expandafter\let\csname\exxx@\nextiii@ @N\endcsname\relax
  \expandafter\let\csname\exxx@\nextiii@ @F\endcsname\relax
  \expandafter\let\csname\exxx@\nextiii@ @W\endcsname\relax
 \fi
 \expandafter\gdef\csname hl@R#1\endcsname{#2}%
 \expandafter\gdef\csname\exstring@#2@R\endcsname{{hl}{#1}}%
 \iftoc@\write\toc@{\noexpand\Namehl#1\noexpand#2^^J}\fi
 \rightadd@#2\to\nofrillslist@%
 \edef\next@{\let\csname\exstring@#2@C\endcsname\expandafter\noexpand
  \csname hl@C#1\endcsname}\next@
 \edef\next@{\let\csname\exstring@#2@P\endcsname\expandafter\noexpand
  \csname hl@P#1\endcsname}\next@
 \edef\next@{\let\csname\exstring@#2@Q\endcsname\expandafter\noexpand
  \csname hl@Q#1\endcsname}\next@
 \edef\next@{\let\csname\exstring@#2@S\endcsname\expandafter\noexpand
  \csname hl@S#1\endcsname}\next@
 \edef\next@{\let\csname\exstring@#2@N\endcsname\expandafter\noexpand
  \csname hl@N#1\endcsname}\next@
 \edef\next@{\let\csname\exstring@#2@F\endcsname\expandafter\noexpand
  \csname hl@F#1\endcsname}\next@
 \edef\next@{\let\csname\exstring@#2@W\endcsname\expandafter\noexpand
  \csname hl@W#1\endcsname}\next@
 \edef\next@{\def\noexpand#2{%
  \def\noexpand\hltype@{\noexpand#2}%
  \def\noexpand\hlname@{\noexpand#2}%
  \gdef\noexpand\hllevel@{#1}%
  \noexpand\FNSS@\noexpand\hl@}}%
 \next@}%
\def\Initialize{\FN@\Init@}
\def\Init@{\ifx\next\HL\let\next@\InitH@\else\ifx\next\hl\let\next@\InitH@
  \else\let\next@\InitS@\fi\fi\next@}
\def\InitH@#1#2{\expandafter\ifx\csname\exstring@#1@C#2\endcsname\relax
 \DN@{\Err@{\noexpand#1level #2 not defined in this style}}\else
 \DN@{\expandafter\gdef\csname\exstring@#1@J#2\endcsname}\fi\next@}
\def\InitC@#1#2{\edef\nextii@{\expandafter\noexpand\csname#1\endcsname{#2}}}
\def\InitS@#1{\expandafter\ifx\csname\exstring@#1@R\endcsname\relax
 \Err@{\noexpand#1not defined in this style}\let\next@\relax\else
 \DN@{\let\next@}\expandafter\next@\csname\exstring@#1@R\endcsname
 \expandafter\InitC@\next@
 \DN@{\expandafter\InitH@\nextii@}\fi\next@}
\def\value#1{\expandafter
 \ifx\csname\exstring@#1@C\endcsname\relax
  \expandafter\ifx\csname\exstring@#1@C1\endcsname\relax
   \DN@{\Err@{\noexpand\value can't be used with \string#1}}%
  \else
   \DN@{\value@#1}%
  \fi
 \else
  \DN@{\number\csname\exstring@#1@C\endcsname\relax}%
 \fi
 \next@}
\def\value@#1#2{\expandafter
 \ifx\csname\exstring@#1@C#2\endcsname\relax
  \DN@{\Err@{\string\value\string#1 can't be followed by \string#2}}%
 \else
  \DN@{\number\csname\exstring@#1@C#2\endcsname\relax}%
 \fi
 \next@}
\newcount\Value
\def\Evaluate#1{\expandafter
 \ifx\csname\exstring@#1@C\endcsname\relax
  \expandafter\ifx\csname\exstring@#1@C1\endcsname\relax
   \DN@{\Err@{\noexpand\Evaluate can't be used with \string#1}}%
  \else
   \DN@{\Evaluate@#1}%
  \fi
 \else
  \DN@{\global\Value\csname\exstring@#1@C\endcsname}%
 \fi
 \next@}
\def\Evaluate@#1#2{\expandafter
 \ifx\csname\exstring@#1@C#2\endcsname\relax
  \DN@{\Err@{\string\Evaluate\string#1 can't be followed by \string#2}}%
 \else
  \DN@{\global\Value\csname\exstring@#1@C#2\endcsname}%
 \fi\next@}
\def\pre#1{\expandafter
 \ifx\csname\exstring@#1@P\endcsname\relax
  \expandafter\ifx\csname\exstring@#1@P1\endcsname\relax
   \DN@{\Err@{\noexpand\pre can't be used with \string#1}}%
  \else
   \DN@{\pre@#1}%
  \fi
 \else
  \DN@{{\csname\exstring@#1@P\endcsname}}%
 \fi
 \next@}
\def\pre@#1#2{\expandafter
 \ifx\csname\exstring@#1@P#2\endcsname\relax
  \DN@{\Err@{\string\pre\string#1 can't be followed by \string#2}}%
 \else
  \DN@{{\csname\exstring@#1@P#2\endcsname}}%
 \fi
 \next@}
\def\post#1{\expandafter
 \ifx\csname\exstring@#1@Q\endcsname\relax
  \expandafter\ifx\csname\exstring@#1@Q1\endcsname\relax
   \DN@{\Err@{\noexpand\post can't be used with \string#1}}%
  \else
   \DN@{\post@#1}%
  \fi
 \else
  \DN@{{\csname\exstring@#1@Q\endcsname}}%
 \fi
 \next@}
\def\post@#1#2{\expandafter
 \ifx\csname\exstring@#1@Q#2\endcsname\relax
  \DN@{\Err@{\string\post\string#1 can't be followed by \string#2}}%
 \else
  \DN@{{\csname\exstring@#1@Q#2\endcsname}}%
 \fi
 \next@}
\def\style#1{\expandafter
 \ifx\csname\exstring@#1@S\endcsname\relax
  \expandafter\ifx\csname\exstring@#1@S1\endcsname\relax
   \DN@{\Err@{\noexpand\style can't be used with \string#1}}%
  \else
   \DN@{\style@#1}%
  \fi
 \else
  \DN@{\csname\exstring@#1@S\endcsname}%
 \fi
 \next@}
\def\style@#1#2{\expandafter
 \ifx\csname\exstring@#1@S#2\endcsname\relax
  \DN@{\Err@{\string\style\string#1 can't be followed by \string#2}}%
 \else
  \DN@{\csname\exstring@#1@S#2\endcsname}%
 \fi
 \next@}
\def\fontstyle#1{\expandafter
 \ifx\csname\exstring@#1@F\endcsname\relax
  \expandafter\ifx\csname\exstring@#1@F1\endcsname\relax
   \DN@{\Err@{\noexpand\fontstyle can't be used with \string#1}}%
  \else
   \DN@{\fontstyle@#1}%
  \fi
 \else
  \DN@##1{{\csname\exstring@#1@F\endcsname##1}}%
 \fi
 \next@}
\def\fontstyle@#1#2{\expandafter
 \ifx\csname\exstring@#1@F#2\endcsname\relax
  \DN@{\Err@{\string\fontstyle\string#1 can't be followed by \string#2}}%
 \else
  \DN@##1{{\csname\exstring@#1@F#2\endcsname##1}}%
 \fi
 \next@}
\def\Reset#1{\expandafter
 \ifx\csname\exstring@#1@C\endcsname\relax
  \expandafter\ifx\csname\exstring@#1@C1\endcsname\relax
   \DN@{\Err@{\noexpand\Reset can't be used with \string#1}}%
  \else
   \DN@{\Reset@#1}%
  \fi
 \else
  \DN@##1{\count@##1\relax\ifx#1\page\else\advance\count@\m@ne\fi
   \global\csname\exstring@#1@C\endcsname\count@}%
 \fi
 \next@}
\def\Reset@#1#2{\expandafter
 \ifx\csname\exstring@#1@C#2\endcsname\relax
  \DN@{\Err@{\string\Reset\string#1 can't be followed by \string#2}}%
 \else
  \DN@##1{\count@##1\relax\advance\count@\m@ne
   \global\csname\exstring@#1@C#2\endcsname\count@}%
 \fi
 \next@}
\def\Offset#1{\expandafter
 \ifx\csname\exstring@#1@C\endcsname\relax
  \expandafter\ifx\csname\exstring@#1@C1\endcsname\relax
   \DN@{\Err@{\noexpand\Offset can't be used with \string#1}}%
  \else
   \DN@{\Offset@#1}%
  \fi
 \else
  \DN@##1{\count@##1\relax\advance\count@\m@ne\global\advance
   \csname\exstring@#1@C\endcsname\count@}%
 \fi
 \next@}
\def\Offset@#1#2{\expandafter
 \ifx\csname\exstring@#1@C#2\endcsname\relax
  \DN@{\Err@{\string\Offset\string#1 can't be followed by \string#2}}%
 \else
  \DN@##1{\count@##1\relax\advance\count@\m@ne
   \global\advance\csname\exstring@#1@C#2\endcsname\count@}%
 \fi
 \next@}
\def\getR@#1#2{\def\nextiv@{\let\nextiii@}\expandafter\nextiv@
 \csname\exstring@#1@R#2\endcsname}
\def\letR@#1#2#3{\expandafter\let\csname#1@#3#2\endcsname\Next@}
\def\letR@@#1#2{\expandafter\let\csname\exstring@#1@#2\endcsname\Next@}
\def\newpre#1{\expandafter
 \ifx\csname\exstring@#1@P\endcsname\relax
  \expandafter\ifx\csname\exstring@#1@P1\endcsname\relax
   \DN@{\Err@{\noexpand\newpre can't be used with \string#1}}%
  \else
   \DN@{\newpre@#1}%
  \fi
 \else
  \DN@{%
   \DNii@{%
    \endgroup
    \expandafter\let\csname\exstring@#1@P\endcsname\Next@
    \expandafter\ifx\csname\exstring@#1@R\endcsname\relax\else
    \getR@#1{}\expandafter\letR@\nextiii@ P\fi
    }%
   \begingroup\noexpands@\afterassignment\nextii@\xdef\Next@}%
 \fi
 \next@}
\def\newpre@#1#2{\expandafter
 \ifx\csname\exstring@#1@P#2\endcsname\relax
  \DN@{\Err@{\string\newpre\string#1 can't be followed by \string#2}}%
 \else
  \DN@{%
   \DNii@{%
    \endgroup
    \expandafter\let\csname\exstring@#1@P#2\endcsname\Next@
    \expandafter\ifx\csname\exstring@#1@R#2\endcsname\relax\else
    \getR@#1{#2}\expandafter\letR@@\nextiii@ P\fi
    }%
   \begingroup\noexpands@\afterassignment\nextii@\xdef\Next@}%
 \fi
 \next@}
\def\newpost#1{\expandafter
 \ifx\csname\exstring@#1@Q\endcsname\relax
  \expandafter\ifx\csname\exstring@#1@Q1\endcsname\relax
   \DN@{\Err@{\noexpand\newpost can't be used with \string#1}}%
  \else
   \DN@{\newpost@#1}%
  \fi
 \else
  \DN@{%
   \DNii@{%
    \endgroup
    \expandafter\let\csname\exstring@#1@Q\endcsname\Next@
    \expandafter\ifx\csname\exstring@#1@R\endcsname\relax\else
    \getR@#1{}\expandafter\letR@\nextiii@ Q\fi
    }%
   \begingroup\noexpands@\afterassignment\nextii@\xdef\Next@}%
 \fi
 \next@}
\def\newpost@#1#2{\expandafter
 \ifx\csname\exstring@#1@Q#2\endcsname\relax
  \DN@{\Err@{\string\newpost\string#1 can't be followed by \string#2}}%
 \else
  \DN@{%
   \DNii@{%
    \endgroup
    \expandafter\let\csname\exstring@#1@Q#2\endcsname\Next@
    \expandafter\ifx\csname\exstring@#1@R#2\endcsname\relax\else
    \getR@#1{#2}\expandafter\letR@@\nextiii@ Q\fi
    }%
   \begingroup\noexpands@\afterassignment\nextii@\xdef\Next@}%
 \fi
 \next@}
\def\newstyle#1{\expandafter
 \ifx\csname\exstring@#1@S\endcsname\relax
  \expandafter\ifx\csname\exstring@#1@S1\endcsname\relax
   \DN@{\Err@{\noexpand\newstyle can't be used
    with \string#1}}%
  \else
   \DN@{\newstyle@#1}%
  \fi
 \else
  \DN@{%
   \DNii@{%
    \expandafter\let\csname\exstring@#1@S\endcsname\Next@
    \expandafter\ifx\csname\exstring@#1@R\endcsname\relax\else
    \getR@#1{}\expandafter\letR@\nextiii@ S\fi
    }%
   \afterassignment\nextii@\gdef\Next@}%
 \fi
 \next@}
\def\newstyle@#1#2{\expandafter
 \ifx\csname\exstring@#1@S#2\endcsname\relax
  \DN@{\Err@{\string\newstyle\string#1 can't be followed by
   \string#2}}%
 \else
  \DN@{%
   \DNii@{%
    \expandafter\let\csname\exstring@#1@S#2\endcsname\Next@
    \expandafter\ifx\csname\exstring@#1@R#2\endcsname\relax\else
    \getR@#1{#2}\expandafter\letR@@\nextiii@ S\fi
    }%
   \afterassignment\nextii@\gdef\Next@}%
 \fi
 \next@}
\def\newnumstyle#1{\expandafter
 \ifx\csname\exstring@#1@N\endcsname\relax
  \expandafter\ifx\csname\exstring@#1@N1\endcsname\relax
   \DN@{\Err@{\noexpand\newnumstyle can't be used with
    \string#1}}%
  \else
   \DN@{\newnumstyle@#1}%
  \fi
 \else
  \DN@##1{%
   \gdef\Next@{##1}%
    \expandafter\let\csname\exstring@#1@N\endcsname\Next@
    \expandafter\ifx\csname\exstring@#1@R\endcsname\relax\else
    \getR@#1{}\expandafter\letR@\nextiii@ N\fi
    }%
 \fi
 \next@}
\def\newnumstyle@#1#2{\expandafter
 \ifx\csname\exstring@#1@N#2\endcsname\relax
  \DN@{\Err@{\string\newnumstyle\string#1 can't be followed by
   \string#2}}%
 \else
  \DN@##1{%
   \gdef\Next@{##1}%
    \expandafter\let\csname\exstring@#1@N#2\endcsname\Next@
    \expandafter\ifx\csname\exstring@#1@R#2\endcsname\relax\else
    \getR@#1{#2}\expandafter\letR@@\nextiii@ N\fi
    }%
  \fi
 \next@}
\def\newfontstyle#1{\expandafter
 \ifx\csname\exstring@#1@F\endcsname\relax
  \expandafter\ifx\csname\exstring@#1@F1\endcsname\relax
   \DN@{\Err@{\noexpand\newfontstyle can't be used with
    \string#1}}%
  \else
   \DN@{\newfontstyle@#1}%
  \fi
 \else
  \DN@##1{%
   \gdef\Next@{##1}%
    \expandafter\let\csname\exstring@#1@F\endcsname\Next@
    \expandafter\ifx\csname\exstring@#1@R\endcsname\relax\else
    \getR@#1{}\expandafter\letR@\nextiii@ F\fi
    }%
 \fi
 \next@}
\def\newfontstyle@#1#2{\expandafter
 \ifx\csname\exstring@#1@F#2\endcsname\relax
  \DN@{\Err@{\string\newfontstyle\string#1 can't be followed by
   \string#2}}%
 \else
  \DN@##1{%
   \gdef\Next@{##1}%
    \expandafter\let\csname\exstring@#1@F#2\endcsname\Next@
    \expandafter\ifx\csname\exstring@#1@R#2\endcsname\relax\else
    \getR@#1{#2}\expandafter\letR@@\nextiii@ F\fi
    }%
 \fi
 \next@}
\def\word#1{\expandafter
 \ifx\csname\exstring@#1@W\endcsname\relax
  \expandafter\ifx\csname\exstring@#1@W1\endcsname\relax
   \DN@{\Err@{\noexpand\word can't be used with \string#1}}%
  \else
   \DN@{\word@#1}%
  \fi
 \else
  \DN@{{\csname\exstring@#1@W\endcsname}}%
 \fi
 \next@}
\def\word@#1#2{\expandafter
 \ifx\csname\exstring@#1@W#2\endcsname\relax
  \DN@{\Err@{\string\word\noexpand#1can't be followed by \string#2}}%
 \else
  \DN@{{\csname\exstring@#1@W#2\endcsname}}%
 \fi
 \next@}
\def\newword#1{\expandafter
 \ifx\csname\exstring@#1@W\endcsname\relax
  \expandafter\ifx\csname\exstring@#1@W1\endcsname\relax
   \DN@{\Err@{\noexpand\newword can't be used  with \string#1}}%
  \else
   \DN@{\newword@#1}%
  \fi
 \else
  \DN@{%
   \DNii@{%
    \expandafter\let\csname\exstring@#1@W\endcsname\Next@
    \expandafter\ifx\csname\exstring@#1@R\endcsname\relax\else
     \getR@#1{}\expandafter\letR@\nextiii@ W\fi
    }%
   \afterassignment\nextii@\gdef\Next@}%
 \fi
 \next@}
\def\newword@#1#2{\expandafter
 \ifx\csname\exstring@#1@W#2\endcsname\relax
  \DN@{\Err@{\string\newword\noexpand#1can't be followed by \string#2}}%
 \else
  \DN@{%
   \DNii@{%
    \expandafter\let\csname\exstring@#1@W#2\endcsname\Next@
    \expandafter\ifx\csname\exstring@#1@R#2\endcsname\relax\else
     \getR@#1{#2}\expandafter\letR@@\nextiii@ W\fi
    }%
   \afterassignment\nextii@\gdef\Next@}%
 \fi
 \next@}
\newif\iffn@
\newcount\footmark@C
\footmark@C\z@
\def\footmark@S#1{$^{#1}$}
\let\footmark@N\arabic
\def\footmark@F{\rm}
\def\foottext@S#1{$^{#1}$}
\def\foottext@F{\rm}
\let\modifyfootnote@\relax
\def\modifyfootnote#1{\def\modifyfootnote@{#1}}
\def\vfootnote@#1{\insert\footins
 \bgroup
 \floatingpenalty\@MM\interlinepenalty\interfootnotelinepenalty
 \leftskip\z@\rightskip\z@\spaceskip\z@\xspaceskip\z@
 \rm\splittopskip\ht\strutbox\splitmaxdepth\dp\strutbox
 \locallabel@\noindent@@{\foottext@F#1}\modifyfootnote@
 \footstrut\FN@\fo@t}
\def\fo@t{\ifcat\bgroup\noexpand\next\expandafter\f@@t\else
 \expandafter\f@t\fi}
\def\f@t#1{#1\@foot}
\def\f@@t{\bgroup\aftergroup\@foot\afterassignment\FNSSP@\let\next@}
\def\@foot{\unskip\lower\dp\strutbox\vbox to\dp\strutbox{}\egroup
 \iffn@\expandafter\fn@false\else
 \expandafter\postvanish@\fi}
\newif\ifplainfn@
\plainfn@true
\def\fancyfootnotes{\plainfn@false}
\newcount\fancyfootmarkcount@
\fancyfootmarkcount@\z@
\newcount\lastfnpage@
\lastfnpage@-\@M
\let\justfootmarklist@\empty
\def\footmark{\let\@sf\empty
 \ifhmode\edef\@sf{\spacefactor\the\spacefactor}\/\fi
 \DN@{\ifx"\next\expandafter\nextii@\else\expandafter\footmark@\fi}%
 \DNii@"##1"{%
  \iffirstchoice@
   {\let\style\footmark@S\let\numstyle\footmark@N
   \footmark@F##1%
   \noexpands@
   \let\style\foottext@S
   \Qlabel@{##1}%
   }%
   \iffn@\else
    {\noexpands@
    \xdef\Next@{{\Thelabel@}{\Thelabel@@}{\Thelabel@@@}{\Thelabel@@@@}}%
    }%
    \expandafter\rightappend@\Next@\to\justfootmarklist@
   \fi
  \fi
  \@sf\relax}%
 \FN@\next@}
\def\footmark@{%
 \iffirstchoice@
  \global\advance\footmark@C\@ne
  \ifplainfn@
   \xdef\adjustedfootmark@{\number\footmark@C}%
  \else
   {\let\\\or\xdef\Next@{\ifcase\number\footmark@C\fnpages@\else
     -\@M\fi}}%
   \ifnum\Next@=-\@M
    \xdef\adjustedfootmark@{\number\footmark@C}%
   \else
    \ifnum\Next@=\lastfnpage@
     \global\advance\fancyfootmarkcount@\@ne
    \else
     \global\fancyfootmarkcount@\@ne
     \global\lastfnpage@\Next@
    \fi
    \xdef\adjustedfootmark@{\number\fancyfootmarkcount@}%
   \fi
  \fi
  {\noexpands@
  \xdef\Thelabel@@@{\adjustedfootmark@}%
  \xdefThelabel@\footmark@N
  \xdef\Thelabel@@@@{\Thelabel@}%
  \xdefThelabel@@\foottext@S
  }%
  \iffn@\else
   {\noexpands@
   \xdef\Next@{{\Thelabel@}{\Thelabel@@}{\Thelabel@@@}{\Thelabel@@@@}}%
   }%
   \expandafter\rightappend@\Next@\to\justfootmarklist@
  \fi
  \ifplainfn@
  \else
   \edef\next@{\write\laxwrite@{F\noexpand\the\pageno}}\next@
  \fi
 \fi
 \footmark@S{\footmark@N{\adjustedfootmark@}}%
 \@sf\relax}
\def\foottext{\prevanish@
 \ifx\justfootmarklist@\empty
  \Err@{There is no \noexpand\footmark for this \string\foottext}\fi
 \DN@\\##1##2\next@{\DN@{##1}\gdef\justfootmarklist@{##2}}%
 \expandafter\next@\justfootmarklist@\next@
 \expandafter\foottext@\next@}
\def\foottext@#1#2#3#4{{\noexpands@
  \xdef\Thelabel@{#1}\xdef\Thelabel@@{#2}%
  \xdef\Thelabel@@@{#3}\xdef\Thelabel@@@@{#4}}%
  \vfootnote@{\thelabel@@}}
\rightadd@\foottext\to\vanishlist@
\def\footnote{\fn@true
 \let\@sf\empty
 \ifhmode\edef\@sf{\spacefactor\the\spacefactor}\/\fi
 \DN@{\ifx"\next\expandafter\nextii@\else\expandafter\nextiii@\fi}%
 \DNii@"##1"{\footmark"##1"\vfootnote@{\let\style\foottext@S
  \let\numstyle\footmark@N##1}}%
 \def\nextiii@{\footmark\vfootnote@{\foottext@S{\footmark@N
  {\adjustedfootmark@}}}}%
 \FN@\next@}
\newdimen\litindent
\litindent20\p@
\newbox\litbox@
\newbox\Litbox@
\newcount\interlitpenalty@
\interlitpenalty@\@M
\newcount\litlines@
{\obeyspaces\gdef\defspace@{\def {\allowbreak\hskip.5emminus.15em}}}
{\obeylines\gdef\letM@{\let^^M\CtrlM@}}
\def\CtrlM@{\egroup
 \ifcase\litlines@\advance\litlines@\@ne\or
 \box\litbox@\advance\litlines@\@ne\else
 \penalty\interlitpenalty@\box\litbox@\fi
 \Lit@}
\def\Lit@{\setbox\litbox@\hbox\bgroup\litdefs@\hskip\litindent}
\newcount\littab@
\littab@8
\def\littab#1{\littab@#1\relax}
{\catcode`\^^I=\active\gdef\letTAB@{\let^^I\TAB@}}
\def\TAB@{\egroup
 \dimen@\wd\litbox@
 \advance\dimen@-\litindent
 \setboxz@h{\tt0}%
 \dimen@ii\littab@\wdz@
 \divide\dimen@\dimen@ii
 \multiply\dimen@\dimen@ii
 \advance\dimen@\littab@\wdz@
 \advance\dimen@\litindent
 \setbox\litbox@\hbox\bgroup\litdefs@\hbox to\dimen@{\unhbox\litbox@\hfil}}
{\catcode`\`=\active\gdef`{\relax\lq}}
\let\litbs@\relax
\let\litbs@@\relax
\def\litbackslash#1{%
 \edef\litbs@{\catcode`\string#1=\z@
 \def\noexpand\litbs@@{\def\expandafter\noexpand\csname\string#1\endcsname
  {\char`\string#1}}}}
\def\litcodes@{\catcode`\\=12
 \catcode`\{=12 \catcode`\}=12
 \catcode`\$=12 \catcode`\&=12
 \catcode`\#=12
 \catcode`\^=12 \catcode`\_=12
 \catcode`\@=12 \catcode`\~=12 \catcode`\"=12
 \catcode`\;=12 \catcode`\:=12 \catcode`\!=12 \catcode`\?=12
 \catcode`\%=12 \litbs@\catcode`\`=\active\obeyspaces\defspace@}
\def\activate@#1#2{{\lccode`\~=`#2%
 \lowercase{%
  \if0#1%
  \gdef\Next@{\def~{\egroup\endgroup\bigskip\vskip-\parskip
   \def\next@{\noindent@@\FN@\pretendspace@}\FNSS@\next@}}\else
  \gdef\Next@{\def~{\egroup\egroup\endgroup}}\fi
  }%
 }}
\def\litdefs@{\let\0\empty\let\1\litdelim@\def\ {\char32 }\litbs@@}%
\def\litdelimiter#1{%
 \edef\litdelim@{\char`#1}%
 \def\lit#1{\leavevmode\begingroup\litcodes@\litdefs@
  \tt\hyphenchar\tentt\m@ne\lit@}%
 \def\lit@##1#1{##1\endgroup\null}%
 \def\Lit#1{\ifhmode$$\abovedisplayskip\bigskipamount
  \abovedisplayshortskip\bigskipamount
  \belowdisplayskip\z@\belowdisplayshortskip\z@
  \postdisplaypenalty\@M
  $$\vskip-\baselineskip\else\bigskip\fi
  \begingroup\litlines@\z@
  \catcode`#1=\active\activate@0#1\Next@
  \def\displaybreak{\egroup\break\litlines@\z@\Lit@}%
  \def\allowdisplaybreak{\egroup\allowbreak\litlines@\z@\Lit@}%
  \def\allowdisplaybreaks{\egroup\allowbreak\interlitpenalty@\z@
   \litlines@\z@\Lit@}%
  \litcodes@\tt\catcode`\^^I=\active\letTAB@
  \obeylines\letM@\Lit@}%
 \def\Litbox##1=#1{\begingroup\ifodd##1\relax\aftergroup\global\fi
  \aftergroup\setbox\aftergroup##1\aftergroup\box\aftergroup\Litbox@
  \def\allowdisplaybreak{\egroup\allowbreak\litlines@\z@\Lit@}%
  \def\allowdisplaybreaks{\egroup\allowbreak\interlitpenalty@\z@
   \litlines@\z@\Lit@}%
  \catcode`#1=\active\activate@1#1\Next@
  \litcodes@\tt\catcode`\^^I=\active\letTAB@
  \obeylines\letM@\global\setbox\Litbox@\vbox\bgroup\litindent\z@%
  \litlines@\z@\Lit@}%
}
\newbox\titlebox@
\setbox\titlebox@\vbox{}
\rightadd@\title\to\overlonglist@
\def\title{\begingroup\Let@
 \global\setbox\titlebox@\vbox\bgroup\tabskip\hss@
 \halign to\hsize\bgroup\bf\hfil\ignorespaces##\unskip\hfil\cr}
\def\endtitle{\crcr\egroup\egroup\endgroup\overlong@false}
\newbox\authorbox@
\rightadd@\author\to\overlonglist@
\def\author{\begingroup\Let@
 \global\setbox\authorbox@\vbox\bgroup\tabskip\hss@
 \halign to\hsize\bgroup\rm\hfil\ignorespaces##\unskip\hfil\cr}
\def\endauthor{\crcr\egroup\egroup\endgroup\overlong@false}
\newbox\affilbox@
\def\affil{\begingroup\Let@
 \global\setbox\affilbox@\vbox\bgroup\tabskip\hss@
 \halign to\hsize\bgroup\rm\hfil\ignorespaces##\unskip\hfil\cr}%
\def\endaffil{\crcr\egroup\egroup\endgroup\overlong@false}
\let\date@\relax
\def\date#1{\gdef\date@{\ignorespaces#1\unskip}}
\def\today{\ifcase\month\or January\or February\or March\or April\or May\or
 June\or July\or August\or September\or October\or November\or December\fi
 \space\number\day, \number\year}
\def\maketitle{\hrule\height\z@\vskip-\topskip
 \vskip24\p@ plus12\p@ minus12\p@
 \unvbox\titlebox@
 \ifvoid\authorbox@\else\vskip12\p@ plus6\p@ minus3\p@\unvbox\authorbox@\fi
 \ifvoid\affilbox@\else\vskip10\p@ plus5\p@ minus2\p@\unvbox\affilbox@\fi
 \ifx\date@\relax\else\vskip6\p@ plus2\p@ minus\p@\centerline{\rm\date@}\fi
 \vskip18\p@ plus12\p@ minus6\p@}
\def\cite{%
 \DNii@(##1)##2{{\rm[}{##2}, {##1\/}{\rm]}}%
 \def\nextiii@##1{{\rm[}{##1\/}{\rm]}}%
 \DN@{\ifx\next(\expandafter\nextii@\else\expandafter\nextiii@\fi}%
 \FN@\next@}
\def\makebib@W{Bibliography}
\def\makebibp@W{References}
\def\makebib{\begingroup\rm\bigbreak\centerline{\smc\makebib@W}%
 \nobreak\medskip
 \sfcode`\.=\@m\everypar{}\parindent\z@
 \def\nopunct{\nopunct@true}\def\nospace{\nospace@true}%
 \nopunct@false\nospace@false
 \def\lkerns@{\null\kern\m@ne sp\kern\@ne sp}%
 \def\nkerns@{\null\kern-\tw@ sp\kern\tw@ sp}%
}

\newif\ifnoprepunct@
\newif\ifnoprespace@
\newif\ifnoquotes@
\def\noprepunct{\noprepunct@true}
\def\noprespace{\noprespace@true}
\def\noquotes{\noquotes@true}
\newbox\nobox@
\newbox\keybox@
\newbox\bybox@
\newbox\paperbox@
\newbox\paperinfobox@
\newbox\jourbox@
\newbox\volbox@
\newbox\issuebox@
\newbox\yrbox@
\newbox\pgbox@
\newbox\ppbox@
\newbox\bookbox@
\newbox\inbookbox@
\newbox\bookinfobox@
\newbox\publbox@
\newbox\publaddrbox@
\newbox\edbox@
\newbox\edsbox@
\newbox\langbox@
\newbox\translbox@
\newbox\finalinfobox@
\def\setbibinfo@#1{\edef\next@{\ifnopunct@1\else0\fi
 \ifnospace@1\else0\fi\ifnoprepunct@1\else0\fi\ifnoprespace@1\else0\fi
 \ifnoquotes@1\else0\fi}%
 \DNii@{00000}%
 \ifx\next@\nextii@\else\xdef\bibinfo@{\bibinfo@\the#1,\next@}%
 \fi}
\def\getbibinfo@#1{\ifx\bibinfo@\empty
 \let\next@0\let\nextii@0\let\nextiii@0\let\nextiv@0\let\nextv@0\else
 \edef\next@{\def
  \noexpand\next@####1\the#1,####2####3####4####5####6####7\noexpand\next@
  {\let\noexpand\next@####2\let\noexpand\nextii@####3%
  \let\noexpand\nextiii@####4\let\noexpand\nextiv@####5%
  \let\noexpand\nextv@####6}%
  \noexpand\next@\bibinfo@\the#1,00000\noexpand\next@}\next@
 \fi}
\newif\ifbookinquotes@
\def\bookinquotes{\bookinquotes@true}
\newif\ifpaperinquotes@
\def\paperinquotes{\paperinquotes@true}
\newif\ifininbook@
\def\ininbook{\ininbook@true}
\newif\ifopenquotes@
\def\closequotes@{\ifopenquotes@''\openquotes@false\fi}
\newif\ifbeginbib@
\newif\ifendbib@
\newif\ifprevjour@
\newif\ifprevbook@
\newdimen\bibindent@
\bibindent@20\p@
\def\bib{\global\let\bibinfo@\empty\global\let\translinfo@\relax\beginbib@true
 \begingroup\noindent@
 \hangindent\bibindent@\hangafter\@ne\bib@}
\def\v@id#1{\setbox#1\box\voidb@x}
\def\bib@{\v@id\nobox@\v@id\keybox@\v@id\bybox@\v@id\paperbox@
 \v@id\paperinfobox@\v@id\jourbox@\v@id\volbox@\v@id\issuebox@
 \v@id\yrbox@\v@id\pgbox@\v@id\ppbox@\v@id\bookbox@\v@id\inbookbox@
 \v@id\bookinfobox@\v@id\publbox@\v@id\publaddrbox@\v@id\edbox@
 \v@id\edsbox@\v@id\langbox@\v@id\translbox@\v@id\finalinfobox@
 \bgroup}
\def\Setnonemptybox@#1#2{\unskip\setbibinfo@#1\egroup#2%
 \def\aftergroup@{\ifdim\wd#1=\z@\setbox#1\box\voidb@x\fi}%
 \setbox#1\vbox\bgroup\aftergroup\aftergroup@\hsize\maxdimen\leftskip\z@
 \rightskip\z@\hbadness\@M\hfuzz\maxdimen\noindent}
\def\setnonemptybox@#1{\Setnonemptybox@#1\relax}
\def\no{\setnonemptybox@\nobox@}
\def\key{\setnonemptybox@\keybox@\bf}
\def\by{\setnonemptybox@\bybox@}
\def\bysame{\setnonemptybox@\bybox@\leaders\hrule\hskip3em\null}
\def\paper{\setnonemptybox@\paperbox@
 \ifpaperinquotes@\getbibinfo@\paperbox@
 \if\nextv@1\else``\fi\else\it\fi}
\def\paperinfo{\setnonemptybox@\paperinfobox@}
\def\jour{\Setnonemptybox@\jourbox@\prevjour@true}
\def\vol{\setnonemptybox@\volbox@\bf}
\def\issue{\setnonemptybox@\issuebox@}
\def\yr{\setnonemptybox@\yrbox@}
\def\toappear{\noprepunct\finalinfo(to appear)}
\def\pg{\setnonemptybox@\pgbox@}
\def\pp{\setnonemptybox@\ppbox@}
\def\book{\Setnonemptybox@\bookbox@\prevbook@true
 \ifbookinquotes@\getbibinfo@\bookbox@
 \if\nextv@1\else``\fi\else\it\fi}
\def\inbook{\Setnonemptybox@\inbookbox@\prevbook@true
 \ifininbook@ in \fi\ifbookinquotes@\getbibinfo@\inbookbox@
 \if\nextv@1\else``\fi\fi}
\def\bookinfo{\setnonemptybox@\bookinfobox@}
\def\publ{\setnonemptybox@\publbox@}
\def\publaddr{\setnonemptybox@\publaddrbox@}
\def\ed{\setnonemptybox@\edbox@}
\def\eds{\setnonemptybox@\edsbox@}
\def\lang{\setnonemptybox@\langbox@}
\def\finalinfo{\setnonemptybox@\finalinfobox@}
\def\setboxzl@{\setbox\z@\lastbox}
\def\getbox@#1{\setbox\z@\vbox{\vskip-\@M\p@
 \unvbox#1%
 \setboxzl@
 \global\setbox\@ne\hbox{\unhbox\z@\unskip\unskip\unpenalty}%
 \ifdim\lastskip=-\@M\p@\else
 \loop\ifdim\lastskip=-\@M\p@
 \else\unskip\unpenalty\setboxzl@
 \global\setbox\@ne\hbox{\unhbox\z@\unhbox\@ne}%
 \repeat\fi}%
 \unhbox\@ne}
\def\adjustpunct@#1{\count@\lastkern
 \ifnum\count@=\z@#1\closequotes@\else
 \ifnum\count@>\tw@#1\closequotes@\else
 \ifnum\count@<-\tw@#1\closequotes@\else
  \unkern\unkern\setboxzl@
  \skip@\lastskip\unskip
  \count@@\lastpenalty\unpenalty
  \ifnum\count@=\tw@\unskip\setboxzl@\fi
  \ifdim\skip@=\z@\else\hskip\skip@\fi
  #1\closequotes@
  \ifnum\count@=\tw@\null\hfill\fi
  \penalty\count@@
 \fi\fi\fi}
\def\prepunct@#1#2{\getbibinfo@#2%
 \ifnopunct@
 \else
  \if\nextiii@0\adjustpunct@#1\fi
 \fi
 \closequotes@
 \ifnospace@
 \else
  \if\nextiv@0\space\else\fi
 \fi
 \nopunct@false\nospace@false
 \if\next@1\nopunct@true\fi
 \if\nextii@1\nospace@true\fi}
\def\ppunbox@#1#2{\prepunct@{#1}#2%
 \getbox@#2}
\let\semicolon@;
\def\endbib@{%
 \ifbeginbib@
  \ifvoid\nobox@
   \ifvoid\keybox@\else\hbox to\bibindent@{[\getbox@\keybox@]\hss}\fi
  \else\hbox to\bibindent@{\hss\getbox@\nobox@. }\fi
  \ifvoid\bybox@\else\getbox@\bybox@\fi
 \else
  \nopunct@true
  \ifvoid\bybox@\else\ppunbox@\relax\bybox@\fi
 \fi
 \ifvoid\translbox@\else\ppunbox@,\translbox@\fi
 \ifvoid\paperbox@\else\ppunbox@,\paperbox@\ifpaperinquotes@
  \if\nextv@1\else\openquotes@true\fi\fi
 \fi
 \ifvoid\paperinfobox@\else\ppunbox@,\paperinfobox@\fi
 \test@false
 \ifvoid\jourbox@\else\test@true\ppunbox@,\jourbox@\fi
 \ifprevjour@\test@true\fi
 \iftest@
  \ifvoid\volbox@\else\ppunbox@\relax\volbox@\fi
  \ifvoid\issuebox@
   \else\prepunct@\relax\issuebox@ no.~\getbox@\issuebox@\fi
  \ifvoid\yrbox@\else\prepunct@\relax\yrbox@(\getbox@\yrbox@)\fi
  \ifvoid\ppbox@\else\ppunbox@,\ppbox@\fi
  \ifvoid\pgbox@\else\prepunct@,\pgbox@ p.~\getbox@\pgbox@\fi
 \fi
 \test@false
 \ifvoid\bookbox@\else\test@true\ppunbox@,\bookbox@\ifbookinquotes@
  \if\nextv@1\else\openquotes@true\fi\fi\fi
 \ifvoid\inbookbox@\else\test@true\ppunbox@,\inbookbox@\ifbookinquotes@
  \if\nextv@1\else\openquotes@true\fi\fi\fi
 \ifprevbook@\test@true\fi
 \iftest@
  \ifvoid\edbox@\else\prepunct@\relax\edbox@(\getbox@\edbox@, ed.)\fi
  \ifvoid\edsbox@\else\prepunct@\relax\edsbox@(\getbox@\edsbox@, eds.)\fi
  \ifvoid\bookinfobox@\else\ppunbox@,\bookinfobox@\fi
  \ifvoid\publbox@\else\ppunbox@,\publbox@\fi
  \ifvoid\publaddrbox@\else\ppunbox@,\publaddrbox@\fi
  \ifvoid\yrbox@\else\ppunbox@,\yrbox@\fi
  \ifvoid\ppbox@\else\prepunct@,\ppbox@ pp.~\getbox@\ppbox@\fi
  \ifvoid\pgbox@\else\prepunct@,\pgbox@ p.~\getbox@\pgbox@\fi
 \fi
 \ifvoid\finalinfobox@
  \ifendbib@
   \ifnopunct@\else.\closequotes@\fi
  \else
  \ifvoid\langbox@\else\space(\getbox@\langbox@)\fi
   \/\semicolon@\closequotes@
  \fi
 \else
  \ifendbib@
   \ppunbox@{.\spacefactor3000\relax}\finalinfobox@
    \ifnopunct@\else.\fi
  \else
   \ppunbox@,\finalinfobox@\/\semicolon@\fi
 \fi
 \ifvoid\langbox@\else\space(\getbox@\langbox@)\fi
}
\def\endbib{\unskip\egroup\endbib@true\endbib@\par\endgroup}
\def\morebib{\unskip\egroup
 \endbib@false\endbib@
 \global\let\bibinfo@\empty\beginbib@false
 \bib@}
\def\anotherbib{\unskip\egroup
 \endbib@false\endbib@
 \global\let\bibinfo@\empty\beginbib@false
 \prevjour@false\prevbook@false\bib@}
\def\transl{\unskip
 \xdef\translinfo@{\the\translbox@,\ifnopunct@1\else0\fi
 \ifnospace@1\else0\fi\ifnoprepunct@1\else0\fi\ifnoprespace@1\else0\fi0}%
 \egroup\endbib@false\endbib@
 \global\let\bibinfo@\translinfo@\beginbib@false
 \bib@
 \egroup
 \def\aftergroup@{\ifdim\wd\translbox@=\z@\setbox\translbox@\box\voidb@x\fi}%
 \setbox\translbox@\vbox\bgroup\aftergroup\aftergroup@
 \hsize\maxdimen\leftskip\z@\rightskip\z@\hbadness\@M\hfuzz\maxdimen
 \noindent}
\newwrite\auxwrite@
\newread\bbl@
\def\UseBibTeX{\immediate\openout\auxwrite@=\jobname.aux
 \let\cite\BTcite@
 \def\nocite##1{\immediate\write\auxwrite@{\string\citation{##1}}}%
 \def\bibliographystyle##1{\immediate\write\auxwrite@{\string
  \bibstyle{##1}}}%
 \def\bibliography@W{Bibliography}%
 \def\bibliography##1{\immediate\write\auxwrite@{\string\bibdata{##1}}%
  \immediate\openin\bbl@=\jobname.bbl
  \ifeof\bbl@
   \W@{No .bbl file}%
  \else
   \immediate\closein\bbl@
   \begingroup\input bibtex \input\jobname.bbl \endgroup
  \fi}%
 }
\def\BTcite@{%
 \DNii@(##1)##2{{\rm[}\BTcite@@##2,\BTcite@@{\rm, }{##1\/}{\rm]}%
  \immediate\write\auxwrite@{\string\citation{##2}}}%
 \def\nextiii@##1{{\rm[}\BTcite@@##1,\BTcite@@\/{\rm]}%
  \immediate\write\auxwrite@{\string\citation{##1}}}%
 \DN@{\ifx\next(\expandafter\nextii@\else\expandafter\nextiii@\fi}%
 \FN@\next@}%
\def\BTcite@@#1,{\BTcite@@@{#1}\FN@\BTcite@@@@}
\def\BTcite@@@@{\ifx\next\BTcite@@
 \expandafter\eat@\else{\rm, }\expandafter\BTcite@@\fi}
\catcode`\~=11
\def\BTcite@@@#1{\nolabel@\cite{#1}\relax
 \DNii@##1~##2\nextii@{##1}%
 \csL@{#1}\expandafter\nextii@\Next@\nextii@\fi}
\catcode`\~=\active

\def\beginthebibliography@#1{\rm\setboxz@h{#1\ }\bibindent@\wdz@
 \bigbreak\centerline{\smc\bibliography@W}\nobreak\medskip
 \sfcode`\.=\@m\everypar{}\parindent\z@}

\newif\iffigproofing@
\def\Figureproofing{\figproofing@true}
\def\noFigureproofing{\figproofing@false}
\newif\ifHby@
\def\Hbyw#1{\global\Hby@true\hbyw\vsize{#1}}
\def\hbyw#1#2{%
 \hbox{%
  \ifHby@
  \else
   \iffigproofing@
    \setbox\z@\vbox{\hrule\width5\p@}\ht\z@\z@
    \vbox to#1{\hrule\height5\p@\width.4\p@\vfil\hrule\height5\p@\width.4\p@}%
    \kern-.4\p@\rlap{\copy\z@}\raise#1\hbox{\rlap{\copy\z@}}%
   \fi
  \fi
  \vbox to#1{\hbox to#2{}\vfil}%
  \ifHby@
  \else
   \iffigproofing@
    \vbox to#1{\hrule\height5\p@\width.4\p@\vfil\hrule\height5\p@\width.4\p@}%
    \kern-.4\p@\llap{\copy\z@}\raise#1\hbox{\llap{\boxz@}}%
   \fi
  \fi}}
\newcount\island@C
\let\island@P\empty
\let\island@Q\empty
\def\island@S#1{#1\null.}
\let\island@N\arabic
\def\island@F{\rm}
\def\island@@@P{\csname\exxx@\islandtype@ @P\endcsname}
\def\island@@@Q{\csname\exxx@\islandtype@ @Q\endcsname}
\def\island@@@S{\csname\exxx@\islandtype@ @S\endcsname}
\def\island@@@N{\csname\exxx@\islandtype@ @N\endcsname}
\def\island@@@F{\csname\exxx@\islandtype@ @F\endcsname}
\def\island@@@C{\csname island@C\islandclass@\endcsname}
\newif\ifplace@
\newif\ifisland@
\def\island{%
 \ifplace@
  \DN@{\let\islandclass@\empty\def\islandtype@{\island}\FN@\island@}%
 \else
  \long\DN@##1\endisland{\Err@{\noexpand\island must be used after some
   type of \string\...place}}%
 \fi
 \next@}
\def\island@{\ifx\next\c\let\next@\island@c\else
 \DN@{\FN@\island@@}\fi\next@}
\def\island@@{\ifcat\bgroup\noexpand\next\let\next@\island@@@\else
 \DN@{\Err@{\noexpand\island must be followed by a {prefix} for
 \string\caption's}}\fi\next@}
\newbox\islandbox@
\newcount\captioncount@
\def\island@@@#1{\def\captionprefix@{#1}\captioncount@\z@
 \global\setbox\islandbox@\vbox\bgroup}
\def\island@c\c#1{%
 \ifplace@
 \DN@{\def\islandclass@{#1}%
  \expandafter\ifx\csname island@C#1\endcsname\relax
  \expandafter\newcount@\csname island@C#1\endcsname
   \global\csname island@C#1\endcsname\z@\fi
  \FNSS@\island@c@}%
 \else
 \DN@{\edef\next@{\long\def\noexpand\next@########1\expandafter\noexpand
  \csname end\exxx@\islandtype@\endcsname{\noexpand\Err@{\noexpand\noexpand
  \expandafter\noexpand
  \islandtype@ must be used after some type of \noexpand\string
   \noexpand\...place}}}\next@\next@}%
 \fi
 \next@}
\def\island@c@{%
 \ifcat\bgroup\noexpand\next
  \let\next@\island@c@@
 \else
  \DN@{\Err@{\noexpand\island\string\c{\expandafter\string\islandclass@} must
   be followed by a {prefix} for \string\caption's}}%
 \fi\next@}
\def\island@c@@#1{\def\captionprefix@{#1}%
 \captioncount@\z@\global\setbox\islandbox@\vbox\bgroup}
\rightadd@\caption\to\nofrillslist@
\newbox\captionbox@
\newbox\Captionbox@
\def\caption{%
 \ifnum\captioncount@=\z@
  \ifnopunct@
   \DN@{\egroup\nopunct@true}%
  \else
   \let\next@\egroup
  \fi
 \else
  \let\next@\relax
 \fi
 \next@
 \advance\captioncount@\@ne
 \FN@\caption@}
\def\caption@{\ifx\next"\expandafter\caption@q\else\expandafter\caption@@\fi}
\def\caption@q"#1"{\quoted@true
 {\noexpands@
 \let\pre\island@@@P\let\post\island@@@Q
 \let\style\island@@@S\let\numstyle\island@@@N
 \Qlabel@{#1}\let\style\relax\xdef\Qlabel@@@@{#1}}%
 \finishcaption@}
\def\caption@@{\quoted@false
 \global\advance\island@@@C\@ne
 {\noexpands@
 \xdef\Thelabel@@@{\number\island@@@C}%
 \xdefThelabel@\island@@@N
 \xdef\Thelabel@@@@{\island@@@P\Thelabel@\island@@@Q}%
 \xdefThelabel@@\island@@@S
 \xdef\Thepref@{\Thelabel@@@@}}%
 \finishcaption@}
\long\def\captionformat@#1#2#3{\rm\strut#1 {\island@@@F#2} #3%
 \punct@.\strut}
\long\def\widerthanisland@#1#2#3{\test@true\setbox\z@\vbox{\hsize\maxdimen
 \noindent@@\captionformat@{#1}{#2}{#3}\par\setboxzl@}%
 \ifdim\wdz@=\z@
  \global\setbox\captionbox@\hbox{\noset@\unlabel@
   \captionformat@{#1}{#2}{#3}}%
  \ifdim\wd\captionbox@>\wd\islandbox@\else\test@false\fi
 \fi}
\long\def\captionformat@@#1#2#3{\widerthanisland@{#1}{#2}{#3}%
 \iftest@
  \global\setbox\captionbox@\vbox{\hsize\wd\islandbox@
   \vskip-\parskip\noindent@@\noset@\unlabel@
   \captionformat@{#1}{#2}{#3}\par}%
 \else
  \global\setbox\captionbox@
   \hbox to\wd\islandbox@{\hfil\box\captionbox@\hfil}%
 \fi}
\long\def\finishcaption@#1{\def\entry@{#1}%
 {\locallabel@
 \captionformat@@
  {\expandafter\ignorespaces\captionprefix@\unskip}%
  {\ifx\thelabel@@\empty\unskip\else\thelabel@@\fi}%
  {\ignorespaces#1\unskip}%
 \ifnum\captioncount@=\@ne
  \global\setbox\islandbox@\vbox{\ticwrite@\vbox{\box\islandbox@}}%
  \global\setbox\Captionbox@\vbox{\box\captionbox@}%
 \else
  \global\setbox\islandbox@\vbox{\unvbox\islandbox@\setboxzl@
   \ticwrite@\boxz@}%
  \global\setbox\Captionbox@\vbox{\unvbox\Captionbox@
   \smallskip\box\captionbox@}%
 \fi}%
 \nopunct@false\nospace@false\ignorespaces}
\def\Sixtic@{\ifx\macdef@\empty\else
 \DN@##1##2\next@{\def\macdef@{##1##2}}%
 \expandafter\next@\macdef@\next@
 \edef\next@
  {\noexpand\six@\tic@\macdef@
  \space\space\space\space\space\space\space\space\space\space\space\space
  \noexpand\six@}%
 \next@\let\macdef@\relax\fi}
\def\ticwrite@{%
 \iftoc@
  {\noexpands@\let\style\relax
  \DN@{\island}%
  \edef\next@{\write\tic@{%
   \ifnopunct@\noexpand\noexpand\noexpand\nopunct\fi
   \ifx\islandtype@\next@\noexpand\noexpand\noexpand\island
    \noexpand\string\noexpand\c{\islandclass@}{\captionprefix@}%
     {\QorThelabel@@@@}\else\noexpand\noexpand\expandafter\noexpand
     \islandtype@{\QorThelabel@@@@}}\fi}%
  \next@}%
  \expandafter\unmacro@\meaning\entry@\unmacro@
  \Sixtic@
  \write\tic@{\noexpand\Page{\number\pageno}{\page@N}{\page@P}{\page@Q}^^J}%
 \fi}
\def\Htrim@#1{%
 \ifHby@
  \dimen@\vsize
  \ifnum\captioncount@=\z@
  \else
   \advance\dimen@-\ht\Captionbox@
   \advance\dimen@-#1%
  \fi
  \global\Hby@false
  \dimen@ii\wd\islandbox@
  \global\setbox\islandbox@\vbox
   {\unvbox\islandbox@\setboxzl@
   \vbox to\z@{\vss\boxz@}\nointerlineskip\hbyw\dimen@\dimen@ii}%
  \global\Hby@true
 \fi}
\newif\ifdata@
\def\iclasstest@#1{\DN@{#1}\ifx\next@\islandclass@
 \test@true\else\test@false\fi}
\skipdef\skipi@=1
\def\endisland{\ifnum\captioncount@=\z@\expandafter\egroup\fi
 \ifdata@
 \else
  \iclasstest@{T}%
  \iftest@
   {\rm\global\skipi@-\dp\strutbox}\global\advance\skipi@\bigskipamount
   \Htrim@\skipi@
   \global\setbox\islandbox@\vbox
    {\ifnum\captioncount@=\z@\else
     \box\Captionbox@
     \nointerlineskip
     \vskip\skipi@\fi
     \box\islandbox@}%
  \else
   {\rm\global\skipi@\dp\strutbox}\global\advance\skipi@\medskipamount
   \Htrim@\skipi@
   \global\setbox\islandbox@\vbox
    {\box\islandbox@
     \ifnum\captioncount@=\z@\else
     \nointerlineskip
     \vskip\skipi@
     \box\Captionbox@
     \fi}%
  \fi
  \ifHby@
  \else
   \dimen@\ht\islandbox@\advance\dimen@\dp\islandbox@
   \ifdim\dimen@>\vsize
    \DN@{\island}%
    \Err@{%
     \ifx\islandtype@\next@\noexpand\island\else
      \expandafter\noexpand\islandtype@\fi
     \ifnum\captioncount@=\z@\else
       with \noexpand\caption\fi
      is larger than page}%
     \ht\islandbox@=\vsize
   \fi
  \fi
 \fi
 \global\Hby@false\island@true}
\def\newisland#1\c#2#3{\define#1{}%
 \iftoc@\immediate\write\tic@{\noexpand\newisland\noexpand#1%
  \string\c{#2}{#3}^^J}\fi
 \expandafter\def\csname\exstring@#1@S\endcsname{\island@S}%
 \expandafter\def\csname\exstring@#1@N\endcsname{\island@N}%
 \expandafter\def\csname\exstring@#1@P\endcsname{\island@P}%
 \expandafter\def\csname\exstring@#1@Q\endcsname{\island@Q}%
 \expandafter\def\csname\exstring@#1@F\endcsname{\island@F}%
 \expandafter\def\csname end\exstring@#1\endcsname{\endisland}%
 \expandafter
 \ifx\csname island@C#2\endcsname\relax
  \expandafter\newcount@\csname island@C#2\endcsname
  \global\csname island@C#2\endcsname\z@
 \fi
 \edef\next@{\noexpand\expandafter\noexpand\let\noexpand
  \csname\exstring@#1@C\noexpand\endcsname
  \csname island@C#2\endcsname}%
 \next@
 \def#1{\def\islandtype@{#1}\island@c\c{#2}{#3}}}
\newisland\Figure\c{F}{Figure}
\newisland\Table\c{T}{Table}
\newbox\islandboxi
\newbox\islandboxii
\newbox\islandboxiii
\newbox\captionboxi
\newbox\captionboxii
\newbox\captionboxiii
\long\def\islandpairdata#1#2{{\data@true
 \place@true
 #1%
 \global\setbox\islandboxi\box\islandbox@
 \global\setbox\captionboxi\box\Captionbox@
 #2%
 \global\setbox\islandboxii\box\islandbox@
 \global\setbox\captionboxii\box\Captionbox@
 }}
\long\def\islandpairbox#1#2{\islandpairdata{#1}{#2}%
 \dimen@\ht\captionboxi
 \ifdim\ht\captionboxii>\dimen@\dimen@\ht\captionboxii\fi
 \ifdim\dimen@>\z@
  \ifdim\ht\captionboxi<\dimen@
   \global\setbox\captionboxi\vbox to\dimen@{\unvbox\captionboxi\vfil}\fi
  \ifdim\ht\captionboxii<\dimen@
   \global\setbox\captionboxii\vbox to\dimen@{\unvbox\captionboxii\vfil}\fi
 \fi
 \global\setbox\islandbox@\vbox
 {\hbox to\hsize{\hfil\box\islandboxi\hfil\box\islandboxii\hfil}%
 \ifdim\dimen@>\z@\nointerlineskip
 {\rm\global\skipi@\dp\strutbox}\global\advance\skipi@\medskipamount
  \vskip\skipi@
  \hbox to\hsize{\hfil\box\captionboxi\hfil\box\captionboxii\hfil}\fi}}
\long\def\islandpairboxa#1#2{\islandpairdata{#1}{#2}%
 \dimen@\ht\captionboxi
 \ifdim\ht\captionboxii>\dimen@\dimen@\ht\captionboxii\fi
 \ifdim\dimen@>\z@
  \ifdim\ht\captionboxi<\dimen@
   \global\setbox\captionboxi\vbox to\dimen@{\vfil\unvbox\captionboxi}\fi
  \ifdim\ht\captionboxii<\dimen@
   \global\setbox\captionboxii\vbox to\dimen@{\vfil\unvbox\captionboxii}\fi
 \fi
 \dimen@ii\ht\islandboxi
 \ifdim\ht\islandboxii>\dimen@ii \dimen@ii\ht\islandboxii\fi
 \ifdim\dimen@ii>\z@
  \ifdim\ht\islandboxi<\dimen@ii
   \global\setbox\islandboxi\vbox to\dimen@ii{\box\islandboxi\vfil}\fi
  \ifdim\ht\islandboxii<\dimen@ii
   \global\setbox\islandboxii\vbox to\dimen@ii{\box\islandboxii\vfil}\fi
 \fi
 \global\setbox\islandbox@\vbox{\ifdim\dimen@>\z@
  \hbox to\hsize{\hfil\box\captionboxi\hfil\box\captionboxii\hfil}%
  \nointerlineskip{\rm\global\skipi@-\dp\strutbox}%
  \global\advance\skipi@\bigskipamount\vskip\skipi@\fi
  \hbox to\hsize{\hfil\box\islandboxi\hfil\box\islandboxii\hfil}}}
\long\def\islandtripledata#1#2#3{{\data@true\place@true
 #1%
 \global\setbox\islandboxi\box\islandbox@
 \global\setbox\captionboxi\box\Captionbox@
 #2%
 \global\setbox\islandboxii\box\islandbox@
 \global\setbox\captionboxii\box\Captionbox@
 #3%
 \global\setbox\islandboxiii\box\islandbox@
 \global\setbox\captionboxiii\box\Captionbox@
 }}
\long\def\islandtriplebox#1#2#3{\islandtripledata{#1}{#2}{#3}%
 \dimen@\ht\captionboxi
 \ifdim\ht\captionboxii>\dimen@ \dimen@\ht\captionboxii\fi
 \ifdim\ht\captionboxiii>\dimen@ \dimen@\ht\captionboxiii\fi
 \ifdim\dimen@>\z@
  \ifdim\ht\captionboxi<\dimen@
   \global\setbox\captionboxi\vbox to\dimen@{\unvbox\captionboxi\vfil}\fi
  \ifdim\ht\captionboxii<\dimen@
   \global\setbox\captionboxii\vbox to\dimen@{\unvbox\captionboxii\vfil}\fi
  \ifdim\ht\captionboxiii<\dimen@
   \global\setbox\captionboxiii\vbox to\dimen@{\unvbox\captionboxiii\vfil}\fi
 \fi
 \global\setbox\islandbox@\vbox
  {\hbox to\hsize{\hfil\box\islandboxi\hfil\box\islandboxii\hfil
   \box\islandboxiii\hfil}%
 \ifdim\dimen@>\z@\nointerlineskip
  {\rm\global\skipi@\dp\strutbox}\global\advance\skipi@\medskipamount
  \vskip\skipi@
  \hbox to\hsize{\hfil\box\captionboxi\hfil\box\captionboxii\hfil
   \box\captionboxiii\hfil}\fi}}
\def\islandtripleboxa#1#2#3{\islandtripledata{#1}{#2}{#3}%
 \dimen@\ht\captionboxi
 \ifdim\ht\captionboxii>\dimen@ \dimen@\ht\captionboxii\fi
 \ifdim\ht\captionboxiii>\dimen@ \dimen@\ht\captionboxiii\fi
 \ifdim\dimen@>\z@
  \ifdim\ht\captionboxi<\dimen@
   \global\setbox\captionboxi\vbox to\dimen@{\vfil\unvbox\captionboxi}\fi
  \ifdim\ht\captionboxii<\dimen@
   \global\setbox\captionboxii\vbox to\dimen@{\vfil\unvbox\captionboxii}\fi
  \ifdim\ht\captionboxiii<\dimen@
   \global\setbox\captionboxiii\vbox to\dimen@{\vfil\unvbox\captionboxiii}\fi
 \fi
 \dimen@ii\ht\islandboxi
 \ifdim\ht\islandboxii>\dimen@ii \dimen@ii\ht\islandboxii\fi
 \ifdim\ht\islandboxiii>\dimen@ii \dimen@ii\ht\islandboxiii\fi
 \ifdim\dimen@ii>\z@
  \ifdim\ht\islandboxi<\dimen@ii
   \global\setbox\islandboxi\vbox to\dimen@ii{\box\islandboxi\vfil}\fi
  \ifdim\ht\islandboxii<\dimen@ii
   \global\setbox\islandboxii\vbox to\dimen@ii{\box\islandboxii\vfil}\fi
  \ifdim\ht\islandboxiii<\dimen@ii
   \global\setbox\islandboxiii\vbox to\dimen@ii{\box\islandboxiii\vfil}\fi
 \fi
 \global\setbox\islandbox@\vbox
  {\ifdim\dimen@>\z@
  \hbox to\hsize{\hfil\box\captionboxi\hfil\box\captionboxii\hfil
   \box\captionboxiii\hfil}%
  \nointerlineskip{\rm\global\skipi@-\dp\strutbox}%
  \global\advance\skipi@\bigskipamount\vskip\skipi@\fi
  \hbox to\hsize{\hfil\box\islandboxi\hfil\box\islandboxii\hfil
   \box\islandboxiii\hfil}}}
\def\Figurepair#1\and#2\endFigurepair{\island@true
 \islandpairbox{\Figure#1\endFigure}{\Figure#2\endFigure}}
\def\Figuretriple#1\and#2\and#3\endFiguretriple{\island@true
 \islandtriplebox{\Figure#1\endFigure}{\Figure#2\endFigure}%
  {\Figure#3\endFigure}}
\def\Tablepair#1\and#2\endTablepair{\island@true
 \islandpairboxa{\Table#1\endTable}{\Table#2\endTable}}
\def\Tabletriple#1\and#2\and#3\endTabletriple{\island@true
 \islandtripleboxa{\Table#1\endTable}{\Table#2\endTable}%
 {\Table#3\endTable}}
\def\place#1{\place@true\island@false
 #1%
 \ifisland@
  \box\islandbox@
 \else
  \Err@{Whoa ... there's no \string\Figure, \string\Table,
   etc., here}%
 \fi
 \place@false}
\newskip\belowtopfigskip
\belowtopfigskip 15\p@ plus 5\p@ minus5\p@
\newskip\abovebotfigskip
\abovebotfigskip 18\p@ plus 6\p@ minus6\p@
\newdimen\minpagesize
\minpagesize 5pc
\dimen@\belowtopfigskip
\advance\dimen@-\abovebotfigskip
\skip\topins\dimen@
\dimen\topins\z@
\newcount\topinscount@
\newbox\topinsdims@
\def\storedim@{\global\setbox\topinsdims@
 \vbox{\hbox to\dimen@{}\unvbox\topinsdims@}}
\def\advancedimtopins@{%
 \ifnum\pageno=\@ne
 \else
   \advance\dimen@\dimen\topins
   \global\dimen\topins\dimen@
 \fi}
\newcount\flipcount@
\def\fliptopins@{%
 \global\flipcount@\z@
 \ifvoid\topins\else
 \setbox\z@\vbox
  {\vskip\p@
   \unvbox\topins
   \global\setbox\topins\vbox{}%
   \loop
    \test@false
    \ifdim\lastskip=\z@\unskip
     \ifdim\lastskip=\z@
      \test@true\fi\fi
    \iftest@
    \global\advance\flipcount@\@ne
    \setboxzl@
    \global\setbox\topins\vbox{\unvbox\topins\boxz@}%
    \unpenalty
   \repeat}\fi}
\newif\ifPar@
\newcount\Parcount@
\newbox\Parbox@
\expandafter\newbox\csname Parfigbox1\endcsname
\expandafter\newbox\csname Parfigbox2\endcsname
\expandafter\newbox\csname Parfigbox3\endcsname
\expandafter\newbox\csname Parfigbox4\endcsname
\expandafter\newbox\csname Parfigbox5\endcsname
\expandafter\newdimen\csname Parprev1\endcsname
\expandafter\newdimen\csname Parprev2\endcsname
\expandafter\newdimen\csname Parprev3\endcsname
\expandafter\newdimen\csname Parprev4\endcsname
\expandafter\newdimen\csname Parprev5\endcsname
\expandafter\newdimen\csname Parprev6\endcsname
\def\Par{\par\global\csname Parprev1\endcsname\prevdepth
 \global\Parcount@\@ne
 \global\Par@true\global\let\Parlist@\empty
 \global\setbox\Parbox@\vbox\bgroup\break}
\def\place@#1#2{%
 \ifisland@
  \ifhmode
   \ifPar@
    \ifnum\Parcount@>5
     \Err@{Only 5 \string\place's allowed per
      \string\Par...\noexpand\endPar paragraph}%
    \else
     \expandafter\expandafter\expandafter
      \global\expandafter\setbox
       \csname Parfigbox\number\Parcount@\endcsname\box\islandbox@
     \global\advance\Parcount@\@ne
     \xdef\Parlist@{\Parlist@#1}%
    \fi
   \else
    \vadjust{#2}%
   \fi
  \else
   #2%
  \fi
 \else
  \Err@{Whoa ... there's no \string\Figure,
   \string\Table, etc., here}%
 \fi
 \place@false}
\long\def\Aplace#1{\prevanish@
 \place@true\island@false
 #1%
 \place@ a\Aplace@
 \postvanish@}
\long\def\AAplace#1{\prevanish@\place@true\island@false
 #1%
 \place@ A\AAplace@
 \postvanish@}
\newif\ifAA@
\def\AAplace@{\AA@true\Aplace@\AA@false}
\let\AAlist@\empty
\def\Aplace@{\allowbreak
 \dimen@=\ht\islandbox@
 \advance\dimen@\abovebotfigskip
 \ht\islandbox@\dimen@
 \advance\dimen@\dp\islandbox@
 \storedim@
 \ifAA@
  \xdef\AAlist@{\AAlist@1}%
  \advancedimtopins@
 \else
  \xdef\AAlist@{\AAlist@0}%
  \ifnum\topinscount@>\@ne\else\advancedimtopins@\fi
 \fi
 \insert\topins{\penalty\z@\splittopskip\z@\floatingpenalty\z@
  \box\islandbox@}%
 \global\advance\topinscount@\@ne}
\long\def\Bplace#1{\prevanish@\place@true\island@false
 #1%
 \place@ b\Bplace@
 \postvanish@}
\def\Bplace@{\allowbreak
 \ifnum\topinscount@=\z@
  \setbox\z@\vbox{\vbox to-\belowtopfigskip{}}%
  \dimen@-\skip\topins
  \ht\z@\dimen@
  \storedim@
  \advancedimtopins@
  \insert\topins{\boxz@}%
  \global\advance\topinscount@\@ne
  \xdef\AAlist@{\AAlist@0}%
 \fi
 \dimen@\ht\islandbox@
 \advance\dimen@\abovebotfigskip
 \ht\islandbox@\dimen@
 \advance\dimen@\dp\islandbox@
 \storedim@
 \xdef\AAlist@{\AAlist@0}%
 \ifnum\topinscount@>\@ne\else\advancedimtopins@\fi
 \insert\topins{\penalty\z@\splittopskip\z@
  \floatingpenalty\z@
  \box\islandbox@}%
 \global\advance\topinscount@\@ne}
\def\breakisland@{\global\setbox\@ne\lastbox\global\skipi@\lastskip\unskip
 \global\setbox\thr@@\lastbox}%
\def\printisland@{\centerline{\box\thr@@}\nobreak\nointerlineskip
 \vskip\skipi@
 \ifdim\ht\@ne<\z@\box\@ne\else\centerline{\box\@ne}\fi}
\def\bottomfigs@{%
 \count@\@ne
 \loop
  \ifnum\count@<\flipcount@
  \nointerlineskip
  \vskip\abovebotfigskip
  \global\setbox\topins\vbox{\unvbox\topins\setboxzl@
   \unvbox\z@
   \breakisland@}%
  \printisland@
  \advance\count@\@ne
  \repeat}
\def\resetdimtopins@{%
 \global\advance\topinscount@-\flipcount@
 \global\setbox\topinsdims@\vbox
  {\unvbox\topinsdims@
   \count@\z@
   \DN@##1##2\next@{\gdef\AAlist@{##2}}%
   \loop
    \ifnum\count@<\flipcount@\setboxzl@
    \expandafter\next@\AAlist@\next@
    \advance\count@\@ne
    \repeat
   \dimen@\z@
   \count@\z@
   \setbox\tw@\vbox{}%
   \edef\nextiii@{\AAlist@}%
   \DN@##1##2\next@{\DNii@{##1}\def\nextiii@{##2}}%
   \loop
    \test@false
    \ifnum\count@<\topinscount@
    \expandafter\next@\nextiii@\next@
     \ifnum\count@<\tw@
      \test@true
     \else
      \if\nextii@ 1\test@true\fi
     \fi
    \fi
    \iftest@
     \setboxzl@
     \advance\dimen@\wdz@
     \setbox\tw@\vbox{\boxz@\unvbox\tw@}%
     \advance\count@\@ne
    \repeat
    \unvbox\tw@
    \global\dimen\topins\dimen@}}
\def\Place@#1#2{%
 \ifisland@
  \ifhmode
   \ifPar@
    \ifnum\Parcount@>5
     \Err@{Only 5 \string\place's allowed per
       \string\Par...\noexpand\endPar paragraph}%
    \else
     \expandafter\expandafter\expandafter\global\expandafter\setbox
      \csname Parfigbox\number\Parcount@\endcsname\box\islandbox@
     \global\advance\Parcount@\@ne
     \xdef\Parlist@{\Parlist@#1}%
     \vadjust{\break}%
    \fi
   \else
    \Err@{\noexpand#2allowed only in a \string\Par...\noexpand\endPar
     paragraph}%
   \fi
  \else
   #2%
  \fi
 \else
  \Err@{Who ... there's no \string\Figure, \string\Table,
   etc., here}%
 \fi
 \place@false}
\newif\ifC@
\newdimen\Cdim@
\long\def\Cplace#1{\prevanish@\place@true\island@false
 #1%
 \Place@ c\Cplace@
 \postvanish@}
\def\Cplace@{\allowbreak
 \ifnum\topinscount@>\z@\else
  \global\C@true\global\Cdim@\pagetotal\fi
 \Aplace@}
\long\def\Mplace#1{\prevanish@\place@true\island@false
 #1%
 \Place@ m\Mplace@
 \postvanish@}
\long\def\MXplace#1{\prevanish@\place@true\island@false
 #1%
 \Place@ M\MXplace@
 \postvanish@}
\newif\ifMX@
\def\MXplace@{\MX@true\Mplace@\MX@false}
\def\Mplace@{\allowbreak
 \dimen@\ht\islandbox@\advance\dimen@\dp\islandbox@
 \ifdim\pagetotal=\z@\else
  \ifdim\lastskip<\abovebotfigskip\advance\dimen@\abovebotfigskip
  \advance\dimen@-\lastskip\fi
 \fi
 \advance\dimen@\pagetotal
 \ifdim\dimen@>\pagegoal
  \Aplace@
 \else
  \nointerlineskip
  \ifdim\lastskip<\abovebotfigskip\removelastskip\vskip\abovebotfigskip\fi
  \setbox\z@\vbox{\unvbox\islandbox@
   \breakisland@}%
  \printisland@
  \ifnum\topinscount@=\z@
   \setbox\z@\vbox{\vbox to-\belowtopfigskip{}}%
   \dimen@-\skip\topins
   \ht\z@\dimen@
   \storedim@
   \advancedimtopins@
   \insert\topins{\boxz@}%
   \global\advance\topinscount@\@ne
   \xdef\AAlist@{\AAlist@0}%
  \fi
  \ifMX@
   \ifnum\topinscount@=\@ne
    \setbox\z@\vbox{\vbox to-\abovebotfigskip{}}%
    \ht\z@\z@
    \dimen@\z@
    \storedim@
    \advancedimtopins@
    \insert\topins{\boxz@}%
    \global\advance\topinscount@\@ne
    \xdef\AAlist@{\AAlist@0}%
   \fi
  \fi
  \nointerlineskip
  \vskip\belowtopfigskip
 \fi}
\expandafter\newbox\csname Parbox1\endcsname
\expandafter\newbox\csname Parbox2\endcsname
\expandafter\newbox\csname Parbox3\endcsname
\expandafter\newbox\csname Parbox4\endcsname
\expandafter\newbox\csname Parbox5\endcsname
\def\endPar{\egroup
 \count@\@ne
 {\vbadness\@M\vfuzz\maxdimen\splitmaxdepth\maxdimen\splittopskip\ht\strutbox
 \setbox\z@\vsplit\Parbox@ to\ht\Parbox@
 \loop
  \ifnum\count@<\Parcount@
  \expandafter\expandafter\expandafter\global\expandafter\setbox
   \csname Parbox\number\count@\endcsname\vsplit\Parbox@ to\ht\Parbox@
  \count@@\count@\advance\count@@\@ne
  \global\csname Parprev\number\count@@\endcsname
   \dp\csname Parbox\number\count@\endcsname
  \advance\count@\@ne
  \repeat}%
 \vskip\parskip
 \count@\@ne
 \def\nextv@##1##2\nextv@{\DN@{##1}\gdef\Parlist@{##2}}%
 \loop
  \ifnum\count@<\Parcount@
   \dimen@\csname Parprev\number\count@\endcsname
   \advance\dimen@\ht\strutbox
   \ifdim\dimen@<\baselineskip
    \advance\dimen@-\baselineskip\vskip-\dimen@
   \else
    \vskip\lineskip
   \fi
   \unvbox\csname Parbox\number\count@\endcsname
   \global\setbox\islandbox@\box\csname Parfigbox\number\count@\endcsname
   \expandafter\nextv@\Parlist@\nextv@
   \if a\next@\Aplace@\else
   \if A\next@\AAplace@\else
   \if b\next@\Bplace@\else
   \if c\next@\Cplace@\else
   \if m\next@\Mplace@\else
   \if M\next@\MXplace@\fi\fi\fi\fi\fi\fi
  \advance\count@\@ne
  \repeat
 \global\Par@false
 \ifvoid\Parbox@
  \prevdepth\csname Parprev\number\count@\endcsname
 \else
  \dimen@\csname Parprev\number\count@\endcsname\advance\dimen@\ht\strutbox
  \ifdim\dimen@<\baselineskip
    \advance\dimen@-\baselineskip\vskip-\dimen@
  \else
    \vskip\lineskip
  \fi
  \dimen@\dp\Parbox@
  \unvbox\Parbox@
  \prevdepth\dimen@
 \fi}
\def\folio{{\page@F\page@S{\page@P\page@N{\number\page@C}\page@Q}}}
\def\advancepageno{\global\advance\pageno\@ne}
\newif\ifspecialsplit@
\newbox\outbox@
\let\shipout@\shipout
\def\plainoutput{\specialsplit@false\ifvoid\topins\else\ifdim\ht\topins=\z@
 \specialsplit@true\advance\minpagesize-\skip\topins\fi\fi
 \fliptopins@
 \setbox\outbox@\vbox{\makeheadline\pagebody\makefootline}%
 {\noexpands@\let\style\relax
 \shipout@\box\outbox@}%
 \advancepageno
 \resetdimtopins@
 \ifvoid\@cclv\else\unvbox\@cclv\penalty\outputpenalty\fi
 \ifnum\outputpenalty>-\@MM\else\dosupereject\fi}
\def\pagebody{\vbox to\vsize{\boxmaxdepth\maxdepth
 \ifvoid\margin@\else
 \rlap{\kern\hsize\vbox to\z@{\kern4\p@\box\margin@\vss}}\fi
 \pagecontents}}
\newif\ifonlytop@
\def\pagecontents{%
 \onlytop@false
 \ifdim\ht\@cclv<\minpagesize\ifnum\flipcount@<\tw@\ifvoid\footins
  \onlytop@true\fi\fi\fi
 \test@false
 \ifC@
  \ifnum\flipcount@=\@ne
   \global\multiply\Cdim@\tw@
   \ifdim\Cdim@>\ht\@cclv
    \test@true
   \fi
  \fi
 \fi
 \global\C@false
 \iftest@
  \dimen@\ht\@cclv
  \advance\dimen@\skip\topins
  {\vfuzz\maxdimen\vbadness\@M
  \splitmaxdepth\maxdepth\splittopskip\topskip
  \setbox\z@\vsplit\@cclv to\dimen@
  \unvbox\z@}%
  \global\setbox\topins\vbox{\unvbox\topins
   \global\setbox\@ne\lastbox}%
  \setbox\z@\vbox{\unvbox\@ne
   \breakisland@}%
  \nointerlineskip
  \vskip\abovebotfigskip
  \printisland@
 \else
  \ifnum\flipcount@>\z@
   \global\setbox\topins\vbox{\unvbox\topins\global\setbox\@ne\lastbox}%
   \setbox\z@\vbox{\unvbox\@ne
    \breakisland@}%
   \printisland@
   \ifonlytop@\kern-\prevdepth\vfill\else\vskip\belowtopfigskip\fi
  \fi
 \fi
 \ifdim\ht\@cclv<\minpagesize
  \ifonlytop@\else\vfill\fi
 \else
  \ifspecialsplit@
   {\vfuzz\maxdimen\vbadness\@M
   \splitmaxdepth\maxdepth\splittopskip\topskip
   \dimen@ii\ht\@cclv \advance\dimen@ii\skip\topins
   \setbox\z@\vsplit\@cclv to\dimen@ii
   \unvbox\z@}%
  \else
   \unvbox\@cclv
  \fi
 \fi
 \bottomfigs@
 \ifvoid\footins\else\vskip\skip\footins\footnoterule\unvbox\footins\fi}
\newread\readdata@
\def\readthedata@#1{\expandafter
 \ifx\csname#1@D\endcsname\relax
  \immediate\openin\readdata@=#1.dat
  \ifeof\readdata@
   \Err@{No file #1.dat}%
  \else
   {\endlinechar\m@ne\gdef\Next@{}%
   \DNii@##1 ##2 ##3pt{\global\data@ht##1\global\data@dp##2%
    \global\data@wd##3pt}%
   \loop
    \ifeof\readdata@
    \else
    \read\readdata@ to\next@
    \ifx\next@\empty\else
     \edef\next@{\expandafter\nextii@\next@}%
     \expandafter\rightadd@\next@\to\Next@
    \fi
    \repeat}%
   \immediate\closein\readdata@
   \expandafter\expandafter\expandafter\global\expandafter
    \let\csname#1@D\endcsname\Next@\global\let\Next@\relax
  \fi
 \fi}
\newdimen\data@ht
\newdimen\data@dp
\newdimen\data@wd
\newif\ifgetdata@
\def\getdata@#1#2{\global\getdata@true\count@#2\relax
 {\let\\\or\xdef\Next@{\ifcase\number\count@#1\else
 \global\noexpand\getdata@false\fi}}\Next@}
\def\paste#1#2{\readthedata@{#1}%
 \getdata@{\csname#1@D\endcsname}{#2}%
 \ifgetdata@
 \dimen@\data@ht \advance\dimen@\data@dp
  \hbox{\special{dvipaste: #1 #2}%
   \lower\data@dp\vbox to\dimen@{\hbox to\data@wd{}\vfil}}%
 \else
  {\lccode`\Z=`\#\lccode`\N=`\N\lccode`\F=`\F%
   \lowercase{\Err@{No data for File [#1], Z#2}}}%
 \fi}
\newdimen\httable
\newdimen\dptable
\newdimen\wdtable
\def\measuretable#1#2{\readthedata@{#1}%
 \getdata@{\csname#1@D\endcsname}{#2}%
 \ifgetdata@
  \httable\data@ht \dptable\data@dp \wdtable\data@wd
 \else
  {\lccode`\Z=`\#\lccode`\N=`\N\lccode`\F=`\F%
  \lowercase{\Err@{No data for File [#1], Z#2}}}%
 \fi}
\def\East#1#2{\setboxz@h{$\m@th\ssize\;{#1}\;\;$}%
 \setbox\tw@\hbox{$\m@th\ssize\;{#2}\;\;$}\setbox4=\hbox{$\m@th#2$}%
 \dimen@\minaw@
 \ifdim\wdz@>\dimen@\dimen@\wdz@\fi\ifdim\wd\tw@>\dimen@\dimen@\wd\tw@\fi
 \ifdim\wd4 >\z@
  \mathrel{\mathop{\hbox to\dimen@{\rightarrowfill}}\limits^{#1}_{#2}}%
 \else
  \mathrel{\mathop{\hbox to\dimen@{\rightarrowfill}}\limits^{#1}}%
 \fi}
\def\West#1#2{\setboxz@h{$\m@th\ssize\;\;{#1}\;$}%
 \setbox\tw@\hbox{$\m@th\ssize\;\;{#2}\;$}\setbox4=\hbox{$\m@th#2$}%
 \dimen@\minaw@
 \ifdim\wdz@>\dimen@\dimen@\wdz@\fi\ifdim\wd\tw@>\dimen@\dimen@\wd\tw@\fi
 \ifdim\wd4 >\z@
  \mathrel{\mathop{\hbox to\dimen@{\leftarrowfill}}\limits^{#1}_{#2}}%
 \else
  \mathrel{\mathop{\hbox to\dimen@{\leftarrowfill}}\limits^{#1}}%
 \fi}
\font\arrow@i=lams1
\font\arrow@ii=lams2
\font\arrow@iii=lams3
\font\arrow@iv=lams4
\font\arrow@v=lams5
\newdimen\standardcgap
\standardcgap40\p@
\newdimen\hunit
\hunit\tw@\p@
\newdimen\standardrgap
\standardrgap32\p@
\newdimen\vunit
\vunit1.6\p@
\def\Cgaps#1{\RIfM@
 \standardcgap#1\standardcgap\relax\hunit#1\hunit\relax
 \else\nonmatherr@\Cgaps\fi}
\def\Rgaps#1{\RIfM@
 \standardrgap#1\standardrgap\relax\vunit#1\vunit\relax
 \else\nonmatherr@\Rgaps\fi}
\newdimen\getdim@
\def\getcgap@#1{\ifcase#1\or\getdim@\z@\else\getdim@\standardcgap\fi}
\def\getrgap@#1{\ifcase#1\getdim@\z@\else\getdim@\standardrgap\fi}
\def\cgaps{\RIfM@\expandafter\cgaps@\else\expandafter\nonmatherr@
 \expandafter\cgaps\fi}
\def\cgaps@{\ifnum\catcode`\;=\active\expandafter\cgapsA@\else
 \expandafter\cgapsO@\fi}
\def\cgapsO@#1{\toks@{\ifcase\i@\or\getdim@=\z@}%
 \gaps@@\standardcgap#1;\gaps@@\gaps@@
 \edef\next@{\the\toks@\noexpand\else\noexpand\getdim@\noexpand\standardcgap
  \noexpand\fi}%
 \toks@=\expandafter{\next@}%
 \edef\getcgap@##1{\i@##1\relax\the\toks@}\toks@{}}
{\catcode`\;=\active
 \gdef\cgapsA@#1{\toks@{\ifcase\i@\or\getdim@=\z@}%
 \gaps@@\standardcgap#1;\gaps@@\gaps@@
 \edef\next@{\the\toks@\noexpand\else\noexpand\getdim@\noexpand\standardcgap
  \noexpand\fi}%
 \toks@=\expandafter{\next@}%
 \edef\getcgap@##1{\i@##1\relax\the\toks@}\toks@{}}
}
\def\Gaps@@{\gaps@@}
\def\gaps@@#1#2;#3{\mgaps@#1#2\mgaps@
 \edef\next@{\the\toks@\noexpand\or\noexpand\getdim@
  \noexpand#1\the\mgapstoks@@}%
 \toks@\expandafter{\next@}%
 \DN@{#3}%
 \ifx\next@\Gaps@@\def\next@##1\gaps@@{}\else
  \def\next@{\gaps@@#1#3}\fi\next@}
{\catcode`\;=\active
 \gdef\rgaps#1{\RIfM@{\ifnum\catcode`\;=\active\def;{\string;}\fi
   \xdef\Next@{\noexpand\rgaps@{#1}}}%
  \Next@\edef\getrgap@##1{\i@##1\relax\the\toks@}\toks@{}\else
  \nonmatherr@\rgaps\fi}
}
\def\rgaps@#1{\toks@{\ifcase\i@\getdim@=\z@}%
 \gaps@@\standardrgap#1;\gaps@@\gaps@@
 \edef\next@{\the\toks@\noexpand\else\noexpand\getdim@\noexpand\standardrgap
  \noexpand\fi}%
 \toks@=\expandafter{\next@}}
\newbox\ZER@
\def\mgaps@#1{\let\mgapsnext@#1\FNSS@\mgaps@@}
\def\mgaps@@{\ifx\next\w\expandafter\mgaps@@@\else
 \expandafter\mgaps@@@@\fi}
\newtoks\mgapstoks@@
\def\mgaps@@@@#1\mgaps@{\getdim@\mgapsnext@\getdim@#1\getdim@
 \edef\next@{\noexpand\getdim@\the\getdim@}%
 \mgapstoks@@\expandafter{\next@}}
\def\mgaps@@@\w#1#2\mgaps@{\mgaps@@@@#2\mgaps@
 \setbox\ZER@\hbox{$\m@th\hskip15\p@\tsize@#1$}%
 \dimen@\wd\ZER@
 \ifdim\dimen@>\getdim@\getdim@\dimen@\fi
 \edef\next@{\noexpand\getdim@\the\getdim@}%
 \mgapstoks@@\expandafter{\next@}}
\def\changewidth#1#2{\setbox\ZER@{$\m@th#2}%
 \hbox to\wd\ZER@{\hss$\m@th#1$\hss}}
\atdef@({\FN@\ARROW@}
\def\ARROW@{\ifx\next)\let\next@\OPTIONS@\else
 \DN@{\csname\string @(\endcsname}\fi\next@}
\newif\ifoptions@
\def\OPTIONS@){\ifoptions@\let\next@\relax\else
 \DN@{\global\options@true\begingroup\optioncodes@}\fi\next@}
\newif\ifN@
\newif\ifE@
\newif\ifNESW@
\newif\ifH@
\newif\ifV@
\newif\ifHshort@
\expandafter\def\csname\string @(\endcsname #1,#2){%
 \ifoptions@\expandafter\endgroup\fi
 \N@false\E@false\H@false\V@false\Hshort@false
 \ifnum#1>\z@\E@true\fi
 \ifnum#1=\z@\V@true\global\tX@false\global\tY@false\global\a@false\fi
 \ifnum#2>\z@\N@true\fi
 \ifnum#2=\z@\H@true\global\tX@false\global\tY@false\global\a@false
  \ifshort@\Hshort@true\fi\fi
 \NESW@false
 \ifN@\ifE@\NESW@true\fi\else\ifE@\else\NESW@true\fi\fi
 \arrow@{#1}{#2}%
 \global\options@false
 \global\scount@\z@\global\tcount@\z@\global\arrcount@\z@
 \global\s@false\global\sxdimen@\z@\global\sydimen@\z@
 \global\tX@false\global\tXdimen@i\z@\global\tXdimen@ii\z@
 \global\tY@false\global\tYdimen@i\z@\global\tYdimen@ii\z@
 \global\a@false\global\exacount@\z@
 \global\x@false\global\xdimen@\z@
 \global\X@false\global\Xdimen@\z@
 \global\y@false\global\ydimen@\z@
 \global\Y@false\global\Ydimen@\z@
 \global\p@false\global\pdimen@\z@
 \global\label@ifalse\global\label@iifalse
 \global\dl@ifalse\global\ldimen@i\z@
 \global\dl@iifalse\global\ldimen@ii\z@
 \global\short@false\global\unshort@false}
\newif\iflabel@i
\newif\iflabel@ii
\newcount\scount@
\newcount\tcount@
\newcount\arrcount@
\newif\ifs@
\newdimen\sxdimen@
\newdimen\sydimen@
\newif\iftX@
\newdimen\tXdimen@i
\newdimen\tXdimen@ii
\newif\iftY@
\newdimen\tYdimen@i
\newdimen\tYdimen@ii
\newif\ifa@
\newcount\exacount@
\newif\ifx@
\newdimen\xdimen@
\newif\ifX@
\newdimen\Xdimen@
\newif\ify@
\newdimen\ydimen@
\newif\ifY@
\newdimen\Ydimen@
\newif\ifp@
\newdimen\pdimen@
\newif\ifdl@i
\newif\ifdl@ii
\newdimen\ldimen@i
\newdimen\ldimen@ii
\newif\ifshort@
\newif\ifunshort@
\def\zero@#1{\ifnum\scount@=\z@
 \if#1e\global\scount@\m@ne\else
 \if#1t\global\scount@\tw@\else
 \if#1h\global\scount@\thr@@\else
 \if#1'\global\scount@6 \else
 \if#1`\global\scount@7 \else
 \if#1(\global\scount@8 \else
 \if#1)\global\scount@9 \else
 \if#1s\global\scount@12 \else
 \if#1H\global\scount@13 \else
 \Err@{\Invalid@@ option \string\0}\fi\fi\fi\fi\fi\fi\fi\fi\fi
 \fi}
\def\one@#1{\ifnum\tcount@=\z@
 \if#1e\global\tcount@\m@ne\else
 \if#1h\global\tcount@\tw@\else
 \if#1t\global\tcount@\thr@@\else
 \if#1'\global\tcount@4 \else
 \if#1`\global\tcount@5 \else
 \if#1(\global\tcount@\ten@ \else
 \if#1)\global\tcount@11 \else
 \if#1s\global\tcount@12 \else
 \if#1H\global\tcount@13 \else
 \Err@{\Invalid@@ option \string\1}\fi\fi\fi\fi\fi\fi\fi\fi\fi
 \fi}
\def\a@#1{\ifnum\arrcount@=\z@
 \if#10\global\arrcount@\m@ne\else
 \if#1+\global\arrcount@\@ne\else
 \if#1-\global\arrcount@\tw@\else
 \if#1=\global\arrcount@\thr@@\else
 \Err@{\Invalid@@ option \string\a}\fi\fi\fi\fi
 \fi}
\def\ds@{\ifnum\catcode`\;=\active\expandafter\dsA@\else
 \expandafter\dsO@\fi}
\def\dsO@(#1;#2){\ds@@{#1}{#2}}
\def\ds@@#1#2{\ifs@\else
 \global\s@true
 \global\sxdimen@\hunit\global\sxdimen@#1\sxdimen@\relax
 \global\sydimen@\vunit\global\sydimen@#2\sydimen@\relax
 \fi}
\def\dtX@{\ifnum\catcode`\;=\active\expandafter\dtXA@\else
 \expandafter\dtXO@\fi}
\def\dtXO@(#1;#2){\dtX@@{#1}{#2}}
\def\dtX@@#1#2{\iftX@\else
 \global\tX@true
 \global\tXdimen@i\hunit\global\tXdimen@i#1\tXdimen@i\relax
 \global\tXdimen@ii\vunit\global\tXdimen@ii#2\tXdimen@ii\relax
 \fi}
\def\dtY@{\ifnum\catcode`\;=\active\expandafter\dtYA@\else
 \expandafter\dtYO@\fi}
\def\dtYO@(#1;#2){\dtY@@{#1}{#2}}
\def\dtY@@#1#2{\iftY@\else
 \global\tY@true
 \global\tYdimen@i\hunit\global\tYdimen@i#1\tYdimen@i\relax
 \global\tYdimen@ii\vunit\global\tYdimen@ii#2\tYdimen@ii\relax
 \fi}
{\catcode`\;=\active
 \gdef\dsA@(#1;#2){\ds@@{#1}{#2}}
 \gdef\dtXA@(#1;#2){\dtX@@{#1}{#2}}
 \gdef\dtYA@(#1;#2){\dtY@@{#1}{#2}}
}
\def\da@#1{\ifa@\else\global\a@true\global\exacount@#1\relax\fi}
\def\dx@#1{\ifx@\else
 \global\x@true
 \global\xdimen@\hunit\global\xdimen@#1\xdimen@\relax
 \fi}
\def\dX@#1{\ifX@\else
 \global\X@true
 \global\Xdimen@\hunit\global\Xdimen@#1\Xdimen@\relax
 \fi}
\def\dy@#1{\ify@\else
 \global\y@true
 \global\ydimen@\vunit\global\ydimen@#1\ydimen@\relax
 \fi}
\def\dY@#1{\ifY@\else
 \global\Y@true
 \global\Ydimen@\vunit\global\Ydimen@#1\Ydimen@\relax
 \fi}
\def\p@@#1{\ifp@\else
 \global\p@true
 \global\pdimen@\hunit\global\divide\pdimen@\tw@
 \global\pdimen@#1\pdimen@\relax
 \fi}
\def\L@#1{\iflabel@i\else
 \global\label@itrue\gdef\label@i{#1}%
 \fi}
\def\l@#1{\iflabel@ii\else
 \global\label@iitrue\gdef\label@ii{#1}%
 \fi}
\def\dL@#1{\ifdl@i\else
 \global\dl@itrue\global\ldimen@i\hunit\global\ldimen@i#1\ldimen@i\relax
 \fi}
\def\dl@#1{\ifdl@ii\else
 \global\dl@iitrue\global\ldimen@ii\hunit\global\ldimen@ii#1\ldimen@ii\relax
 \fi}
\def\s@{\ifunshort@\else\global\short@true\fi}
\def\uns@{\ifshort@\else\global\unshort@true\global\short@false\fi}
\def\optioncodes@{\let\0\zero@\let\1\one@\let\a\a@\let\ds\ds@\let\dtX\dtX@
 \let\dtY\dtY@\let\da\da@\let\dx\dx@\let\dX\dX@\let\dY\dY@\let\dy\dy@
 \let\p\p@@\let\L\L@\let\l\l@\let\dL\dL@\let\dl\dl@\let\s\s@\let\uns\uns@}
\def\slopes@{\\161\\152\\143\\134\\255\\126\\357\\238\\349\\45{10}\\56{11}%
 \\11{12}\\65{13}\\54{14}\\43{15}\\32{16}\\53{17}\\21{18}\\52{19}\\31{20}%
 \\41{21}\\51{22}\\61{23}}
\newcount\tan@i
\newcount\tan@ip
\newcount\tan@ii
\newcount\tan@iip
\newdimen\slope@i
\newdimen\slope@ip
\newdimen\slope@ii
\newdimen\slope@iip
\newcount\angcount@
\newcount\extracount@
\def\slope@{{\slope@i\secondy@\advance\slope@i-\firsty@
 \ifN@\else\multiply\slope@i\m@ne\fi
 \slope@ii\secondx@\advance\slope@ii-\firstx@
 \ifE@\else\multiply\slope@ii\m@ne\fi
 \ifdim\slope@ii<\z@
  \global\tan@i6 \global\tan@ii\@ne\global\angcount@23
 \else
  \dimen@\slope@i\multiply\dimen@6
  \ifdim\dimen@<\slope@ii
   \global\tan@i\@ne\global\tan@ii6 \global\angcount@\@ne
  \else
   \dimen@\slope@ii\multiply\dimen@6
   \ifdim\dimen@<\slope@i
    \global\tan@i6 \global\tan@ii\@ne\global\angcount@23
   \else
    \global\tan@ip\z@\global\tan@iip\@ne
    \def\\##1##2##3{\global\angcount@##3\relax
     \slope@ip\slope@i\slope@iip\slope@ii
     \multiply\slope@iip##1\relax\multiply\slope@ip##2\relax
     \ifdim\slope@iip<\slope@ip
      \global\tan@ip##1\relax\global\tan@iip##2\relax
     \else
      \global\tan@i##1\relax\global\tan@ii##2\relax
      \def\\####1####2####3{}%
     \fi}%
    \slopes@
    \slope@i\secondy@\advance\slope@i-\firsty@
    \ifN@\else\multiply\slope@i\m@ne\fi
    \multiply\slope@i\tan@ii\multiply\slope@i\tan@iip\multiply\slope@i\tw@
    \count@\tan@i\multiply\count@\tan@iip
    \extracount@\tan@ip\multiply\extracount@\tan@ii
    \advance\count@\extracount@
    \slope@ii\secondx@\advance\slope@ii-\firstx@
    \ifE@\else\multiply\slope@ii\m@ne\fi
    \multiply\slope@ii\count@
    \ifdim\slope@i<\slope@ii
     \global\tan@i\tan@ip\global\tan@ii\tan@iip
     \global\advance\angcount@\m@ne
    \fi
   \fi
  \fi
 \fi}%
}
\def\slope@a#1{{\def\\##1##2##3{\ifnum##3=#1\global\tan@i##1\relax
 \global\tan@ii##2\relax\fi}\slopes@}}
\newcount\i@
\newcount\j@
\newcount\colcount@
\newcount\Colcount@
\newcount\tcolcount@
\newdimen\rowht@
\newdimen\rowdp@
\newcount\rowcount@
\newcount\Rowcount@
\newcount\maxcolrow@
\newtoks\colwidthtoks@
\newtoks\Rowheighttoks@
\newtoks\Rowdepthtoks@
\newtoks\widthtoks@
\newtoks\Widthtoks@
\newtoks\heighttoks@
\newtoks\Heighttoks@
\newtoks\depthtoks@
\newtoks\Depthtoks@
\newif\iffirstCDcr@
\def\dotoks@i{%
 \global\widthtoks@\expandafter{\the\widthtoks@\else\getdim@\z@\fi}%
 \global\heighttoks@\expandafter{\the\heighttoks@\else\getdim@\z@\fi}%
 \global\depthtoks@\expandafter{\the\depthtoks@\else\getdim@\z@\fi}}
\def\dotoks@ii{%
 \global\widthtoks@{\ifcase\j@}%
 \global\heighttoks@{\ifcase\j@}%
 \global\depthtoks@{\ifcase\j@}}
\def\preCD@#1\endCD{\setbox\ZER@
 \vbox{%
  \def\arrow@##1##2{{}}%
  \global\rowcount@\m@ne\global\colcount@\z@\global\Colcount@\z@
  \global\firstCDcr@true\toks@{}%
  \global\widthtoks@{\ifcase\j@}%
  \global\Widthtoks@{\ifcase\i@}%
  \global\heighttoks@{\ifcase\j@}%
  \global\Heighttoks@{\ifcase\i@}%
  \global\depthtoks@{\ifcase\j@}%
  \global\Depthtoks@{\ifcase\i@}%
  \global\Rowheighttoks@{\ifcase\i@}%
  \global\Rowdepthtoks@{\ifcase\i@}%
  \Let@
  \everycr{%
   \noalign{%
    \global\advance\rowcount@\@ne
    \ifnum\colcount@<\Colcount@
    \else
     \global\Colcount@\colcount@\global\maxcolrow@\rowcount@
    \fi
    \global\colcount@\z@
    \iffirstCDcr@
     \global\firstCDcr@false
    \else
     \edef\next@{\the\Rowheighttoks@\noexpand\or\noexpand\getdim@\the\rowht@}%
      \global\Rowheighttoks@\expandafter{\next@}%
     \edef\next@{\the\Rowdepthtoks@\noexpand\or\noexpand\getdim@\the\rowdp@}%
      \global\Rowdepthtoks@\expandafter{\next@}%
     \global\rowht@\z@\global\rowdp@\z@
     \dotoks@i
     \edef\next@{\the\Widthtoks@\noexpand\or\the\widthtoks@}%
      \global\Widthtoks@\expandafter{\next@}%
     \edef\next@{\the\Heighttoks@\noexpand\or\the\heighttoks@}%
      \global\Heighttoks@\expandafter{\next@}%
     \edef\next@{\the\Depthtoks@\noexpand\or\the\depthtoks@}%
      \global\Depthtoks@\expandafter{\next@}%
     \dotoks@ii
    \fi}}%
  \tabskip\z@
  \halign{&\setbox\ZER@\hbox{\vrule\height\ten@\p@\width\z@\depth\z@     
   $\m@th\displaystyle{##}$}\copy\ZER@
   \ifdim\ht\ZER@>\rowht@\global\rowht@\ht\ZER@\fi
   \ifdim\dp\ZER@>\rowdp@\global\rowdp@\dp\ZER@\fi
   \global\advance\colcount@\@ne
   \edef\next@{\the\widthtoks@\noexpand\or\noexpand\getdim@\the\wd\ZER@}%
    \global\widthtoks@\expandafter{\next@}%
   \edef\next@{\the\heighttoks@\noexpand\or\noexpand\getdim@\the\ht\ZER@}%
    \global\heighttoks@\expandafter{\next@}%
   \edef\next@{\the\depthtoks@\noexpand\or\noexpand\getdim@\the\dp\ZER@}%
    \global\depthtoks@\expandafter{\next@}%
   \cr#1\crcr}}%
 \Rowcount@\rowcount@
 \global\Widthtoks@\expandafter{\the\Widthtoks@\fi\relax}%
 \edef\Width@##1##2{\i@##1\relax\j@##2\relax\the\Widthtoks@}%
 \global\Heighttoks@\expandafter{\the\Heighttoks@\fi\relax}%
 \edef\Height@##1##2{\i@##1\relax\j@##2\relax\the\Heighttoks@}%
 \global\Depthtoks@\expandafter{\the\Depthtoks@\fi\relax}%
 \edef\Depth@##1##2{\i@##1\relax\j@##2\relax\the\Depthtoks@}%
 \edef\next@{\the\Rowheighttoks@\noexpand\fi\relax}%
 \global\Rowheighttoks@\expandafter{\next@}%
 \edef\Rowheight@##1{\i@##1\relax\the\Rowheighttoks@}%
 \edef\next@{\the\Rowdepthtoks@\noexpand\fi\relax}%
 \global\Rowdepthtoks@\expandafter{\next@}%
 \edef\Rowdepth@##1{\i@##1\relax\the\Rowdepthtoks@}%
 \global\colwidthtoks@{\fi}%
 \setbox\ZER@\vbox{%
  \unvbox\ZER@
  \count@\rowcount@
  \loop
   \unskip\unpenalty
   \setbox\ZER@\lastbox
   \ifnum\count@>\maxcolrow@\advance\count@\m@ne
   \repeat
  \hbox{%
   \unhbox\ZER@
   \count@\z@
   \loop
    \unskip
    \setbox\ZER@\lastbox
    \edef\next@{\noexpand\or\noexpand\getdim@\the\wd\ZER@\the\colwidthtoks@}%
     \global\colwidthtoks@\expandafter{\next@}%
    \advance\count@\@ne
    \ifnum\count@<\Colcount@
    \repeat}}%
 \edef\next@{\noexpand\ifcase\noexpand\i@\the\colwidthtoks@}%
  \global\colwidthtoks@\expandafter{\next@}%
 \edef\Colwidth@##1{\i@##1\relax\the\colwidthtoks@}%
 \global\colwidthtoks@{}\global\Rowheighttoks@{}\global\Rowdepthtoks@{}%
 \global\widthtoks@{}\global\Widthtoks@{}\global\heighttoks@{}%
 \global\Heighttoks@{}\global\depthtoks@{}\global\Depthtoks@{}%
}
\newcount\xoff@
\newcount\yoff@
\newcount\endcount@
\newcount\rcount@
\newdimen\firstx@
\newdimen\firsty@
\newdimen\secondx@
\newdimen\secondy@
\newdimen\tocenter@
\newdimen\charht@
\newdimen\charwd@
\def\outside@{\Err@{This arrow points outside the \string\CD}}
\newif\ifsvertex@
\newif\iftvertex@
\def\arrow@#1#2{\global\xoff@#1\relax\global\yoff@#2\relax
 \count@\rowcount@\advance\count@-\yoff@
 \ifnum\count@<\@ne\outside@\else\ifnum\count@>\Rowcount@\outside@\fi\fi
 \count@\colcount@\advance\count@\xoff@
 \ifnum\count@<\@ne\outside@\else\ifnum\count@>\Colcount@\outside@\fi\fi
 \tcolcount@\colcount@\advance\tcolcount@\xoff@
 \Width@\rowcount@\colcount@\divide\getdim@\tw@\tocenter@-\getdim@
 \ifdim\getdim@=\z@
  \firstx@\z@\firsty@\mathaxis@\svertex@true
 \else
  \svertex@false
  \ifHshort@
   \Colwidth@\colcount@\divide\getdim@\tw@
   \ifE@ \firstx@\getdim@ \else \firstx@-\getdim@ \fi
  \else
   \ifE@ \firstx@\getdim@ \else \firstx@-\getdim@ \fi
  \fi
  \ifE@
   \ifH@ \advance\firstx@\thr@@\p@ \else \advance\firstx@-\thr@@\p@ \fi  
  \else
   \ifH@ \advance\firstx@-\thr@@\p@ \else \advance\firstx@\thr@@\p@ \fi  
  \fi
  \ifN@
   \Height@\rowcount@\colcount@ \firsty@\getdim@                         
   \ifV@ \advance\firsty@\thr@@\p@ \fi                                   
  \else
   \ifV@
    \Depth@\rowcount@\colcount@ \firsty@-\getdim@                        
    \advance\firsty@-\thr@@\p@                                           
   \else
    \firsty@\z@                                                          
   \fi
  \fi
 \fi
 \ifV@
 \else
  \Colwidth@\colcount@\divide\getdim@\tw@
  \ifE@\secondx@\getdim@\else\secondx@-\getdim@\fi
  \ifE@\else\getcgap@\colcount@\advance\secondx@-\getdim@\fi
  \endcount@\colcount@\advance\endcount@\xoff@
  \count@\colcount@
  \ifE@
   \advance\count@\@ne
   \loop
    \ifnum\count@<\endcount@
    \Colwidth@\count@\advance\secondx@\getdim@
    \getcgap@\count@\advance\secondx@\getdim@
    \advance\count@\@ne
    \repeat
  \else
   \advance\count@\m@ne
   \loop
    \ifnum\count@>\endcount@
    \Colwidth@\count@\advance\secondx@-\getdim@
    \getcgap@\count@\advance\secondx@-\getdim@
    \advance\count@\m@ne
    \repeat
  \fi
  \Colwidth@\count@\divide\getdim@\tw@
  \ifHshort@
  \else
   \ifE@\advance\secondx@\getdim@\else\advance\secondx@-\getdim@\fi
  \fi
  \ifE@\getcgap@\count@\advance\secondx@\getdim@\fi
  \rcount@\rowcount@\advance\rcount@-\yoff@
  \Width@\rcount@\count@\divide\getdim@\tw@
  \tvertex@false
  \ifH@\ifdim\getdim@=\z@\tvertex@true\Hshort@false\fi\fi
  \ifHshort@
  \else
   \ifE@\advance\secondx@-\getdim@\else\advance\secondx@\getdim@\fi
  \fi
  \iftvertex@
   \advance\secondx@.4\p@
  \else
   \ifE@\advance\secondx@-\thr@@\p@\else\advance\secondx@\thr@@\p@\fi    
  \fi
 \fi
 \ifH@
 \else
  \ifN@
   \Rowheight@\rowcount@\secondy@\getdim@
  \else
   \Rowdepth@\rowcount@\secondy@-\getdim@
   \getrgap@\rowcount@\advance\secondy@-\getdim@
  \fi
  \endcount@\rowcount@\advance\endcount@-\yoff@
  \count@\rowcount@
  \ifN@
   \advance\count@\m@ne
   \loop
    \ifnum\count@>\endcount@
    \Rowheight@\count@\advance\secondy@\getdim@
    \Rowdepth@\count@\advance\secondy@\getdim@
    \getrgap@\count@\advance\secondy@\getdim@
    \advance\count@\m@ne
    \repeat
  \else
   \advance\count@\@ne
   \loop
    \ifnum\count@<\endcount@
    \Rowheight@\count@\advance\secondy@-\getdim@
    \Rowdepth@\count@\advance\secondy@-\getdim@
    \getrgap@\count@\advance\secondy@-\getdim@
    \advance\count@\@ne
    \repeat
  \fi
  \tvertex@false
  \ifV@\Width@\count@\colcount@\ifdim\getdim@=\z@\tvertex@true\fi\fi
  \ifN@
   \getrgap@\count@\advance\secondy@\getdim@
   \Rowdepth@\count@\advance\secondy@\getdim@
   \iftvertex@
    \advance\secondy@\mathaxis@
   \else
    \Depth@\count@\tcolcount@\advance\secondy@-\getdim@
    \advance\secondy@-\thr@@\p@                                          
   \fi
  \else
   \Rowheight@\count@\advance\secondy@-\getdim@
   \iftvertex@
    \advance\secondy@\mathaxis@
   \else
    \Height@\count@\tcolcount@\advance\secondy@\getdim@
    \advance\secondy@\thr@@\p@                                           
   \fi
  \fi
 \fi
 \ifV@\else\advance\firstx@\sxdimen@\fi
 \ifH@\else\advance\firsty@\sydimen@\fi
 \iftX@
  \advance\secondy@\tXdimen@ii
  \advance\secondx@\tXdimen@i
  \slope@
 \else
  \iftY@
   \advance\secondy@\tYdimen@ii
   \advance\secondx@\tYdimen@i
   \slope@
   \secondy@\secondx@\advance\secondy@-\firstx@
   \ifNESW@\else\multiply\secondy@\m@ne\fi
   \multiply\secondy@\tan@i\divide\secondy@\tan@ii\advance\secondy@\firsty@
  \else
   \ifa@
    \slope@
    \ifNESW@\global\advance\angcount@\exacount@\else
     \global\advance\angcount@-\exacount@\fi
    \ifnum\angcount@>23 \global\angcount@23 \fi
    \ifnum\angcount@<\@ne\global\angcount@\@ne\fi
    \slope@a\angcount@
    \ifY@
     \advance\secondy@\Ydimen@
    \else
     \ifX@
      \advance\secondx@\Xdimen@
      \dimen@\secondx@\advance\dimen@-\firstx@
      \ifNESW@\else\multiply\dimen@\m@ne\fi
      \multiply\dimen@\tan@i\divide\dimen@\tan@ii
      \advance\dimen@\firsty@\secondy@\dimen@
     \fi
    \fi
   \else
    \ifH@\else\ifV@\else\slope@\fi\fi
   \fi
  \fi
 \fi
 \ifH@\else\ifV@\else\ifsvertex@\else
  \dimen@6\p@\multiply\dimen@\tan@ii
  \count@\tan@i\advance\count@\tan@ii\divide\dimen@\count@
  \ifE@\advance\firstx@\dimen@\else\advance\firstx@-\dimen@\fi
  \multiply\dimen@\tan@i\divide\dimen@\tan@ii
  \ifN@\advance\firsty@\dimen@\else\advance\firsty@-\dimen@\fi
 \fi\fi\fi
 \ifp@
  \ifH@\else\ifV@\else
   \getcos@\pdimen@\advance\firsty@\dimen@\advance\secondy@\dimen@
   \ifNESW@\advance\firstx@-\dimen@ii\else\advance\firstx@\dimen@ii\fi
  \fi\fi
 \fi
 \ifH@\else\ifV@\else
  \ifnum\tan@i>\tan@ii
   \charht@\ten@\p@\charwd@\ten@\p@
   \multiply\charwd@\tan@ii\divide\charwd@\tan@i
  \else
   \charwd@\ten@\p@\charht@\ten@\p@
   \divide\charht@\tan@ii\multiply\charht@\tan@i
  \fi
  \ifnum\tcount@=\thr@@
   \ifN@\advance\secondy@-.3\charht@\else\advance\secondy@.3\charht@\fi
  \fi
  \ifnum\scount@=\tw@
   \ifE@\advance\firstx@.3\charht@\else\advance\firstx@-.3\charht@\fi
  \fi
  \ifnum\tcount@=12
   \ifN@\advance\secondy@-\charht@\else\advance\secondy@\charht@\fi
  \fi
  \iftY@
  \else
   \ifa@
    \ifX@
    \else
     \secondx@\secondy@\advance\secondx@-\firsty@
     \ifNESW@\else\multiply\secondx@\m@ne\fi
     \multiply\secondx@\tan@ii\divide\secondx@\tan@i
     \advance\secondx@\firstx@
    \fi
   \fi
  \fi
 \fi\fi
 \ifH@\harrow@\else\ifV@\varrow@\else\arrow@@\fi\fi}
\newdimen\mathaxis@
\mathaxis@90\p@\divide\mathaxis@36
\def\harrow@b{\ifE@\hskip\tocenter@\hskip\firstx@\fi}
\def\harrow@bb{\ifE@\hskip\xdimen@\else\hskip\Xdimen@\fi}
\def\harrow@e{\ifE@\else\hskip-\firstx@\hskip-\tocenter@\fi}
\def\harrow@ee{\ifE@\hskip-\Xdimen@\else\hskip-\xdimen@\fi}
\def\harrow@{\dimen@\secondx@\advance\dimen@-\firstx@
 \ifE@\let\next@\rlap\else\multiply\dimen@\m@ne\let\next@\llap\fi
 \next@{%
  \harrow@b
  \smash{\raise\pdimen@\hbox to\dimen@
   {\harrow@bb\arrow@ii
    \ifnum\arrcount@=\m@ne\else\ifnum\arrcount@=\thr@@\else
     \ifE@
      \ifnum\scount@=\m@ne
      \else
       \ifcase\scount@\or\or\char118 \or\char117 \or\or\or\char119 \or
       \char120 \or\char121 \or\char122 \or\or\or\arrow@i\char125 \or
       \char117 \hskip\thr@@\p@\char117 \hskip-\thr@@\p@\fi
      \fi
     \else
      \ifnum\tcount@=\m@ne
      \else
       \ifcase\tcount@\char117 \or\or\char117 \or\char118 \or\char119 \or
       \char120 \or\or\or\or\or\char121 \or\char122 \or\arrow@i\char125
       \or\char117 \hskip\thr@@\p@\char117 \hskip-\thr@@\p@\fi
      \fi
     \fi
    \fi\fi
    \dimen@\mathaxis@\advance\dimen@.2\p@
    \dimen@ii\mathaxis@\advance\dimen@ii-.2\p@
    \ifnum\arrcount@=\m@ne
     \let\leads@\null
    \else
     \ifcase\arrcount@
      \def\leads@{\hrule\height\dimen@\depth-\dimen@ii}\or
      \def\leads@{\hrule\height\dimen@\depth-\dimen@ii}\or
      \def\leads@{\hbox to\ten@\p@{%
       \leaders\hrule\height\dimen@\depth-\dimen@ii\hfil
       \hfil
      \leaders\hrule\height\dimen@\depth-\dimen@ii\hskip\z@ plus2fil\relax
       \hfil
       \leaders\hrule\height\dimen@\depth-\dimen@ii\hfil}}\or
     \def\leads@{\hbox{\hbox to\ten@\p@{\dimen@\mathaxis@\advance\dimen@1.2\p@
       \dimen@ii\dimen@\advance\dimen@ii-.4\p@
       \leaders\hrule\height\dimen@\depth-\dimen@ii\hfil}%
       \kern-\ten@\p@
       \hbox to\ten@\p@{\dimen@\mathaxis@\advance\dimen@-1.2\p@
       \dimen@ii\dimen@\advance\dimen@ii-.4\p@
       \leaders\hrule\height\dimen@\depth-\dimen@ii\hfil}}}\fi
    \fi
    \cleaders\leads@\hfil
    \ifnum\arrcount@=\m@ne\else\ifnum\arrcount@=\thr@@\else
     \arrow@i
     \ifE@
      \ifnum\tcount@=\m@ne
      \else
       \ifcase\tcount@\char119 \or\or\char119 \or\char120 \or\char121 \or
       \char122 \or \or\or\or\or\char123 \or\char124 \or
       \char125 \or\char119 \hskip-\thr@@\p@\char119 \hskip\thr@@\p@\fi
      \fi
     \else
      \ifcase\scount@\or\or\char120 \or\char119 \or\or\or\char121 \or\char122
      \or\char123 \or\char124 \or\or\or\char125 \or
      \char119 \hskip-\thr@@\p@\char119 \hskip\thr@@\p@\fi
     \fi
    \fi\fi
    \harrow@ee}}%
  \harrow@e}%
 \iflabel@i
  \dimen@ii\z@\setbox\ZER@\hbox{$\m@th\tsize@@\label@i$}%
  \ifnum\arrcount@=\m@ne
  \else
   \advance\dimen@ii\mathaxis@
   \advance\dimen@ii\dp\ZER@\advance\dimen@ii\tw@\p@
   \ifnum\arrcount@=\thr@@\advance\dimen@ii\tw@\p@\fi
  \fi
  \advance\dimen@ii\pdimen@
  \next@{\harrow@b\smash{\raise\dimen@ii\hbox to\dimen@
   {\harrow@bb\hskip\tw@\ldimen@i\hfil\box\ZER@\hfil\harrow@ee}}\harrow@e}%
 \fi
 \iflabel@ii
  \ifnum\arrcount@=\m@ne
  \else
   \setbox\ZER@\hbox{$\m@th\tsize@\label@ii$}%
   \dimen@ii-\ht\ZER@\advance\dimen@ii-\tw@\p@
   \ifnum\arrcount@=\thr@@\advance\dimen@ii-\tw@\p@\fi
   \advance\dimen@ii\mathaxis@\advance\dimen@ii\pdimen@
   \next@{\harrow@b\smash{\raise\dimen@ii\hbox to\dimen@
    {\harrow@bb\hskip\tw@\ldimen@ii\hfil\box\ZER@\hfil\harrow@ee}}\harrow@e}%
  \fi
 \fi}
\let\tsize@\tsize
\def\tsizeCDlabels{\let\tsize@\tsize}
\def\ssizeCDlabels{\let\tsize@\ssize}
\def\tsize@@{\ifnum\arrcount@=\m@ne\else\tsize@\fi}
\def\varrow@{\dimen@\secondy@\advance\dimen@-\firsty@
 \ifN@\else\multiply\dimen@\m@ne\fi
 \setbox\ZER@\vbox to\dimen@
  {\ifN@\vskip-\Ydimen@\else\vskip\ydimen@\fi
   \ifnum\arrcount@=\m@ne\else\ifnum\arrcount@=\thr@@\else
    \hbox{\arrow@iii
     \ifN@
      \ifnum\tcount@=\m@ne
      \else
       \ifcase\tcount@\char117 \or\or\char117 \or\char118 \or\char119 \or
       \char120 \or\or\or\or\or\char121 \or\char122 \or\char123 \or
       \vbox{\hbox{\char117}\nointerlineskip\vskip\thr@@\p@
       \hbox{\char117}\vskip-\thr@@\p@}\fi
      \fi
     \else
      \ifcase\scount@\or\or\char118 \or\char117 \or\or\or\char119 \or
      \char120 \or\char121 \or\char122 \or\or\or\char123 \or
      \vbox{\hbox{\char117}\nointerlineskip\vskip\thr@@\p@
      \hbox{\char117}\vskip-\thr@@\p@}\fi
     \fi}%
    \nointerlineskip
   \fi\fi
   \ifnum\arrcount@=\m@ne
    \let\leads@\null
   \else
    \ifcase\arrcount@\let\leads@\vrule\or\let\leads@\vrule\or
    \def\leads@{\vbox to\ten@\p@{%
     \hrule\height1.67\p@\depth\z@\width.4\p@
     \vfil
     \hrule\height3.33\p@\depth\z@\width.4\p@
     \vfil
     \hrule\height1.67\p@\depth\z@\width.4\p@}}\or
    \def\leads@{\hbox{\vrule\height\p@\hskip\tw@\p@\vrule}}\fi
   \fi
  \cleaders\leads@\vfill\nointerlineskip
   \ifnum\arrcount@=\m@ne\else\ifnum\arrcount@=\thr@@\else
    \hbox{\arrow@iv
     \ifN@
      \ifcase\scount@\or\or\char118 \or\char117 \or\or\or\char119 \or
      \char120 \or\char121 \or\char122 \or\or\or\arrow@iii\char123 \or
      \vbox{\hbox{\char117}\nointerlineskip\vskip-\thr@@\p@
      \hbox{\char117}\vskip\thr@@\p@}\fi
     \else
      \ifnum\tcount@=\m@ne
      \else
       \ifcase\tcount@\char117 \or\or\char117 \or\char118 \or\char119 \or
       \char120 \or\or\or\or\or\char121 \or\char122 \or\arrow@iii\char123 \or
       \vbox{\hbox{\char117}\nointerlineskip\vskip-\thr@@\p@
       \hbox{\char117}\vskip\thr@@\p@}\fi
      \fi
     \fi}%
   \fi\fi
   \ifN@\vskip\ydimen@\else\vskip-\Ydimen@\fi}%
 \ifN@
  \dimen@ii\firsty@
 \else
  \dimen@ii-\firsty@\advance\dimen@ii\ht\ZER@\multiply\dimen@ii\m@ne
 \fi
 \rlap{\smash{\hskip\tocenter@\hskip\pdimen@\raise\dimen@ii\box\ZER@}}%
 \iflabel@i
  \setbox\ZER@\vbox to\dimen@{\vfil
   \hbox{$\m@th\tsize@@\label@i$}\vskip\tw@\ldimen@i\vfil}%
  \rlap{\smash{\hskip\tocenter@\hskip\pdimen@
  \ifnum\arrcount@=\m@ne\let\next@\relax\else\let\next@\llap\fi
  \next@{\raise\dimen@ii\hbox{\ifnum\arrcount@=\m@ne\hskip-.5\wd\ZER@\fi
   \box\ZER@\ifnum\arrcount@=\m@ne\else\hskip\tw@\p@\fi}}}}%
 \fi
 \iflabel@ii
  \ifnum\arrcount@=\m@ne
  \else
   \setbox\ZER@\vbox to\dimen@{\vfil
    \hbox{$\m@th\tsize@\label@ii$}\vskip\tw@\ldimen@ii\vfil}%
   \rlap{\smash{\hskip\tocenter@\hskip\pdimen@
   \rlap{\raise\dimen@ii\hbox{\ifnum\arrcount@=\thr@@\hskip4.5\p@\else
    \hskip2.5\p@\fi\box\ZER@}}}}%
  \fi
 \fi
}
\newdimen\goal@
\newdimen\shifted@
\newcount\Tcount@
\newcount\Scount@
\newbox\shaft@
\newcount\slcount@
\def\getcos@#1{%
 \ifnum\tan@i<\tan@ii
  \dimen@#1%
  \ifnum\slcount@<8 \count@9 \else \ifnum\slcount@<12 \count@8 \else
   \count@7 \fi\fi
  \multiply\dimen@\count@\divide\dimen@\ten@
  \dimen@ii\dimen@\multiply\dimen@ii\tan@i\divide\dimen@ii\tan@ii
 \else
  \dimen@ii#1%
  \count@-\slcount@\advance\count@24
  \ifnum\count@<8 \count@9 \else \ifnum\count@<12 \count@8
   \else\count@7 \fi\fi
  \multiply\dimen@ii\count@\divide\dimen@ii\ten@
  \dimen@\dimen@ii\multiply\dimen@\tan@ii\divide\dimen@\tan@i
 \fi}
\newdimen\adjust@
\def\Nnext@{\ifN@\let\next@\raise\else\let\next@\lower\fi}
\def\arrow@@{\slcount@\angcount@
 \ifNESW@
  \ifnum\angcount@<\ten@
   \let\arrowfont@\arrow@i\global\advance\angcount@\m@ne
   \global\multiply\angcount@13
  \else
   \ifnum\angcount@<19
    \let\arrowfont@\arrow@ii\global\advance\angcount@-\ten@
    \global\multiply\angcount@13
   \else
    \let\arrowfont@\arrow@iii\global\advance\angcount@-19
    \global\multiply\angcount@13
  \fi\fi
  \Tcount@\angcount@
 \else
  \ifnum\angcount@<5
   \let\arrowfont@\arrow@iii\global\advance\angcount@\m@ne
   \global\multiply\angcount@13 \global\advance\angcount@65
  \else
   \ifnum\angcount@<14
    \let\arrowfont@\arrow@iv\global\advance\angcount@-5
    \global\multiply\angcount@13
   \else
    \ifnum\angcount@<23
     \let\arrowfont@\arrow@v\global\advance\angcount@-14
     \global\multiply\angcount@13
    \else
     \let\arrowfont@\arrow@i\global\angcount@117
  \fi\fi\fi
  \ifnum\angcount@=117 \Tcount@115 \else\Tcount@\angcount@\fi
 \fi
 \Scount@\Tcount@
 \ifE@
  \ifnum\tcount@=\z@\advance\Tcount@\tw@\else\ifnum\tcount@=13
   \advance\Tcount@\tw@\else\advance\Tcount@\tcount@\fi\fi
  \ifnum\scount@=\z@\else\ifnum\scount@=13 \advance\Scount@\thr@@\else
   \advance\Scount@\scount@\fi\fi
 \else
  \ifcase\tcount@\advance\Tcount@\thr@@\or\or\advance\Tcount@\thr@@\or
  \advance\Tcount@\tw@\or\advance\Tcount@6 \or\advance\Tcount@7
  \or\or\or\or\or\advance\Tcount@8 \or\advance\Tcount@9 \or
  \advance\Tcount@12 \or\advance\Tcount@\thr@@\fi
  \ifcase\scount@\or\or\advance\Scount@\thr@@\or\advance\Scount@\tw@\or
  \or\or\advance\Scount@4 \or\advance\Scount@5 \or\advance\Scount@\ten@
  \or\advance\Scount@11 \or\or\or\advance\Scount@12 \or\advance
  \Scount@\tw@\fi
 \fi
 \ifcase\arrcount@\or\or\global\advance\angcount@\@ne\else\fi
 \ifN@\shifted@\firsty@\else\shifted@-\firsty@\fi
 \ifE@\else\advance\shifted@\charht@\fi
 \goal@\secondy@\advance\goal@-\firsty@
 \ifN@\else\multiply\goal@\m@ne\fi
 \setbox\shaft@\hbox{\arrowfont@\char\angcount@}%
 \ifnum\arrcount@=\thr@@
  \getcos@{1.5\p@}%
  \setbox\shaft@\hbox to\wd\shaft@{\arrowfont@
   \rlap{\hskip\dimen@ii
    \smash{\ifNESW@\let\next@\lower\else\let\next@\raise\fi
     \next@\dimen@\hbox{\arrowfont@\char\angcount@}}}%
   \rlap{\hskip-\dimen@ii
    \smash{\ifNESW@\let\next@\raise\else\let\next@\lower\fi
      \next@\dimen@\hbox{\arrowfont@\char\angcount@}}}\hfil}%
 \fi
 \rlap{\smash{\hskip\tocenter@\hskip\firstx@
  \ifnum\arrcount@=\m@ne
  \else
   \ifnum\arrcount@=\thr@@
   \else
    \ifnum\scount@=\m@ne
    \else
     \ifnum\scount@=\z@
     \else
      \setbox\ZER@\hbox{\ifnum\angcount@=117 \arrow@v\else\arrowfont@\fi
       \char\Scount@}%
      \ifNESW@
       \ifnum\scount@=\tw@
        \dimen@\shifted@\advance\dimen@-\charht@
        \ifN@\hskip-\wd\ZER@\fi
        \Nnext@
        \next@\dimen@\copy\ZER@
        \ifN@\else\hskip-\wd\ZER@\fi
       \else
        \Nnext@
        \ifN@\else\hskip-\wd\ZER@\fi
        \next@\shifted@\copy\ZER@
        \ifN@\hskip-\wd\ZER@\fi
       \fi
       \ifnum\scount@=12
        \advance\shifted@\charht@\advance\goal@-\charht@
        \ifN@\hskip\wd\ZER@\else\hskip-\wd\ZER@\fi
       \fi
       \ifnum\scount@=13
        \getcos@{\thr@@\p@}%
        \ifN@\hskip\dimen@\else\hskip-\wd\ZER@\hskip-\dimen@\fi
        \adjust@\shifted@\advance\adjust@\dimen@ii
        \Nnext@
        \next@\adjust@\copy\ZER@
        \ifN@\hskip-\dimen@\hskip-\wd\ZER@\else\hskip\dimen@\fi
       \fi
      \else
       \ifN@\hskip-\wd\ZER@\fi
       \ifnum\scount@=\tw@
        \ifN@\hskip\wd\ZER@\else\hskip-\wd\ZER@\fi
        \dimen@\shifted@\advance\dimen@-\charht@
        \Nnext@
        \next@\dimen@\copy\ZER@
        \ifN@\hskip-\wd\ZER@\fi
       \else
        \Nnext@
        \next@\shifted@\copy\ZER@
        \ifN@\else\hskip-\wd\ZER@\fi
       \fi
       \ifnum\scount@=12
        \advance\shifted@\charht@\advance\goal@-\charht@
        \ifN@\hskip-\wd\ZER@\else\hskip\wd\ZER@\fi
       \fi
       \ifnum\scount@=13
        \getcos@{\thr@@\p@}%
        \ifN@\hskip-\wd\ZER@\hskip-\dimen@\else\hskip\dimen@\fi
        \adjust@\shifted@\advance\adjust@\dimen@ii
        \Nnext@
        \next@\adjust@\copy\ZER@
        \ifN@\hskip\dimen@\else\hskip-\dimen@\hskip-\wd\ZER@\fi
       \fi
      \fi
  \fi\fi\fi\fi
  \ifnum\arrcount@=\m@ne
  \else
   \loop
    \ifdim\goal@>\charht@
    \ifE@\else\hskip-\charwd@\fi
    \Nnext@
    \next@\shifted@\copy\shaft@
    \ifE@\else\hskip-\charwd@\fi
    \advance\shifted@\charht@\advance\goal@-\charht@
    \repeat
   \ifdim\goal@>\z@
    \dimen@\charht@\advance\dimen@-\goal@
    \divide\dimen@\tan@i\multiply\dimen@\tan@ii
    \ifE@\hskip-\dimen@\else\hskip-\charwd@\hskip\dimen@\fi
    \adjust@\shifted@\advance\adjust@-\charht@\advance\adjust@\goal@
    \Nnext@
    \next@\adjust@\copy\shaft@
    \ifE@\else\hskip-\charwd@\fi
   \else
    \adjust@\shifted@\advance\adjust@-\charht@
   \fi
  \fi
  \ifnum\arrcount@=\m@ne
  \else
   \ifnum\arrcount@=\thr@@
   \else
    \ifnum\tcount@=\m@ne
    \else
     \setbox\ZER@
      \hbox{\ifnum\angcount@=117 \arrow@v\else\arrowfont@\fi\char\Tcount@}%
     \ifnum\tcount@=\thr@@
      \advance\adjust@\charht@
      \ifE@\else\ifN@\hskip-\charwd@\else\hskip-\wd\ZER@\fi\fi
     \else
      \ifnum\tcount@=12
       \advance\adjust@\charht@
       \ifE@\else\ifN@\hskip-\charwd@\else\hskip-\wd\ZER@\fi\fi
      \else
       \ifE@\hskip-\wd\ZER@\fi
     \fi\fi
     \Nnext@
     \next@\adjust@\copy\ZER@
     \ifnum\tcount@=13
      \hskip-\wd\ZER@
      \getcos@{\thr@@\p@}%
      \ifE@\hskip-\dimen@\else\hskip\dimen@\fi
      \advance\adjust@-\dimen@ii
      \Nnext@
      \next@\adjust@\box\ZER@
     \fi
  \fi\fi\fi}}%
 \iflabel@i
  \rlap{\hskip\tocenter@
  \dimen@\firstx@\advance\dimen@\secondx@\divide\dimen@\tw@
  \advance\dimen@\ldimen@i
  \dimen@ii\firsty@\advance\dimen@ii\secondy@\divide\dimen@ii\tw@
  \global\multiply\ldimen@i\tan@i\global\divide\ldimen@i\tan@ii
  \ifNESW@\advance\dimen@ii\ldimen@i\else\advance\dimen@ii-\ldimen@i\fi
  \setbox\ZER@\hbox{\ifNESW@\else\ifnum\arrcount@=\thr@@\hskip4\p@\else
   \hskip\tw@\p@\fi\fi
   $\m@th\tsize@@\label@i$\ifNESW@\ifnum\arrcount@=\thr@@\hskip4\p@\else
   \hskip\tw@\p@\fi\fi}%
  \ifnum\arrcount@=\m@ne
   \ifNESW@\advance\dimen@.5\wd\ZER@\advance\dimen@\p@\else
    \advance\dimen@-.5\wd\ZER@\advance\dimen@-\p@\fi
   \advance\dimen@ii-.5\ht\ZER@
  \else
   \advance\dimen@ii\dp\ZER@
   \ifnum\slcount@<6 \advance\dimen@ii\tw@\p@\fi
  \fi
  \hskip\dimen@
  \ifNESW@\let\next@\llap\else\let\next@\rlap\fi
  \next@{\smash{\raise\dimen@ii\box\ZER@}}}%
 \fi
 \iflabel@ii
  \ifnum\arrcount@=\m@ne
  \else
   \rlap{\hskip\tocenter@
   \dimen@\firstx@\advance\dimen@\secondx@\divide\dimen@\tw@
   \ifNESW@\advance\dimen@\ldimen@ii\else\advance\dimen@-\ldimen@ii\fi
   \dimen@ii\firsty@\advance\dimen@ii\secondy@\divide\dimen@ii\tw@
   \global\multiply\ldimen@ii\tan@i\global\divide\ldimen@ii\tan@ii
   \advance\dimen@ii\ldimen@ii
   \setbox\ZER@\hbox{\ifNESW@\ifnum\arrcount@=\thr@@\hskip4\p@\else
    \hskip\tw@\p@\fi\fi
    $\m@th\tsize@\label@ii$\ifNESW@\else\ifnum\arrcount@=\thr@@\hskip4\p@
    \else\hskip\tw@\p@\fi\fi}%
   \advance\dimen@ii-\ht\ZER@
   \ifnum\slcount@<9 \advance\dimen@ii-\thr@@\p@\fi
   \ifNESW@\let\next@\rlap\else\let\next@\llap\fi
   \hskip\dimen@\next@{\smash{\raise\dimen@ii\box\ZER@}}}%
  \fi
 \fi
}
\def\outCD@#1{\def#1{\Err@{\noexpand#1must not be used within \string\CD}}}
\newskip\preCDskip@
\newskip\postCDskip@
\preCDskip@\z@
\postCDskip@\z@
\def\preCDspace#1{\RIfMIfI@
 \onlydmatherr@\preCDspace\else\advance\preCDskip@#1\relax\fi\else
 \onlydmatherr@\preCDspace\fi}
\def\postCDspace#1{\RIfMIfI@
 \onlydmatherr@\postCDspace\else\advance\postCDskip@#1\relax\fi\else
 \onlydmatherr@\postCDspace\fi}
\def\predisplayspace#1{\RIfMIfI@
 \onlydmatherr@\predisplayspace\else
 \advance\abovedisplayskip#1\relax
 \advance\abovedisplayshortskip#1\relax\fi
 \else\onlydmatherr@\preCDspace\fi}
\def\postdisplayspace#1{\RIfMIfI@
 \onlydmatherr@\postdisplayspace\else
 \advance\belowdisplayskip#1\relax
 \advance\belowdisplayshortskip#1\relax\fi
 \else\onlydmatherr@\postdisplayspace\fi}
\def\PreCDSpace#1{\global\preCDskip@#1\relax}
\def\PostCDSpace#1{\global\postCDskip@#1\relax}
\def\CD#1\endCD{%
 \outCD@\cgaps\outCD@\rgaps\outCD@\Cgaps\outCD@\Rgaps
 \preCD@#1\endCD
 \advance\abovedisplayskip\preCDskip@
 \advance\abovedisplayshortskip\preCDskip@
 \advance\belowdisplayskip\postCDskip@
 \advance\belowdisplayshortskip\postCDskip@
 \vcenter{\offinterlineskip
  \vskip\preCDskip@\Let@\global\colcount@\@ne\global\rowcount@\z@
  \everycr{%
   \noalign{%
    \ifnum\rowcount@=\Rowcount@
    \else
     \getrgap@\rowcount@\vskip\getdim@
     \global\advance\rowcount@\@ne\global\colcount@\@ne
    \fi}}%
  \tabskip\z@
  \halign{&\global\xoff@\z@\global\yoff@\z@
   \getcgap@\colcount@\hskip\getdim@
   \hfil\vrule\height\ten@\p@\width\z@\depth\z@
   $\m@th\displaystyle{##}$\hfil
   \global\advance\colcount@\@ne\cr
   #1\crcr}\vskip\postCDskip@}%
 \preCDskip@\z@\postCDskip@\z@
 \def\getcgap@##1{\ifcase##1\or\getdim@\z@\else\getdim@\standardcgap\fi}%
 \def\getrgap@##1{\ifcase##1\getdim@\z@\else\getdim@\standardrgap\fi}%
 \let\Width@\relax\let\Height@\relax\let\Depth@\relax\let\Rowheight@\relax
 \let\Rowdepth@\relax\let\Colwidth@\relax
}
\let\enddocument\bye
\def\alloc@#1#2#3#4#5{\global\advance\count1#1by\@ne
  \ch@ck#1#4#2%
  \allocationnumber=\count1#1%
  \global#3#5=\allocationnumber
  \wlog{\string#5=\string#2\the\allocationnumber}}
\catcode`\@=\active

\catcode`\"=12   
\font\blackital=cmmib10  \skewchar\blackital='177 \font\sblackital=cmmib7
\skewchar\sblackital='177 \font\ssblackital=cmmib5  \skewchar\ssblackital='177
\font\sanss=cmss10 \font\ssanss=cmss8 scaled 900 \font\sssanss=cmss8 scaled 600
\font\blackboard=msbm10 \font\sblackboard=msbm7 \font\ssblackboard=msbm5
\font\caligr=eusm10 \font\scaligr=eusm7 \font\sscaligr=eusm5 
\font\fraktur=eufm10 \font\sfraktur=eufm7 \font\ssfraktur=eufm5

\def\tx#1{{\fam0\relax#1}}

\newfam\bifam
\textfont\bifam=\blackital \scriptfont\bifam=\sblackital
\scriptscriptfont\bifam=\ssblackital
\def\bi#1{{\fam\bifam\relax#1}}


\newfam\bbfam
\textfont\bbfam=\blackboard \scriptfont\bbfam=\sblackboard
\scriptscriptfont\bbfam=\ssblackboard
\def\bb#1{{\fam\bbfam\relax#1}}

\newfam\ssfam
\textfont\ssfam=\sanss \scriptfont\ssfam=\ssanss \scriptscriptfont\ssfam=\sssanss
\def\ss#1{{\fam\ssfam\relax#1}}

\newfam\clfam
\textfont\clfam=\caligr \scriptfont\clfam=\scaligr \scriptscriptfont\clfam=\sscaligr
\def\cl#1{{\fam\clfam\relax#1}}

\newfam\frfam
\textfont\frfam=\fraktur \scriptfont\frfam=\sfraktur
\scriptscriptfont\frfam=\ssfraktur

\mathchardef\za="710B  
\mathchardef\zb="710C  
\mathchardef\zg="710D  
\mathchardef\zd="710E  
\mathchardef\zve="710F 
\mathchardef\zz="7110  
\mathchardef\zh="7111  
\mathchardef\zvy="7112 
\mathchardef\zi="7113  
\mathchardef\zk="7114  
\mathchardef\zl="7115  
\mathchardef\zm="7116  
\mathchardef\zn="7117  
\mathchardef\zx="7118  
\mathchardef\zp="7119  
\mathchardef\zr="711A  
\mathchardef\zs="711B  
\mathchardef\zt="711C  
\mathchardef\zu="711D  
\mathchardef\zvf="711E 
\mathchardef\zq="711F  
\mathchardef\zc="7120  
\mathchardef\zw="7121  
\mathchardef\ze="7122  
\mathchardef\zy="7123  
\mathchardef\zvp="7124 
\mathchardef\zvr="7125 
\mathchardef\zvs="7126 
\mathchardef\zf="7127  
\mathchardef\zG="7000  
\mathchardef\zD="7001  
\mathchardef\zY="7002  
\mathchardef\zL="7003  
\mathchardef\zX="7004  
\mathchardef\zP="7005  
\mathchardef\zS="7006  
\mathchardef\zU="7007  
\mathchardef\zF="7008  
\mathchardef\zC="7009  
\mathchardef\zW="700A  

\catcode`\"=\active
\loadmsam
\newsymbol\leqslant 1336
\newsymbol\geqslant 133E
\newsymbol\boxtimes 1202

\def\leqs{\leqslant}

\def\*{{\textstyle *}}

\def\tT{\widetilde{\sT}}

\def\proof{\demo{Proof}}
\def\endproof{\hfill \vrule height4pt width6pt depth2pt \enddemo}

\def\R{{\bb R}}

\def\bA{{\bi A}}
\def\bB{{\bi B}}

\def\bY{{\bi Y}}

\def\bV{{\bi V}}

\def\bT{{\bi T}}
\def\bC{{\bi C}}

\def\bZ{{\bi Z}}

\def\sC{{\ss C}}

\def\sL{{\ss L}}

\def\sP{{\ss P}}

\def\sT{{\ss T}}
\def\sV{{\ss V}}

\def\sv{{\ss v}}

\def\xd{\tx{d}}
\def\xi{\tx{i}}

\def\dim{\operatorname{dim}}

\def\cA{{\cl A}}

\def\ot{\boxtimes}
\def\bA{{\bi A}}
\def\bB{{\bi B}}
\def\bI{{\bi I}}
\def\bZ{{\bi Z}}
\def\xc{\tx{c}}
\font\eurm=eurm10
\def\eu#1{\hbox{$\eurm #1$}}
\def\ezw{\eu \zw}
\def\ezp{\eu\zp}
\def\ezt{\eu\zt}
\def\ezk{\eu\zk}

\def\ezd{\eu\zd}
\def\eza{\eu\za}
\def\tT{\widetilde{\sT}}
\def\bT{\overline{\sT}}

\def\dt{\xd_\sT}
\def\dl{\xd_\sL}
\def\dtt{\xd_{\tT Z}}
\def\dbt{\xd_{\bT Z}}

\input paperbcp.st\relax
\RefWarnings \TagsOnRight \overfullrule=0pt
\document
   \title
     AN AFFINE FRAMEWORK \\ FOR ANALYTICAL MECHANICS
   \endtitle
    \abbrevtitle{Affine framework}
   \author
                PAWE\L\ URBA\'NSKI
    \endauthor
    \affil
                Division of Mathematical Methods in Physics \\
                University of Warsaw \\
                Ho\D{z}a 74, 00-682 Warszawa, Poland \\
                E-mail: urbanski\@fuw.edu.pl
   \endaffil
    \abbrevauthors{P. Urba\'nski}
   \thanks{Supported by KBN, grant No 2 PO3A 041 18}
    \subjclass{17B55, 58C11}
    \keywords{affine bundle, symplectic manifold, contact manifold,
Legendre trans\-formation.}
   \abstract
An affine Cartan calculus is developed. The concepts of special affine
bundles and special affine duality are introduced.  The canonical
isomorphisms, fundamental for Lagrangian and Hamiltonian formulations of
the dynamics in the affine setting are proved.
\endabstract

   \maketitle
        \hl1{Introduction}
    Gauge independence of the Langrangian formulation of dynamics of
charged particles  can be achieved by increasing the dimension of the
configuration space of the particle (\Ref{MT}]). The four dimensional
space-time of general relativity is replaced by the five dimensional
space-time-phase of Kaluza. The phase space of the particle is the
cotangent bundle of the Kaluza space and the gauge independent Lagrangian
is a function on the tangent bundle of the Kaluza space.
    A similar construction makes possible a frame-independent formulation
of the Newtonian analytical mechanics (see [\Ref{Tu1}] for details).

    There is an alternate approach, based on ideas and results of
Tulczyjew, in which  the four dimensional space-time is used as the
configuration space of the charged particle ([\Ref{TU}]). The phase space
is no longer a cotangent bundle and not even a vector bundle. It is an
affine bundle modeled on the cotangent bundle $\sT ^{\textstyle *}M $
of the space-time manifold $M$.
    The dynamics is a generalized Hamiltonian system, but the Lagrangian
generating object is not a function. It is a section of an affine bundle
modeled on $\sT M\times \R$.

    Also time-dependent mechanics in the inhomogeneous formulation requires
affine objects. Here infinitesimal configurations are first jets of
motions interpreted as sections of the space-time fibered over time, and
are elements of an affine bundle over the space-time. The phase manifold is
a vector bundle, however  Hamiltonian is not a function but a section of a
bundle over the phase manifold ([\Ref{U1}]).

    We  develope the geometric framework for these approaches. The
standard geometric constructions based on the algebra of functions on a
manifold $M$ are replaced by constructions   based on the
affine space of sections of an affine bundle $\zz\colon Z\rightarrow M$,
modeled on the trivial bundle $M\times \R$.

    In Section~2 we adopt the definition of a covector as an equivalence
class of functions to the affine case. The manifold $\sP Z$ of affine covectors is
an affine bundle modeled on $\sT^{\textstyle *} M$ and carries a canonical
symplectic structure. Sections of $\zz$  generate Lagrangian submanifolds
of $\sP Z$.

    In analytical mechanics of a particle the phase space is considered a
vector bundle, dual to  the vector bundle of infinitesimal configurations
(velocities). In the presented approach vector bundles are replaced by affine bundles,
vector spaces by affine spaces. Affine functions on an affine space $\underline{A}$
are replace by affine sections of a bundle $\zz_A\colon A\rightarrow \underline{A}$ in
the category of affine space, modeled on $\underline{A}\times \R$. This means that $A,
\underline{A}$ are affine spaces and $\zz_A$ is an affine surjective mapping.  Such a
bundle is uniquely determined be the affine space $A$ and a distinguished vector $v_1$
in the model vector space $\sV(A)$. We call it a special affine space (Section~3).

    The dual object to a special affine space $\bA=(A,v_1)$ is the set of
affine sections of $\zz_A$.  It carries a natural structure of an affine
space with the vector space of affine functions on $\underline{A}$ as the
model  vector space. With the constant function $1$  the dual affine space
becomes a special affine space. We introduce in an obvious way the concept
of a special affine bundle. Principal examples of special affine bundles
are $\sP Z\times \R\rightarrow M$, the first jet (contact) bundle $\sC
Z\rightarrow M$, and their duals.

    It is convenient to represent an element of the dual object as a
morphism in the category. For this purpose we identify a  section $\zs$ of
$\zz_A$ with a morphism $\widetilde{\zf}_\zs\colon \bA \rightarrow \bI$,
where $\bI= (\R,1)$ and $\widetilde{\zf}(a) = a- \zf(\zz_A(a))$.
    We must be aware of the fact that with this representation the
canonical pairing is not symmetric, but skew-symmetric.

The theory of special affine spaces (bundles) and duality is developed in
Sections~3, 4 and 5. In Section~3 we present main algebraic constructions in
the category of special affine spaces.  In Section~4 we interpret the
bundles dual to $\sP Z\times \R$ and $\sC Z $ in terms of the tangent
bundle $\sT Z$. In Section~5 we prove the existence of canonical
isomorphisms between  $\sC\bA$ and $\sC \bA^\#$, where $\bA$ is a special
affine bundle and $\bA^\#$ is the special affine dual bundle. This isomorphism
corresponds to the canonical isomorphism between  $\sT^{\textstyle *} E$
and $\sT^{\textstyle *} E^{\textstyle *} $ for a vector bundle $E$, and
makes possible the Legendre transformation.

    The Lagrange formulation of the dynamics of a particle is possible
because of the canonical Tulczyjew isomorphism (symplectomorphism) of
$\sT\sT^{\textstyle *} M$ and $\sT^{\textstyle *} \sT M$. In Section~5  we
give a proof of an affine counterpart for this isomorphism. As a
consequence we obtain an affine setting for the  Lagrangian formulation of
the dynamics. As an example we discuss the case of a relativistic charged
particle (Section~6).

    A similar approach to time-dependent non-relativistic mechanics has
been recently developed be E.~Massa with collaborators ([\Ref{MP},
\Ref{MV}, \Ref{Vi}]), W. Sarlet, and E. Mart\'\i nez ([\Ref{Sa2},
\Ref{Sa}, \Ref{Sa1}]).

    This work is a contribution to a programme of study of geometric
foundations of physical theories conducted jointly with Professor
Tulczyjew at the University of Camerino.
        \hl1{Affine phase and contact spaces}
    \hl2{Affine spaces  and  affine bundles}
  An {\it affine space} is a triple $(A,V,\za)$, where $A$ is a set, $V$ is a real
vector space of finite dimension and $\za$ is a mapping $\za \colon A \times A
\rightarrow V$ such that
      \list
   \item $\za(a_3,a_2) + \za(a_2,a_1) + \za(a_1,a_3) = 0$;
   \item the mapping $\za(\cdot,a) \colon A \rightarrow V$ is bijective for each $a
\in A$.
      \endlist
   We shall also write simply $A$ to denote the affine space $(A,V,\za)$ and
$\sV(A)$ to denote $V$.
    If $(A,V,\za)$ then also $(A,V,-\za)$ is an affine space. We
will write for brevity $\overline{A}$ to denote the affine space $(A,V,-\za)$.
   We will write also $a_2 - a_1$ instead of $\za(a_2,a_1)$ and  we will denote
by $a + v$ the unique point $a' \in A$ such that $a' - a = v$.

   Let $\zx \colon E \rightarrow M$ be a vector bundle.  An {\it affine
bundle modeled on} $\zx $ is a differential fibration $\eta \colon A
\rightarrow M$ and a differentiable mapping $\zr \colon A \times _M A
\rightarrow E$ such that
   \list
   \item $\zx \circ \zr = \eta \times _M \eta $,
   \item $\zr (a_3,a_2) + \zr (a_2,a_1) = \rho (a_3,a_1)$ for each triple
$(a_3,a_2,a_1) \in A \times _M A \times _M A$,
   \item for each local section $\zs \colon U \rightarrow A$ of $\eta $, the
mapping $\sigma _0\colon  \rightarrow E$ defined by $\zs _0(m) = \zr (\sigma
(m),\sigma (m))$ is the zero section of $\zx $ over $U$,
   \item for each local section $\sigma \colon U \rightarrow A$ of $\eta $, the
mapping $\zr _{\zs} \colon \eta ^{-1}(U) \rightarrow \zx ^{-1}(U)$ defined by
$\zr _{\sigma}(a) = \zr (a,\zs (\eta (a)))$ is a diffeomorphism.
   \endlist
   We will  write simply $A$ to denote the affine bundle $(A,E,\zr)$ and
$\sV(A)$ to denote $E$.

   \noindent {\bf Remark.} Any section $\zs\in \ss{Sec}(\eta)$ of $A$ induces an
obvious isomorphism $I_\zs\in \ss{Aff}_M(A,\sV(A))$ of affine bundles:
                $$ A_m\ni a\mapsto I_\zs(a)= a-\zs(m) \in \sV(A_m).
                                         \tag \label{Fa25}$$

    Let $\zt_i\colon A_i\rightarrow M$ be an affine bundle modeled on a
vector bundle $\sv(\zt_i)\colon
\sV(A_i)\rightarrow M$, $i=1,2,3$. Note that the space $\cA_i$ of sections
of $\zt_i$ is an affine space modeled on the space $\sV(\cA_i)$ of
sections of $\sv(\zt_i)$.
    For an affine bundle morphism $\zvf\colon A_1\rightarrow A_2$ we
denote by $\zvf_\sv \colon \sV({A_1})\rightarrow \sV({A_2})$ its
linear part, i.e.
        $$\zvf_\sv(v) = \zvf(a+v) - \zvf(a) \ \ \text{for}\ \ a\in A,\ v\in
\sV(A), \ \ \zt_1(a) = \sv(\zt_1)(v).
                                                         \tag \label{Fa25a}$$
\noindent
        We will denote by $\ss{Aff}_M (A_1,A_2)$ (resp.
$\ss{Hom}_M(V_1,V_2)$) the set of  morphisms over the identity on the base in the affine
(resp. vector) case. We shall also write $\ss{Aff}(A,\R)$ instead of
$\ss{Aff}_M(A, M\times \R)$ and $\ss{Lin}(V,\R)$ instead of
$\ss{Hom}_M(V,M\times \R)$.
   For a bi-affine mapping
            $$ F\colon A_1\times _M A_2 \rightarrow A_3
                                                            $$
we denote by $F^\sv$ and ${}^\sv F$, respectively, the mappings
                $$\align
    F^\sv&\colon A_1\times_ M \sV(A_2) \rightarrow \sV(A_3) \\
    &\colon (a_1,v_2)  \mapsto (F(a_1,\cdot ))_\sv(v_2)
                                          \tag \label{Fa25b}\endalign$$
and
            $$\align
    {}^\sv F&\colon \sV(A_1)\times_ M A_2 \rightarrow \sV(A_3) \\
    &\colon (v_1,a_2)  \mapsto(F(\cdot,a_2 ))_\sv(v_1).
                                           \tag \label{Fa25c}\endalign $$

\noindent   These mappings are, respectively, affine-linear and linear-affine in
the obvious sense.
By $F_\sv$ we denote the bilinear part of $F$, i.e.
        $$\align
    F_\sv&\colon \sV(A_1) \times _M \sV(A_2) \rightarrow \sV (A_3) \\
    &\colon  (v_1,v_2) \mapsto(F^\sv(\cdot ,v_2))_\sv(v_1) =
({}^\sv F(v_1,\cdot))_\sv(v_2)\\
    &= F(a_1 +v_1, a_2 +v_2) - F(a_1 + v_2, a_2) + F(a_1,a_2) -
F(a_1, a_2 +v_2).
                           \tag \label{Fa25d}\endalign$$

   \hl2{Affine Cartan calculus}
   The standard Cartan calculus of differential forms is based on the algebra of
differentiable functions on a manifold $M$. In the affine Cartan calculus we replace
functions by sections of an affine bundle $\zz \colon Z\rightarrow M$ modeled on
the trivial bundle $M\times \R$. We can consider $Z$ a principal bundle
with the structure group $(\R,+)$ (\Ref{TUZ}). The space of sections of
$\zz$ is an affine space modeled on the space of sections of $M\times \R$
and a section of $M\times \R$ we identify with a function on $M$.

   \hl2{The phase and the contact fibrations}
   Let $\zz \colon Z \rightarrow M$ be an affine fibration modeled on the
trivial fibration $pr_M\colon M \times \R \rightarrow M$.  We define an equivalence
relation in the set of all pairs $(m,\zs )$, where $m$ is a point in $M$ and $\sigma $
is a section of $\zz $.  Two pairs $(m,\zs )$ and $(m',\zs ')$ are equivalent if $m' =
m$ and $\xd(\zs ' - \zs)(m) = 0$.  We denote by $\sP Z$ the set of equivalence
classes. The class of $(m,\zs )$ will be denoted by $\xd\sigma (m)$ or by $\xd_m\zs$
and will be called the {\it differential} of $\zs $ at $m$. We define a caonical
projection
        $$\align
    {\sP}\zz&\colon \sP Z\rightarrow M \\
    &\colon \xd\zs(m) \mapsto m .
                                                \tag \label{Fa25e}\endalign$$

 We define also a mapping
        $$\align
    \sP\zr&\colon \sP Z\times _M\sP Z\rightarrow \sT^{\textstyle *} M \\
    &\colon (\xd\zs _2(m), \xd\zs _1(m)) \mapsto \xd(\zs _2 - \zs _1)(m)
                                                \tag \label{Fa25f}\endalign$$

 The pair $(\sP\zz ,\sP\zr )$ makes $\sP Z$  an affine fibration
modeled on the cotangent fibration $\zp _M \colon \sT^{\textstyle *} M
\rightarrow M$.  This fibration is called the {\it phase fibration} of
$\zz$. A section of $\sP\zz$ will be called an {\it affine 1-form}.

   Let $\za\colon M\rightarrow \sP Z$ be an affine 1-form and let $\zs$ be a
section of $\zz$. The differential $\xd_m(\za -\xd \zs)$ does not depend on the
choice of $\zs$ and will be called the {\it differential of $\za$ at $m$}. We will
denote it by $\xd\za(m)$ or by $\xd_m\za$. The differential of an affine 1-form is
an ordinary 2-form.

    \noindent
   {\bf Remark} Let us consider  $Z$  a principal bundle with the
structure group $(\R,+)$.   A section of $\sP Z$ defines a connection in
$Z$ and can be interpreted as an affine form of this connection.

   We define another equivalence relation in the set of all pairs $(m,\zs
)$. Two pairs $(m,\zs )$ and $(m',\zs ')$ are equivalent if $m' = m$,
$\zs(m) = \zs'(m)$, and $\xd(\zs ' - \zs)(m) = 0$.  We can identify the
equivalence class of $(m,\zs)$ with the first jet of the section $\zs$
with the source point $m$.  We denote by $\sC Z$ the set of equivalence
classes. The class of $(m,\zs )$ will be denoted by $\xc\sigma (m)$ or by
$\xc_m\zs$ and will be called the {\it contact element} of $\zs $ at $m$.
We define a mapping ${\sC}\zeta \colon \sC Z \rightarrow M$ by $\sC\zz
(\xc\zs (m)) = m$ and a mapping
    $$\sC\zr \colon \sC Z \times _M \sC Z \rightarrow \sT^{\textstyle *}
M\times \R $$
    by
    $$\sC\zr (\xd\zs _2(m), \xc\zs _1(m)) = (\xd(\zs _2 - \zs _1)(m),
\zs_2(m)- \zs_1(m)).$$
    The pair $(\sC\zz ,\sC\zr )$ makes $\sC Z$  an affine fibration
modeled on the fibration $\zg _M \colon \sT^{\textstyle *} M \times \R
\rightarrow M$. This fibration is called the {\it contact fibration} of
$\zz$.  We have an obvious isomorphism of affine bundles

        $$  \sC Z = \sP \times _M Z .
                                                            $$

   Let $Z$ and $Z'$ be special affine bundles modeled on $M\times \R$ and
$M'\times \R$ respectively. A morphism $\zF\colon Z\rightarrow Z'$
induces a mapping $\zF^{\textstyle *}\colon  \ss{Sec}(Z') \rightarrow
\ss{Sec(Z)}$. The mapping $\zF^{\textstyle *} $ and the relation
$\zF^{-1}$ project to relations $\sP\zF \colon \sP Z' \rightarrow \sP Z$
and $\sC\zF \colon \sC Z' \rightarrow \sC Z$.

   \hl2{Canonical structures on $\sP Z$ and $\sC Z$}
   We show first  that $\sP Z$ carries a canonical symplectic structure. For a choosen
section $\zs$ of $\zz$ we have isomorphisms (see \Ref{Fa20})
      $$\align
   I_\zs&\colon Z\rightarrow M\times \R  \\
      I_{\xd\zs}&\colon \sP Z\rightarrow \sT^{\textstyle *} M
                                       \tag \label{Fa29}\endalign$$
   and for two sections $\zs, \zs'$ the mapping
      $$\align
   I_{\xd\zs'}\circ I_{\xd\zs}^{-1}&\colon \sT^{\textstyle *} M\rightarrow
\sT^{\textstyle *} M  \\
   &\colon p \mapsto p +\xd(\zs -\zs')(\ezp_M(p))
                                       \tag \label{Fa30}\endalign$$
   is a symplectomorphism. It follows that a two-form $I_{\xd\zs}^{\textstyle *}
\ezw_M$, where $\ezw_M$ is the canonical symplectic form on $\sT^{\textstyle *}
M$, does not depend on the choice of $\zs$. We will denote this form by
$\ezw_\bZ$.

   There is a canonical projection
   $$\zm\colon \sC Z \rightarrow Z
                              \tag\label{Fa32}$$
   which is a morphism of  affine bundles $\zz_{\sC \bZ}\colon  \sC Z
\rightarrow \sP Z$ and $\zz \colon Z\rightarrow M$, so we have a pull-back of
sections of $\zz$ to sections of $\zz_{\sC\bZ}$. Now we can define a section
$\zvy_\bZ$ of $\sP\zz_{\sC\bZ}\colon \sP \sC Z \rightarrow \sP Z$ by
      $$ \zvy_\bZ(p) = \xd_p \zm^{\textstyle *} \zs_p,
                                                      \tag \label{Fa33}$$
   where $\zs_p$ is a section of $\zz$ which represents $p\in \sP Z$.
   It is easy to see that the canonical symplectic form on $\sP Z$ is equal $\xd
\zvy_\bZ$. The affine 1-form $\zvy_\bZ$ is called the {\it Liouville affine form}
of $\sC Z$ and defines the {\it canonical contact structure} of $\sC Z$.

   \hl2{Generating objects}
   As in the cotangent bundle, the image of an affine 1-form $\za \colon
M\rightarrow \sP Z$ is a Lagrangian submanifold of $(\sP Z, \ezw_\bZ)$ if and only
if it is closed, i.e. if $\xd\za =0$. If $\za =\xd \zs$ for a section $\zs$ of $Z$
then we say that $\zs$ is a generating section of the Lagrangian submanifold
$\za(M)$.

    Now, let $C\subset M$ be a submanifold and $\zs\colon C\rightarrow Z$ a
section of $\zz$ over $C$. We define a submanifold $L\subset \sP Z$ by
        $$ L =\{p\in \sP Z\colon p=\xd_m\zs' \ \text{where}\ m\in C \
\text{and} \ \zs'|C = \zs\}.
                                                    \tag \label{Fa33a}$$
    As in the trivial case, $L$ is a Lagrangian submanifold of $(\sP Z,
\ezw_\bZ)$. We say that $L$ is generated by the section $\zs$ over the
constraints manifold $C$.

    Let $\varsigma\colon Y\rightarrow N$ be an affine bundle modeled on
$N\times \R$ and let $\eta \colon Y\rightarrow Z$ be an epimorphism of
affine bundles such that $\eta(y+r) = \eta(y) +r$ (it is a morphism of
principal bundles). We denote by $\underline{\eta}$ the underlying
epimorphism $\underline{\eta}\colon N\rightarrow M$. It follows that for
$m\in M$ the subbundle $\eta^{-1}(Z_m)$ is trivial (has distinguished
constant sections). Let $\zs$ be a section of $\varsigma$. We say that a
point $n\in N$ is {\it critical} for $\zs$ if
        $$\xd(\zs|_{N_m})(n) = \xd\zs_0(n)
                                                        \tag \label{F33b}$$
    for a constant section $\zs_0$ of $Y|N_m$. Let $S(\zs)$ be the set of
critical points of $\zs$.   The condition \Ref{F33b} implies that there
exists a section $\zs'$ of $\zz$ such that $\xd_n(\zs -\eta^{\textstyle
*}\zs') =0 $. We define a mapping
    $$\chi\colon S\rightarrow \sP Z \colon n\mapsto  \xd_m \zs'. $$

            \claim \c{d}{Definition}{}{\rm               \label{Caa}
        A section $\zs$ of $\varsigma\colon Y\rightarrow N$ is an {\it
affine Morse family} if for a section $\zs'$ of $\zz\colon Z\rightarrow M$
the function $\zs - \eta^{\textstyle *} \zs'$ on $N$ is a Morse family.
        }\endclaim

    If $\zs$ is an affine Morse family then $\chi(S)$ is a Lagrangian submanifold  of
$(\sP Z, \ezw_\bZ)$.

        \hl1{Special affine spaces and duality}

        In the category of vector spaces, the dual object is defined as a
space of linear functions. In the presented here approach  functions are
replaced by sections  of a bundle. In particular, linear functions on a
vector space $\underline{V}$ are replaced by linear sections of a
fibration $\zz \colon V\rightarrow \underline{V}$ (modeled on
$\underline{V}\times \R$) in the category of vector spaces. We observe
that the space  $V^\ddag$ of linear sections  is affine, not linear.
Morever, an affine section can be uniquely represented by a liner section
and a number, i.e. an alement of the trivial bundle $V^\ddag\times \R$.
This observation leads us to the concept of special affine space (and its
model  special vector space) which puts bundles $\zz \colon V\rightarrow
\underline{V}$ and its dual $V^\ddag \times \R$ into  one category.

   A {special vector space} is a pair $\bV =(V,v_1)$, where $V$ is a vector space
and $v_1\in V$ is  a distinguished, non-zero vector.
   A {\it special affine space} is an affine space modeled on a special vector
space.
    Let $\bV =(V,v_1)$ and $\bV' = (V',v_1')$ be special vector spaces. A linear
mapping  $F\colon V\rightarrow V'$  is called a {\it morphism of special vector
spaces} if $F(v_1) = v_1'$.
    A {\it morphism of special affine spaces} is an affine mapping such that its
linear part is a morphism of special vector spaces.
   Let $\bA = (A,v_1)$ be a special affine space with the distinguished vector
$v_1\in \sV(A)$. There is a canonical action of $(\R,+)$ on $A$ given by the
formula
        $$ A\times \R\ni (a,r) \rightarrow a +rv_1.
                                                            \tag \label{Faa}$$
    The space of orbits is an affine space modeled on the quotient vector
space $\underline{V} = V/\{v_1\}$. We denote it by $\underline{A}$ and the
canonical projection $A\rightarrow \underline{A}$ by $\zz_\bA$. With the
action~\Ref{Faa} the fibration $\zz_\bA\colon A\rightarrow \underline{A}$
is an affine bundel modeled on $\underline{A}\times \R$.

    By $\overline \bA $ we denote the special affine space
$(\overline{A},v_1)$.
   Let $\bA=(A,v_1)$ and $\bB=(B,w_1)$ be special affine spaces. The
product $A\times B$ is an affine space modeled on $\sV(A)\times \sV(B)$.
Let us denote by $A\ot B$  the quotient affine space $(A\times B)/L$,
where $L\subset \sV(A) \times \sV(B)$ is a one-dimensional vector subspace
spanned by the vector $(v_1,-w_1)$. We have $\sV(A\ot B) = (\sV(A)\times
\sV(B))/L$ with a distinguished vector $v_1\ot w_1 =[(v_1,0)] = [(0,w_1)]$. The
{\it product} of special affine spaces is a special affine space $\bA\ot
\bB = (A\ot B, v_1\ot w_1)$ modeled on a special vector space $(\sV(A)\ot
\sV(B), v_1\ot w_1)$, where $\sV(A)\times \sV(B) = (\sV(A)\times \sV(B))/L$.

   We have obvious isomorphisms
      $$ \align \bA\ot \bB &\simeq \bB\ot \bA\\
            \bA\ot (\bB\ot \bC) &\simeq (\bA\ot \bB)\ot \bC.
                                    \tag \label{Fa1}\endalign$$

    Since $A$ and $\underline{A}$ are affine spaces, we can distinguished
the set of affine sections of $\zz_\bA$.  It is an affine space modeled on
the vector space of affine functions on $\underline{A}$:
        $$ \zs(\underline{a}) = (\zs -\zs')(\underline{a})\cdot v_1 +
\zs'(\underline{a}) .
                                                            \tag \label{Fa2}$$
 We consider this
space as the {\it dual space}  to $\bA$.

    In order to interpret the dual object in terms of morphisms, we
introduce the distinguished  special affine space $\bI = (\R,1)$.
    There is one-to-one correspondence between affine sections of the
fibration $\zz_A\colon A\rightarrow \underline{A} $ and morphisms from
$\bA$ to $\bI$ given by the formula:
      $$ \ss{Sec}(\zz_A) \ni \zf \rightarrow \widetilde{\zf}
\colon  \widetilde{\zf}(a) = a - \zf(\zz_A(a)).
                                             \tag \label{Fa4}$$

   Let us denote by $A^\#$ the set of morphisms $\zf\colon \bA\rightarrow \bI$.
 The induced from $\ss{AffSec}\zz_\bA $ affine structure on $A^\#$  (see
\Ref{Fa2}) is modeled on the vector space of affine functions on
$\underline{A}$ and is given by
      $$ (\zf -\zf')(\underline{a}) = \zf'(a) - \zf(a),
                                             \tag \label{Fa2a}$$
where $\underline{a} = \zz_\bA(a)$.
   It follows that the dimensions of $\sV(A)$ and $\sV(A^\#)$ are equal.

   The special affine space $\bA^\# =(A^\#, \bold 1_{\underline{A}})$,
where $\bold 1_{\underline{A}}$  denotes the constant function equal $1$,
will be called the {\it special affine dual} to the special affine space $\bA$.

      \claim \c{e}{Example}{}{\rm                           \label{Ca1}
   Let $\underline{A}$ be an affine space. We put $A= \underline{A}\times \R$
with $\sV(A) = \sV(\underline{A})\times \R$ and $\bA = (A,(0,1))$. The dual space
$A^\#$ is identified with the vector space $\underline{A}^\dag$ of affine
functions on $\underline{A}$.

With this identification we have $\bA^\# = (\underline{A}^\dag, \bold 1_{\underline{A}})$.
      }\endclaim

\noindent {\bf Remark.} If we take $\bA$ as in the example, then
$\overline{\bA}=(\overline{A},(0,1))$ and $\overline{A} =
\overline{(\underline{A})}\times \overline{\R}$. Here $\overline{\R}$
denotes the affine space of real numbers with the affine structure dual to
the canonical one. The multiplication by $-1$ gives an isomorphism of $\R$
and $\overline{\R}$ with the identity as the linear part, so we can identify $\overline{\bA}$ with
$(\overline{\underline{A}}\times \R, (0,-1))$.

         \claim \c{p}{Proposition}{}                     \label{Ca2}
   Let $\bA =(A,v_1)$ and $\bB =(B,w_1)$ be special affine spaces. There are the
following canonical isomorphisms:
         $$\align
   \overline{\overline{\bA}}&= \bA  \\
         \overline{\bA\ot \bB}&= \overline{\bA}\ot \overline{\bB}  \\
         \overline{\bA^\#}&= (\overline{\bA})^\# \\
      (\bA\ot \bB)^\#&= \bA^\#\ot \bB^\# .
                                          \tag \label{Fa6}\endalign$$
      \endclaim
      \proof
   The first two identities are obvious.  An affine section of $\zz_\bA$
is also an affine section of $\zz_{\overline{\bA}}$ but the affine
structures are different. In terms of morphism, $\zf\in A^\#$ gives the same
section as $-\zf \in (\overline{A})^\#$.

   Now, for each pair $\zf\in A^\#, \zc\in B^\#$ we define a mapping
   $$\zf+\zc\colon A\times B\rightarrow \R\colon (a,b)\mapsto \zf(a) + \zc(b).$$
   Since
   $$\zf(a+rv) + \zc(b-rw)= \zf(a) + \zc(b),$$
   the mapping $\zf + \zc$ projects to a morphism of special affine spaces
$\zf\oplus\zc\colon \bA\ot \bB \rightarrow \bI$. On the other hand,
$\zf\oplus\zc =\zf'\oplus\zc'$ if and only if $\zf(a) +\zc(b) = \zf'(a)
+\zc'(b)$ for each $a,b$, i.e. if and only if $\zf - \zf' = \zc' - \zc \in
\R$. It is equivalent to say that the pairs  $(\zf,\zc)$ and $(\zf',\zc')$
define the same element of $\bA^\#\ot \bB^\#$.
      \endproof

   We have also  canonical identifications $\overline \bI =\bI$, $\bI\ot \bI
= \bI$ and $\bI^\#=\bI$ given by mappings
      $$\align
   \R\rightarrow \R&\colon r\mapsto -r,  \\
   \R\ot \R\rightarrow \R&\colon [(r,s)] \mapsto r+s \\
   \R^\# \rightarrow \R&\colon  \zf \mapsto -\zf(0)  .
                                         \tag \label{Fa7}\endalign$$

   Let $\zC\colon \bA\rightarrow \bB$ be a morphism of special affine spaces
$\bA =(A,v)$ and $\bB =(B,w)$. The set
      $$G= \{ A\times B \ni (a,b) \colon b=\zC (b) \}
                                             \tag \label{Fa8}$$
   is invariant with respect to the $\R$-action $((a,b),r)\mapsto (a+rv, a+rw)$
and consequently, it is uniquely determined by its image in $A\ot \overline{B}$,
which is the image of a section of $A\ot B \rightarrow \underline{A}\times
\underline{B}$, defined on the graph $\ss{gr}(\underline{\zC})$ of $\underline{\zC}\colon
\underline{A}\rightarrow \underline{B}$.  The corresponding to this section  morphism
of special affine bundles
      $$ \ss{gr}(\zC)\colon (\zz_A\ot \zz_B)^{-1}(\ss{gr}(\underline{\zC}))
\rightarrow I
                                             \tag \label{Fa9}$$
will be called the {\it graph} of $\zC$.

   Let $\zF\colon A\times B \rightarrow \R $ be a bi-affine mapping such that for
each $a\in A$ and each $b\in B$ the corresponding mappings
      $$\align
   \zF(a,\cdot)&\colon B\rightarrow  \R \colon b'\mapsto \zF(a,b') \\
   \zF(\cdot,b)&\colon A\rightarrow  \R \colon a'\mapsto \zF(a',b)
                                       \tag \label{Fa10}\endalign$$
are morphisms of special affine spaces. It follows that
      $$ \zF(a,b+w) = \zF(a+v,b),
                                             \tag \label{Fa11}$$
   i.e. $\zF$ projects to a mapping $\widetilde{\zF}\colon A\ot B\rightarrow
\R$. We say that $\zF$ is special bi-affine.
    On the other hand, we say that a mapping
      $$\widetilde{\zF} \colon \bA\ot \bB\rightarrow \bI$$
is  special bi-affine if its pull-back to a function $\zF\colon A\times
B\rightarrow I$ is  special bi-affine.
   A special bi-affine mapping $\zF\colon  A\times B\rightarrow \R$
induces mappings
      $$\align
   \zF_l&\colon A \rightarrow B^\# \\
   &\colon a  \mapsto \zF(a,\cdot )
                                       \tag \label{Fa12}\endalign$$
   and
      $$\align
   \zF_r&\colon B\rightarrow A^\# \\
   &\colon b \mapsto \zF(\cdot ,b).
                                       \tag \label{Fa13}\endalign$$

      For $r\in \R$ we have
   $\zF_l(a+rv)(b) = \zF_l(a)(b+rw) = \zF_l(a)(b) +r\bold 1 (a)$ and,
according to \Ref{Fa2},  $\zF_l(a+rv) - \zF_l(a) = -r\bold 1$. It follows
that $\zF_l$ is a morphism of special affine spaces
   $$\zF_l \colon \overline{\bA} \rightarrow \bB^\# $$
and similarly,
   $$\zF_r \colon \overline{\bB} \rightarrow \bA^\#.$$
   We say that $\zF$ is {\it nondegenerate} if $\zF_l, \zF_r$ are injective. For
injective $\zF_l, \zF_r$, we have
   $$\dim \sV(A) \leqs \dim \sV(B^\#) = \dim \sV(B) \leqs \dim \sV(A), $$
i.e. they are isomorphisms.
   A {\it pairing} between special affine spaces $\bA=(A,v_1)$ and $\bB =
(B,w_1)$  is a nondegenerate special bi-affine mapping
      $$ \zF\colon  \bA \ot \overline \bB \rightarrow \bI.
                                             \tag \label{Fa14}$$

      \claim \c{e}{Example}{}{\rm                           \label{Ca3}
   Let $\bA=(A,v)$ be a special affine space. The mapping
   $$\align
   \zD_\bA&\colon A^\#\ot \overline A\rightarrow \bI  \\
   &(\zf,a)\colon  \mapsto -\zf(a)
                                       \tag \label{Fa15}\endalign$$
is the canonical pairing between $\bA^\#$ and $\bA$, while the mapping
      $$\align
   \overline{\zD}_\bA&\colon A\ot \overline{A}^\#\rightarrow \bI  \\
   &(a,\zf)\colon  \mapsto \zf(a)
                                       \tag \label{Fa16}\endalign$$ is the
canonical pairing between $\bA$ and $\bA^\#$. Since these pairings are
non-degenerate, we have the canonical isomorphism of $\bA$ and
$(\bA^\#)^\#$.
      }\endclaim

   For a morphism $\zF\colon \bA\rightarrow \bB$  of special affine
spaces, we define  the {\it dual morphism} $\zF^\#\colon \bB^\#\rightarrow
\bA^\#$ by the formula
   $$\zF^\#(\zc) = \zc\circ \zF.$$

      \claim \c{p}{Proposition}{}                              \label{Ca4}
   Let $\zF\colon \bA\rightarrow \bB$ be a morphism of special affine spaces. Then
         $$ (\zF^\#)^\# =\zF,$$
   where we have used the canonical isomorphisms $(\bA^\#)^\# = \bA$ and
$(\bB^\#)^\# = \bB$.
      \endclaim
      \proof
   Let $\zq\in (A^\#)^\#$ and $a\in A$ be related by the canonical
isomorphism, i.e. let $\zq(\zf) = -\zf(a)$ for each $\zf \in A^\# $. By
the definition of the dual morphism, we have for  each $\zc\in B^\#$
      $$\align
   (\zF^{\#\#}(\zq))(\zc)&=\zq(\zF^\#(\zc))  \\
         &=-(\zF^\#(\zc))(a)  \\
         &= -\zc(\zF(a)),
                                          \tag \label{Fa17}\endalign$$
what implies that $\zF^{\#\#}(\zq)$ and $\zF(a)$ are related by the canonical
isomorpism of $\bB$ and $\bB^{\#\#}$.
      \endproof

   For $\bB = \bA^\#$ we have $\zF^\# \colon \bA^{\#\#} = \bA \rightarrow
\bA^\#$.
      \claim \c{d}{Definition}{}{\rm                        \label{Ca5}
   A morphism $\zF\colon \bA \rightarrow \bA^\#$ is called {\it self-dual } if $\zF =
\zF^{\#}$.
      }\endclaim

   As we have already noticed, the model vector space for $A^\#$ is the vector space of affine functions on
$\underline{A}$. The distiguished function is the constant unit function. It
follows that the model vector space for $\underline{A^\#}$ is the space
$\sV(A)^{\textstyle *} $. The canonical projection $A^\#\rightarrow \underline{A^\#}$
corresponds to extracting from an affine function its linear part.
      \claim \c{p}{Proposition}{}                           \label{Ca6}
   Let  $\zF\colon \bA\rightarrow \bB$ be a morphism of special affine spaces.
Then $\sV(\zF^\#) =\underline{\zF}^{\textstyle *} $, where $\sV(\zF^\#)$ is
the linear part of $\zF^\#$.
      \endclaim
      \proof
   Let $a\in A$ and let $f$ be an affine function on $\underline{B}$. Then
for $\underline{a}= \zz_\bA(a)$ and $\zc\in B^\#$
      $$\align
   (\sV(\zF^\#)(f))(\underline{a})&= (\zF^\#(\zc +f) - \zF^\#(\zc))(a)    \\
         &=(\zF^\#(\zc))(a)- (\zF^\#(\zc +f))(a)  \\
         &= \zc(\zF(a)) - \left(\zc(\zF(a)) -
f(\underline{\zF}(\underline{a})) \right) \\
         &= f(\underline{\zF}(\underline{a})) =
(\underline{\zF}^{\textstyle *} f)(a)
                                       \tag \label{Fa18}\endalign$$
      \endproof

      \claim \c{c}{Corollary}{}                           \label{Ca7}
   Let  $\zF\colon \bA\rightarrow \bA^\#$ be a self-dual  morphism, then the
linear part of the induced mapping of affine spaces $\underline{\zF} \colon
\underline{A} \rightarrow \underline{A^\#}$ is skew self-adjoint.
      \endclaim
   \proof
   We have from the previous proposition that
   $$\sV(\zF^\#) = \sV(\zF)^{\textstyle *}\colon \sV(\zF^{\#\#}) \rightarrow
\sV(\zF^\#). $$
   We identify the model space $\sV(A^\#)$ with the space
$\ss{Aff}(\underline{A},\R)$ of affine functions on $\underline{A}$.
Similarly, the space $\sV(A^{\#\#})$ is identified with the space
$\ss{Aff}(\underline{A^\#},\R)$. Using the pairing~\Ref{Fa15} (or~\Ref{Fa16})
we see that an element of $\sV(A)$ can be interpreted as an affine
function on $\overline{\underline{A^\#}}$. With this interpretation the
sets $\sV(A)$ and $\sV(A^{\#\#})$ are equal and the linear part of the
canonical isomorphism $A\rightarrow A^{\#\#}$  is the identity. Let $\zf$
be an affine function on $\underline{A^\#}$ with the linear part $\sv(\zf)
\colon \sV(\underline{A^\#}) = \sV(A)^{\textstyle *} \rightarrow \R$. The
same function, interpreted as an affine function on
$\overline{\underline{A^\#}}$ has the linear part $-\sv(\zf)$. It follows,
that the canonical isomorphism of $\bA$ and $\bA^{\#\#}$ projects to the
minus standard isomorphism of $\sV(A)$ and $\sV(A)^{\textstyle{* *}} $.

   The linear part of a self-dual morphisms projects then to a skew self-adjoint
linear mapping.
   \endproof

      \claim \c{p}{Proposition}{}                              \label{Ca8}
   Let $\zF\colon \bA\rightarrow \bA^\#$ be a morphism of special affine spaces.
Then $\zF$ is self-dual if and only if for each pair $a,b\in A$
      $$ \zF(b)(a) = -\zF(a)(b) .
                                                \tag \label{Fa19}$$
       \endclaim
      \proof
   Let $\zc\in A^{\#\#}$ and  $b\in A$. By the definition of the dual morphism
   $\zF^\#(\zc)(b) = \zc(\zF(b))$. If $\zc$ and $a\in A$ are related by  the
canonical isomorphism $\bA\rightarrow \bA^{\#\#}$ then $\zc(\zF(b)) =
-\zF(b)(a)$ and we get the formula
        $$ \zF^\#(a)(b)= - \zF(b)(a).
                                                            \tag \label{Fa19b}$$
   It follows that $\zF$ is self-dual if and only if $\zF(a)(b) = -\zF(b)(a)$ for
each pair $a,b\in A$.
   \endproof
   Let us choose a point $a_0\in A$ and $\zf_0\in A^\#$. We define an
isomorphism $A\rightarrow V\times \R$, where $V =\sV(\underline{A})$, by

      $$ A\ni a\mapsto (\sv(\zz_A)(a-a_0), \zf_0(a)) \in V\times \R,
                                                \tag \label{Fa20}$$
   where $\zz_A\colon A\rightarrow \underline{A}$ is the canonical projection.
   The affine structure in $V\times \R$ induced by the affine structure of
$A$ coincides with the canonical affine structure of a vector space and
the distinguished vector is $(0,1)$.

   Now, let $\zf\in A^\#$. We have
      $$\align
   \zf(a)&= \zf_0(a) -(\zf -\zf_0)(\zz_A(a))  \\
         &= \zf_0(a) -\sv(\zf-\zf_0)(a-a_0) -(-\zf(a_0)).
                                       \tag \label{Fa21}\endalign$$

   The mapping

      $$ A^\#\ni \zf \mapsto (\sv(\zf-\zf_0), -\zf(a_0)) \in  V^{\textstyle *} \times \R
                                             \tag \label{Fa22}$$
   gives an isomorphism of special affine spaces. The affine structure in
$V^{\textstyle *} \times \R$ induced by the affine structure of $A^\#$ is
the canonical affine structure of a vector space and the distinguished
vector is $(0,1)$.
  With the isomorphisms \Ref{Fa20}, \Ref{Fa22} the pairing $\zD_\bA$ reads
      $$ \zD_\bA\colon ((f,t),(v,r))\mapsto t -r +\langle f,v\rangle .
                                                      \tag \label{Fa23}$$

\noindent{\bf Remark} We must be aware of the fact that in the formula~\Ref{Fa23}
the pair $(v,r)$ represents an element of $\overline{A}$ and the affine structure
in $\overline{A}$ does not correspond to the canonical affine structure on the
vector space $V\times \R$ but to the conjugate one. This implies that for given
$(f,t)$ the linear part of the mapping $(v,r)\mapsto t-r + \langle f,v\rangle $
is $r-\langle f,v\rangle $.

   \claim \c{t}{Theorem}{}                              \label{Ca9}
   Let $\underline{\zF}\colon \underline{A}\rightarrow \underline{A^\#} $ be an
affine mapping such that its linear part is skew self-adjoint. There exists exactly
one lift of $\underline{\zF}$ to a self-dual morphism $\zF\colon \bA\rightarrow \bA^\#$.
      \endclaim
      \proof
   Let $a_0\in A$ and $\zf\in A^\#$ and let $\zF\colon \bA\rightarrow \bA^\#$ be a morphism.
   With the isomorphisms \Ref{Fa20} and \Ref{Fa22} we represent $\zF$ by a quadruple
$(F,f,g,t)$, where $F\in \ss{Lin}(V,V^{\textstyle *} )$, $f,g\in V^\*, t\in \R$:
      $$ \zF(v,r) = (F(v)+f, g(v)+t+r ).
                                                      \tag \label{Fa24}$$
   $\zF$ is self dual if and only if \Ref{Fa19} is satisfied, i.e.
   $$\zD_A((F(v)+f,g(v) +t +r),(w,s)) = -\zD_A((F(w)+f,g(w)+t +s),(v,r))$$
   for each pairs $(v,r),(w,s)$. Using the formula \Ref{Fa23} we obtain
   $$\langle F(v)+f,w\rangle +g(v) +t +r -s = - \langle F(w) +f,v\rangle -g(w) -t
-s +r.$$
   It follows that $\zF$ is self-dual if and only if $F = -F^{\textstyle *} $,
$f=-g$ and $t=0$. This implies that a self-dual morphism is uniquely determined by
$F$ and $f$, i.e. by the mapping $\underline{\zF}\colon \underline{A}\rightarrow
\underline{A^\#} $ with the skew self-adjoint linear part.
         \endproof
    \hl2{Special affine bundles}
   A {\it special vector bundle} is a vector bundle with a distinguished
non-vanishing section and a {\it special affine bundle} is an affine bundle
modeled on a special affine space. A fibre of a special vector bundle is a
special vector space and a fibre of a special afffine bundle is a special affine
space.
   Let $\bA= (A,X)$ and $\bB=(B,Y)$ be special affine bundles. We define in
an obvious way special affine bundles $\overline{\bA}$, $\bA\ot \bB$,
$\bA^\#$ and $\bA\ot _M\bB$ if $A, B$ have the same base manifold $M$. We
have also and affine bundle $\underline{A} = A/\{X\}$ and the canonical
projection $\zz_\bA\colon A\rightarrow \underline{A}$.

    Let $\bV =(V,X)$ and $\bV' = (V',X')$ be special vector bundles. A
morphism  $F\colon V\rightarrow V'$ of vector bundles is called a {\it
morphism of special vector bundles} if
$X$ and $X'$ are $F$-related.
    A {\it morphism of special affine bundles} is a morphism of affine bundles
such that its linear part is a morphism of special vector bundles.
        \claim \c{e}{Example}{}{\rm                         \label{Ca9a}
Let $\zz \colon Z \rightarrow M$ be an affine fibration modeled on the
trivial fibration $pr_M\colon M \times \R \rightarrow M$. A pair $\bZ= (Z,\bold
1_M)$ is a special affine bundle.
Also the pair $\sC\bZ = (\sC Z,(0,\bold 1_M))$ is a special affine
fibration with $\underline{\sC Z} = \sP Z$.
        }\endclaim
      \claim \c{e}{Example}{}{\rm                           \label{Ca10}
   Let  $pr_Z\colon Y=Z\times \R \rightarrow Z$ be a trivial affine bundle and
let $\bY =(Y,\bold 1_Z)$. We define a morphism of special affine bundles
      $$\align
   \zl_Z&\colon \bY\rightarrow \bZ  \\
   &\colon (z,r) \mapsto z+r
                                          \tag \label{Fa26}\endalign$$
   which gives a mapping $\zl_Z^{\textstyle *} \colon \ss{Sec}(Z) \rightarrow
\ss{Sec}(Y) = \Cal C^\infty(Z)$. We have
      $$ \zl_Z^{\textstyle *} \zs(z) = \zs(\zz(z)) -z.
                                                \tag \label{Fa27}$$
   In a trivialization provided by a section, we have
      $$ \zl_Z^{\textstyle *} \zs(m,t) = \zs(m) -t.
                                                \tag \label{Fa28}$$
   The image of the induced relation $\sP \zl_Z$ is a coisotropic submanifold
$K_{-1}$ of $\sT^{\textstyle *} Z$,
   $$K_{-1}=\{a\colon \langle a,X_1\rangle =-1\},$$
where $X_1$ is a vector field on $Z$ represented by the mapping $Z\times \R\ni
(z,r) \rightarrow z+r \in Z$.
      }\endclaim
      \claim \c{p}{Proposition}{}                        \label{Ca11}
   Let $\bZ=(Z,\bold 1_M)$ and  $\bZ'=(Z',\bold 1_{M'})$ be special affine
bundles modeled on $M\times \R$ and $M'\times \R$ respectively. There are
canonical symplectomorphisms
      $$\align
   (\sP(\bZ\ot \bZ'),\ezw_{\bZ\ot \bZ'}) &= (\sP\bZ\times\sP\bZ',
\ezw_\bZ + \ezw_{\bZ'})   \\
       (\sP\overline{Z}, \ezw_{\overline{\bZ}})   &= (\sP Z, -\ezw_\bZ).
                                       \tag \label{Fa31}\endalign$$
      \endclaim
      \proof

    (1)  \ let $\zs,\zc$ be  sections of $Z$ and $\zs', \zc'$ sections of $Z'$.
   We have from the definition of $Z\ot Z'$ that
   $$\zs\ot \zs' - \zc\ot \zc' = (\zs-\zc) +(\zs' - \zc')).$$
   It follows that $(\zs\ot \zs', (m,m'))$ is equivalent to $(\zc\ot
\zc',(n,n'))$ if and only if $(\zs,m)$ is equivalent to $(\zc,n)$ and $(\zs',m')$ is
equivalent to $(\zc',n')$. The correspondence $\zs\ot \zs' \rightarrow
(\xd \zs(m),\xd\zs'(m'))$ gives a diffeomorphism of manifolds $\sP(Z\ot Z')$
and $\sP Z\times \sP Z'$. It is easy task to check that this isomorphism is also
a symplectomorphism.

   \indent (2) \ It is obvious that $\sP Z$ and $\sP \overline{Z}$ are equal
as manifolds. Let $\zs$ be a section of $Z$. The same mapping interpreted as a
section of  $\overline{Z}$ will be denoted by $\overline{\zs}$. Since $\zs -\zs'
= \overline{\zs}' - \overline{\zs}$, the isomorphisms $I_{\xd\zs}\colon \sP Z
\rightarrow \sT^{\textstyle *} M$ and $I_{\xd\overline{\zs}}\colon \sP
\overline{Z}\rightarrow \sT^{\textstyle *} M $ are related by $I_{\xd\zs} = -
I_{\xd\overline{\zs}}$. It follows that
   $$\ezw_{\overline{\bZ}} = I_{\xd\overline{\zs}}^{\textstyle *} \ezw_M = -
I_{\xd\zs}^{\textstyle *} \ezw_M =-\ezw_{\bZ}.$$
      \endproof
 We have also similar equalities for contact fibrations:
      $$\align
    \sC \overline{Z}&= \overline{\sC Z}  \\
         \sC(Z\ot Z')&= \sC Z\ot \sC Z'.
                                                \tag \label{Fa34}\endalign$$
        \hl1{Reduced tangent bundles}
   Let $\bZ=(Z,\bold 1_M)$ be, as before, a special affine bundle modeled
on the trivial bundle $M\times \R$. We denote by $\zvf$ the   $\R$-action on $Z$:
$\zvf(z,r) = z+r$. The $\R$-action $\zvf$ induces an $\R$-action $\zvf_{\textstyle *}
$ on $\sT Z$. In this section we provide an interpretation of special affine bundles,
dual to $\sP Z\times \R $ and $\sC \bZ$ in terms of vectors tangent to $Z$ and the
action $\zvf$.

   We will make use of the correspondence~\Ref{Fa27}  between sections of
$\zz\colon Z\rightarrow M$ and functions on $Z$.  For $\zs\in
\ss{Sec}(\zz)$ we denote by  $f_\zs$ the corresponding function on $Z$.
For each $v\in \sT_z Z$ we have a mapping
   $$\Cal C^\infty (Z)\ni f\mapsto \langle \xd_z f,v\rangle\in \R,$$
   and the induced mapping
   $$\ss{Sec}(\zz) \ni \zs \mapsto \langle \xd_z f_\zs, v\rangle \in \R.
                                 \tag \label{Fa35}$$

   In a trivialization $f_\zs(m,s) = \zs(m) - s$, $v = \underline{v} +
\dot s \dfrac{\partial }{\partial s}$, where $\underline{v} \in
\sT_{\zz(z)} M$, and the formula~\Ref{Fa35} takes the form
      $$ \zs \mapsto \langle \xd_m \zs, \underline{v}\rangle -\dot s.
                                             \tag \label{Fa36}$$
   This shows that $v$ defines an affine function on $\sP_{\zz(z)} Z $ and
the vector $\dfrac{\partial }{\partial s} $ gives the function $-\bold
1_{\sP Z}$. Two vectors $v,v'$ define the same function on $\sP_{\zz(z)}
Z$ if and only if they are in the same orbit of the $\R$-action
$\zvf_{\textstyle *} $. We will denote by $\widetilde{\sT}Z$ the space of
orbits of this action  and we will call it the {\it reduced tangent
bundle}. It is a vector bundle over $M$ and there are a canonical
projection $\tT\zz \colon \tT Z \rightarrow \sT M$ induced by the tangent
projection $\sT \zz \colon \sT Z \rightarrow \sT M$ and the reduced
canonical projection $\widetilde{\ezt}_\bZ \colon \tT Z \rightarrow M$.

   This way we have obtained a bi-affine mapping
      $$ \zF\colon  \sP Z\times _M \tT Z \rightarrow \R \colon
(\xd_{\zz(z)}\zs, [v]) \mapsto \langle \xd_z f_\zs, v\rangle \in \R
                                          \tag \label{Fa37}$$
   which extends to a special bi-affine mapping (denoted by the same letter)
   $$ \zF\colon  (\sP Z\times \R, (0,1))\ot _M (\tT Z,-X_1) \rightarrow \R ,
                                          \tag \label{Fa38}$$
   where $-X_1$ is the fundamental vector field for the action $\zvf_{\textstyle
*} $, i.e. in local trivialization $X_1 = \dfrac{\partial }{\partial s} $.

   The induced mapping~\Ref{Fa13} $\zF_r\colon \tT Z \rightarrow (\sP
Z\times \R)^\# = \sP Z^\dag$ (see Example~\Ref{Ca1}) is bijective and
defines an isomorphism of special affine bundles
      $$\zF_r\colon \overline{(\tT Z, -X_1)} \rightarrow (\sP Z\times
\R, (0,1))^\# = (\sP Z^\dag, \bold 1_{\sP Z}).
                                             \tag \label{Fa39}$$
   Since $\tT Z$ is a vector bundle, we have an involution $v\mapsto -v$, which
is an isomorphisms of special vector (affine) spaces
      $$ \overline{(\tT Z, -X_1)} \simeq (\tT Z,X_1).
                                             \tag \label{Fa40}$$
   In the following we identify the special affine dual $(\sP Z\times \R,
(0,1))^\#$ with $\tT \bZ = (\tT Z, X_1)$.
   In a trivialization $I_\zs$ provided by a section $\zs_0$, we have $Z
\simeq M\times \R$, $\sP Z \simeq \sT ^{\textstyle *} M$ and $\tT Z \simeq
\sT M\times \R$, with
      $$ (p,0)(\underline{v},\dot s) =\dot s- \langle p,\underline{v}\rangle .
                                             \tag \label{Fa41}$$

    \noindent {\bf Remark.} The trivialization $\tT Z\simeq \sT M\times \R$
is given by the decomposition of  a vector tangent to $Z$:
        $$ \sT Z\ni v = \underline{v} + \dot s \frac{\partial }{\partial s}
                                                            \tag \label{Fa41b}$$
     and it is not compatible with the special affine structure $(\tT Z,
-X_1)$. A trivializaton $v\mapsto (\underline{v},r)= (\underline{v},-\dot
s) $ is compatible and gives
        $$ (\underline{v},r) (p,0) = -(p,0)(\underline{v},r) = r + \langle
p, \underline{v}\rangle
                                                            \tag \label{Fa41c}$$
    (see~\Ref{Fa15} and~\Ref{Fa24}).

    \noindent {\bf Remark.} We can identify the manifolds $\sP Z$ and $\sP
\overline{Z}$ and consequently,  bundles $(\sP Z)^\dag$ and $(\sP
\overline{Z})^\dag$ of affine functions. However, if $v\in \sT Z$ defines a
function $g_v$ on $\sP Z$,  then the function defined by $v$ on $\sP
\overline{Z}$ iz $-g_v$. The natural isomorphism of $\sP^\dag Z$ and
$\sP^\dag \overline{Z}$ corresponds to the multiplication by $-1$  in $\tT
Z$.

    For each $v\in \sT_m M$ there is the unique element $\widetilde{v}\in
\tT_mZ$ such that $\tT\zz(\widetilde{v}) =v$ and $(p,0)(\widetilde{v})
=0$. In the trivial case this means that $r= \langle p, v\rangle $. This
justifies a notation
   $$\widetilde{v} =  \langle p,v\rangle  \tag\label{Fa42}$$
   which will be used in the following.

         \claim \c{p}{Proposition}{}                        \label{Ca12}
      Let $\bZ=(Z,\bold 1_M)$ and  $\bZ'=(Z',\bold 1_{M'})$ be special affine
bundles modeled on $M\times \R$ and $M'\times \R$ respectively. There are
canonical isomorphisms
      $$\align
   \tT I&\simeq I,  \\
         \tT (Z\ot Z')&\simeq \tT Z\ot \tT Z'.
                                       \tag \label{Fa43}\endalign$$
      \endclaim
      \proof
   $\sT I$ = $\R\times \R$ and $\zvf_{\textstyle *} ((s,\dot s),r)= (s+r,\dot s)
$. It follows that the space of orbits is parametrized by $\dot s$ and the
distinguished vector is $1$.


   The second isomorphism follows from Proposition~\Ref{Ca1} and
Proposition~\Ref{Ca11}:
   $$ \tT (Z\ot Z') = (\sP(Z\ot Z'))^\dag = (\sP Z \times \sP Z')^\dag
= \sP^\dag Z\ot \sP^\dag Z' = \tT Z\ot \tT Z'. $$
      \endproof

   Now, we  find the dual to the contact bundle $\sC \bZ$. As before, we
associate to a vector $v \in \sT_z Z$ a mapping  $\Cal C^\infty(Z) \rightarrow
\R$, this time  by the formula
      $$ f\mapsto \langle \xd_zf,v\rangle +f
                                             \tag \label{Fa44}$$
   which gives for $f=f_\zs$, in a trivialization provided by a section of $\zz$,
      $$ \ss{Sec}(\zz)\ni \zs \mapsto \langle \xd_m \zs, \underline{v}\rangle
-\dot s + \zs(m) -s .
                                             \tag \label{Fa45}$$
   This mapping projects to a special  affine mapping $\sC_m Z \rightarrow \R $ and
this way a vector $v\in \sT_z Z$ defines an element of $\sC^\#_m Z=(\sC_m Z)^\#$,
where $m =\zz(z)$. Thus we have a mapping $\sT Z\rightarrow \sC^\# Z$. It follows
from the formula~\Ref{Fa45} that  this mapping is surjective and when restricted to
$\sT_zZ$ it is bijective.  Two vectors $v\in \sT_zZ, v'\in \sT_{z'}Z$ define the
same element of $\sC^\# Z$ if $\zz(z) = \zz(z')$ and
   $$  v' = \zvf_{\textstyle *}(v,r) -rX_1(z'),
                                             \tag \label{Fa46}$$
   i.e. they are in the same orbit of an $\R$-action on $\sT Z$:
      $$ (v,r) \mapsto \zvf_{\textstyle *}(v,r) -rX_1(z+r).
                                             \tag \label{Fa47}$$
   The space of orbits of this action will be denoted by $\bT Z$. The vector
bundle structure of $\sT Z$ induces a structure of an affine bundle on $\bT Z$
with the model vector  bundle $\tT Z$. We denote by $\bT\zz$ and
$\overline{\ezt}_\bZ$ the induced projections $ \bT Z\rightarrow \sT M$ and $\bT
Z \rightarrow M$.

   The mapping
      $$((\zs,m),v)\mapsto \langle \xd_z f_\zs,v\rangle +f_\zs(z), \ \ \zz(z) =m$$
   projects then  to a bi-affine mapping
   $$\sC Z\times_M \bT Z \rightarrow \R$$
   and a special bi-affine and nondegenerate mapping
   $$\sC \bZ\ot_M (\bT Z,-X_1) \rightarrow I $$
   which gives an isomorphism
      $$ \sC^\#\bZ \simeq \overline{(\bT Z,-X_1)}.
                                             \tag \label{Fa48}$$

   The bundle $\overline{(\bT Z,-X_1)}$ is, by the definition, equal
$(\overline{\bT Z}, -X_1)$. This means that we have on $\bT Z$ an affine structure
induced from the affine structure on $\sT Z$, conjugate to the canonical one.
   Using the involution  $\widetilde{v}\rightarrow -\widetilde{v}$ on $\tT Z$,
we get  $(\overline{\bT Z},-X_1) = (\bT Z, X_1)$ and the  isomorphism
      $$  \sC^\#\bZ \simeq (\bT Z,X_1) = \bT \bZ.
                                             \tag \label{Fa49}$$
   This isomorphism  can be deduced also from Proposition~1, an obvious isomorphism
$\sC \bZ = (\sP Z, \bold 1_{\sP Z})\ot _MZ$, and the canonical isomorphism
$Z^\# = Z$. The last isomorphism is provided by the pairing
      $$ \bZ\ot_M \overline{\bZ}\ni z\ot z'\mapsto z- z'.
                                             \tag \label{Fa50}$$
   \noindent{\bf Remark.} The involution $v\mapsto -v$ on $\sT Z$ gives
correspondence between the $\R$-action~\Ref{Fa47} and an $\R$-action
given by the formula
      $$ (v,r) \mapsto \zvf_{\textstyle *}(v,r) +rX_1(z+r).
                                              \tag \label{Fa51}$$
   Let us denote by $\bT^cZ$ the manifold of orbits of the action~\Ref{Fa51}. It
is an affine bundle over $M$ with the affine structure induced by the vector
bundle structure on $\sT Z$ (over $Z$) and isomorphic to the affine bundle $\bT
Z$. The model bundle for $\bT^c Z$ is $\tT Z$.  This way we have obtained, via
the involution in $\sT Z$, an isomorphism of special affine bundles
         $$\overline{\bT\bZ} \simeq (\bT^c Z, -X_1).  \tag\label{Fa52}$$

   For every bundle $Z$ we have a canonical projection $\zm\colon \sC Z
\rightarrow Z$ which is a surjective morphism of special affine bundles $\zz_{\sC
\bZ}\colon  \sC Z \rightarrow \sP Z$ and $\zz \colon Z\rightarrow M$. There is
also the dual injective morphism
      $$ \bZ= \bZ^\# \rightarrow \sC^\#\bZ = \bT \bZ
                                             \tag \label{Fa53}$$
   which is the composition of the zero section $Z\rightarrow \sT Z$ and
the canonical projection $\sT Z\rightarrow \bT Z$.

         \claim \c{p}{Proposition}{}                        \label{Ca13}
      Let $\bZ=(Z,\bold 1_M)$ and  $\bZ'=(Z',\bold 1_{M'})$ be special affine
bundles modeled on $M\times \R$ and $M'\times \R$ respectively. There are
canonical isomorphisms
      $$\align
   \bT I&\simeq I,  \\
         \bT \overline{\bZ}&\simeq \overline{\bT \bZ },  \\
         \bT (Z\ot Z')&\simeq \bT Z\ot \bT Z'.
                                       \tag \label{Fa54}\endalign$$
      \endclaim
   \proof
   From~\Ref{Fa53} we  have the canonical injection $I\rightarrow \bT I$ which is
an isomorphism.
   Two remaining isomorphisms follow by the duality from~\Ref{Fa34}
and~\Ref{Fa6}. They can be also obtained directly from the definition of
$\bT Z$ as the quotient of the tangent bundle $\sT Z$.
   \endproof
   Let $\zF \colon \bZ \rightarrow \bZ'$ be a morphism of special affine
bundles. The fundamental vector fields $-X_1$ and $-X_1'$ of the canonical
$\R$-actions are $\zF_{\textstyle *} $-related, i.e. $X'_1 =
\zF_{\textstyle *} X_1 $. It follows that the tangent mapping $\sT\zF
\colon \sT Z \rightarrow \sT Z'$  projects to morphisms of special affine
bundles
      $$\align
   \tT \zF&\colon \tT \bZ \rightarrow \tT \bZ'  \\
   \bT \zF&\colon \bT \bZ \rightarrow \bT \bZ'  .
                                          \tag \label{Fa55}\endalign$$
   In particular, we can apply $\tT, \bT$ to a morphism $\zf\colon \bZ
\rightarrow I $. Such a morphism corresponds in a unique way to a section of
$\zz$. Let us recall that this correspondence is given by the condition
         $$\zf\circ \zs =0 .$$
   We denote by $\zf_\zs$ the morphism corresponding to a section $\zs$.

   \noindent {\bf Remark} A morphism $\zf_\zs\colon \bZ\rightarrow I$
defines a function $\zf_\zs \colon Z\rightarrow \R$ which is related to
the already introduced function $f_\zs$ by $\zf_\zs = -f_\zs$.

The reducet tangent morphism
   $$\tT \zf \colon \tT Z\rightarrow \tT I = I$$
and the morphism
   $$\bT \zf \colon \bT Z \rightarrow \bT I =I  $$
 (see~\Ref{Fa43}, \Ref{Fa54}) correspond to sections of $\tT\zz \colon \tT
Z\rightarrow \sT M$ and $\bT \zz\colon  \bT Z \rightarrow \sT M  $
respectively. If $\zf = \zf_\zs$ then we denote these sections by
$\xd_{\tT} \zs$ and $\dbt \zs$ respectively, and we call them {\it tangent
or complete lifts} of the section $\zs$.

   In the trivialization provided by a section of $\zz$, $\zf_\zs(m,s) = s-\zs(m)$,
$\tT Z =\sT M\times \R$, $\bT Z= \sT M\times \R$,
and
         $$\align
   \tT \zf_\zs(v,\dot s)&= \dot s - \dt \zs(v),  \\
      \dtt\zs   &= \dt \zs , \\
      \bT\zf_\zs (v,t)&= t-\zs(m) -\dt \zs(v),\\
      \dbt\zs&= \dt\zs + \zs,
                                       \tag \label{Fa56}\endalign$$
   where $v\in \sT_m M$ and $\dt$ is the standard tangent lift of a function on $M$
to a function on $\sT M$: $\dt f(v) = \langle \xd f,v\rangle $ (see [\Ref{PT}, \Ref{GU1}]).

        \hl1{Canonical isomorphisms}
   It is known (see [\Ref{KU}]) that for a vector bundle $\zx \colon E\rightarrow M$
there is a canonical isomorphism of $\sT^{\textstyle *} E$ and $\sT^{\textstyle
*} E^{\textstyle *} $. The graph of this isomorphism is a Lagrangian submanifold
of $\sT^{\textstyle *} (E^{\textstyle *}\times E )$ generated by the evaluation
function $E^{\textstyle *} \times _ME \ni (a,f) \mapsto \langle f,a\rangle$.

   In the affine case we replace $E$, or rather $E\times \R$ by a special affine
bundle $\bA = (A, X_A)$ and $E^{\textstyle *} \times \R$ by $\bA^\#$.
   We have the pairing~\Ref{Fa15}  between $\bA^\#$ and $\bA$
      $$\align
   \zD_\bA&\colon \bA^\#\ot_M \overline{\bA}\rightarrow \bI  \\
   &\colon (\zf,a)  \mapsto -\zf(a)
                                       \tag \label{Fa77}\endalign$$
   which corresponds to a section $\ezd_\bA$ of $\zz_{\bA^\#\ot
\overline{\bA}} \colon A^\#\ot_M \overline{A} \rightarrow  \underline{A^\#}
\times_M \underline{\overline{A}}$.  This section generates a Lagrangian
submanifold of $\sP(\bA^\#\ot_M \overline{\bA})$ and, together with the
canonical embedding
    $$\bA^\#\ot_M \overline{\bA} \hookrightarrow \bA^\#\ot
\overline{\bA},$$
    a Lagrangian submanifold $L$ of $\sP(\bA^\#\ot \overline{\bA}) =
\sP(\bA^\#)\times \sP(\overline{\bA})$.  The pull-back of $\ezd_\bA$ with respect
to the projection
   $$\sP(A^\#\ot \overline{A})\rightarrow \underline{A^\#}
\times \underline{\overline{A}} $$
\noindent
   is a section of
   $$\zz_{\sC(A^\#\ot\overline{A})} \colon \sC
(\bA^\#\ot\overline{\bA}) = \sC A^\#\ot \overline{\sC A} \rightarrow
\underline{A^\#} \times \underline{\overline{A}}$$
   over $L$.  This section corresponds to a morphism

      $$   \zm_\bA \colon \zz_{\sC(\bA^\#\ot\overline{\bA})}^{-1}(L) \rightarrow I.
                                             \tag \label{Fa78}$$

      \claim \c{t}{Theorem}{}                              \label{Ca15}
   The morphism $\zm_\bA$ is the graph of an isomorphism $\zC_\bA\colon
\sC \bA \rightarrow \sC \bA^\#$
   of special affine bundles. Moreover, $\zC_A$ preserves Liouville forms, i.e.
   $$\zvy_\bA = \zC^{\textstyle *} _\bA \zvy_{\bA^\#}. $$
      \endclaim
      \proof
   Let $\zc$ be a section of $\zh\colon A\rightarrow M $. It defines an
isomorphism $I_\zc\colon A\rightarrow E$ of $A$ and the model bundle $E$.
It provides also a trivialization of the bundle $\zz_{\bA^\#}\colon A^\#
\rightarrow \underline{A^\#}$. Similarly, a section $\zc^{\textstyle *} $
gives a trivialization of $\zz_\bA \colon A \rightarrow  \underline{A}$
and an isomorphism of $A^\#$ and its model space. Having choosen $\zc$ and
$\zc^{\textstyle *} $, we can identify $A$ with $\underline{E}\times \R$
and $A^\#$ with $\underline{E}^{\textstyle *} \times \R$.
   We have then
   $$ \zD_\bA((f,t),(e,r)) = t-r + \langle f, e\rangle $$
   and
   $$ \zd_\bA(f,e) = - \langle f,e\rangle .$$
   According to~\Ref{Fa31} we identify $\sP\overline{\bA}$ with $\sP \bA$ and
the symplectic structure on $\ezw_{\overline{\bA}}$ with $-\ezw_\bA$. In the
trivialization $\sP\bA = \sT^{\textstyle *} \underline{E}$, $\sP\bA^\#=
\sT^{\textstyle *} \underline{E}^{\textstyle *} $, $\sC\bA = \sT^{\textstyle *}
\underline{E}\times \R$, and $\sC\bA^\#= \sT^{\textstyle *}
\underline{E}^{\textstyle *}\times \R$.

   Let $(x^i,y^a,s)$ be an adapted coordinate system  on $\underline{E}\times \R$ and
$(x^i,f_a,t)$ the dual coordinate system on $\underline{E}^{\textstyle *}\times \R$.
Let $(x^i,y^a,p_j,\zp_a)$ be the induced coordinate system on $\sT^{\textstyle *}
\underline{E}$ and $(x^i,f_a,q_i,\zq^a)$ the induced coordinate system on
$\sT^{\textstyle *} \overline{E}^{\textstyle *} $. We have
   $$\zvy_\bA = p_i\xd x^i + \zp_a \xd y^a, \ \ -\zvy_{\bA^\#} = -q_i\xd x^i +
\zq^a\xd f_a,       \tag\label{Fa80} $$
   and the Lagrangian submanifold $L$, generated by the function $\zd_\bA = -
f_ay^a$ is given by equations
   $$\align
      y^a&= -\zq^a,  \\
      p_j   &= q_j,  \\
      \zp_b   &= f_b .
                                       \tag \label{Fa81}\endalign$$
   It follows that $L$ is the graph of a symplectomorphism $\sT^{\textstyle *}
\underline{E}\rightarrow \sT^{\textstyle *} \underline{E}^{\textstyle *} $.

   The  coordinate systems on $ A$ and $ A^\#$ give the following
local expression for $\zD_\bA$:
   $$ \zm_\bA = t-r + f_a y^a    \tag \label{Fa82}$$
   and for the corresponding section $\zd_\bA$
   $$ t-r = - f_ay^a.$$
   It follows that $\zC_\bA$ is given in the induced local coordinates by

      $$\align
      x^i\circ \zC_\bA&= x^i,\\
      \zq^a\circ \zC_\bA&= -y^a,  \\
      q_j\circ \zC_\bA&= p_j,  \\
      f_b\circ \zC_\bA&= \zp_b, \\
      t\circ \zC_\bA&= r- f_ay^a.
                                          \tag \label{Fa83}\endalign$$
   It follows from these formulae and from~\Ref{Fa80}  that
   $$ \zC_\bA^{\textstyle *}  \zvy_{\bA^\#} = p_i\xd x^i - y^a\xd \zp_a +
\xd(\zp_ay^a) = p_i\xd x^i +\zp_a\xd y^a = \zvy_\bA.$$
   \endproof


   There are canonical contact structures on the contact fibration $\sC
\zz_\bA \colon \sC A\rightarrow \underline{A}$ and on the contact fibration
$\sC\zz_{\bA^\#} \colon \sC A^\#\rightarrow \underline{A}^\#$.
   It follows from the theorem that there are also  fibrations $\sC A\rightarrow
\underline{A}^\#$  and $\sC A^\# \rightarrow \underline{A}$,  which are  special
affine bundles.
   Moreover, we have morphisms~\Ref{Fa32}  of special affine bundles:
         $$\align
      \sC \bA&\rightarrow  \bA,  \\
         \sC \bA&\rightarrow \bA^\#,  \\
         \sC \bA^\#&\rightarrow \bA^\#, \\
      \sC\bA^\#&\rightarrow \bA .
                                          \tag \label{Fa84}\endalign$$

   The projections $\sC \bA \rightarrow \bA^\#$, $\sC\bA^\#\rightarrow \bA$, and
the related structures can be obtained directly, i.e. not using the mapping
$\zC_\bA$. As an example  we give an alternative definition of  the projection $\sC\bA
\rightarrow \bA^\#$.

   First, we define a mapping
      $$ \zq\colon  A\times _MA\rightarrow \sT A
                                                   \tag \label{Fa85}$$
   where $\zq(a,b)$ is a vector represented by a curve $\zg_{a,b}\colon t\mapsto a + t(b-a)$.
   In the case under consideration the action $\zvf$ of~\Ref{Fa51} reads
   $$\zvf(a,r) = a+rX_A(\zh(a))$$
   and the mapping $(\zq(a,b),r)\rightarrow \zq(a+rX_A(\zh(a)))$ coincides with the
action~\Ref{Fa51}. Hence the mapping $\zq$ projects to
      $$ \overline{\zq}\colon \underline{A}\times _MA\rightarrow \bT A
                                                   \tag \label{Fa86}$$
   and for each $\underline{a}\in \underline{A}$ the mapping
$\overline{\zq}(\underline{a},\cdot)$ is a morphism of special affine bundles
      $$ \overline{\zq}(\underline{a},\cdot)\colon \bA \rightarrow \bT \bA .
                                                   \tag \label{Fa87}$$
   The dual morphisms
   $$\overline{\zq}(\underline{a},\cdot)^\#\colon \sC_{\underline{a}} \bA \rightarrow \bA^\#$$
 define the projection
      $$\sC\bA \rightarrow \bA^\#.$$

   \hl2{Tangent affine bundles}
   Let $\zx\colon E\rightarrow M$ be a vector bundle. It is well known that the
tangent manifold $\sT E$ carries in a natural way two different structures of a
vector bundle. One on the  canonical fibration  $\ezt_E\colon \sT E\rightarrow E$
and  the second on the tangent fibration $\sT\zx \colon \sT E\rightarrow \sT M$.

   Let $\zh\colon A\rightarrow M$ be an affine bundle modeled on $\zx\colon
E\rightarrow M$. The tangent manifold $\sT A$ is a vector bundle $\ezt_A\colon
\sT A\rightarrow A$ and it is also an affine bundle with respect to the tangent
fibration $\sT \zh\colon \sT A \rightarrow \sT M$. The model vector bundle is
$\sT\zx\colon \sT E\rightarrow \sT M$ and the affine structure is obtained
applying the tangent functor to mappings which define the affine structure of $A$.

   Let $\bA=(A,X_A)$ be a special affine bundle. The canonical projections

      $$ \tT A\rightarrow \sT \underline{A}\rightarrow \sT M
                                                   \tag \label{Fa88}$$
   give a fibration
      $$ \tT \zh\colon \tT A\rightarrow \sT M.
                                                   \tag \label{Fa89}$$
   and the affine bundle structure on $\sT A \rightarrow \sT M$ induces an
affine bundle structure on this fibration. The model vector bundle is the reduced
tangent bundle $\tT \zx \colon \tT E \rightarrow \sT M$. To avoid ambiguities we
denote by $\dot +$ the operation of addition for this structure. The section $X_A$
can be lifted to a section $\widetilde{X}_A$  of $\tT \zh$ in the following way:
   for $v\in \sT_m M$ the vector $\widetilde{X}_A(v)$ is the equivalence class of
a  vector $\widehat{X}_A(v)$ in $\sT_0(m)E$ defined by the formula
      $$  \widehat{X}_A(v) = \sv_\sT X_A(0_E(m)) +\sT 0_E(v) ,
                                                   \tag \label{Fa90}$$
where $\sv_\sT X_A$ is the standart vertical lift of $X_A$ and $0_E\colon
M\rightarrow E$ is the zero section of $\zh$.

   The pair $(\tT A, \widetilde{X}_A)$ is a special affine bundle over $\sT M$.
   Thus the reduced tangent manifold $\tT A$ carries two special affine
structures: with respect to the canonical projection $\widetilde{\ezt}_\bA\colon
\tT A\rightarrow \underline{A}$ (the distinguished vector field is the vertical
lift of $X_A$) and with respect to the reduced tangent projection $\tT \zh\colon
\tT A\rightarrow \sT M$.

   In local coordinates $(x^i,y^a,s)$ on $A$ and $(x^i,y^a,\dot x^j, \dot y^b,
\dot s)$ on $\tT A$, we have $y^a\circ X_A=0, s\circ X_A=1$ and
      $$\gathered
         y^a\circ \widetilde{X}_A=0, \ \dot y^b\circ \widetilde{X}_A=0, \ \dot s\circ
\widetilde{X}_A =1\\
         \dot x_i\circ \sv_\sT X_A = 0, \ \dot y^a\circ \sv_\sT X_A =0, \ \dot s\circ
\sv_\sT X_A = 1.
            \endgathered                           \tag \label{Fa91}$$
   It follows that the induced $\R$-actions $\tT A \ni v\mapsto v \dot +
r\widetilde{X}_A$ and $\tT A\ni v\rightarrow v + \sv_\sT X_A$ coincide and the
quotient manifold for these actions is $\sT \underline{A}$.

The special affine dual to $(\tT A, \sv_\sT X_A)$ is $\sP\bA\times \R$. We show that
the special affine dual to $(\tT A, \widetilde{X}_A)$ is $(\tT A^\#,
\widetilde{X}_{A^\#})$.

   Let us notice first that the equality  $\tT (A\ot B) = \tT A\ot \tT
B$ implies not only the equality
   $$(\tT (A\ot B), \sv_\sT X_{A\ot B})= (\tT A, \sv_\sT X_A)\ot (\tT
B, \sv_\sT X_B)$$
   which follows from Proposition~\Ref{Ca12}, but also
   $$(\tT (A\ot B), \widetilde{X}_{A\ot B}) = ((\tT A,
\widetilde{X}_A)\ot (\tT B, \widetilde{X}_B) . \tag\label{Fa92}$$
   If $A$ and $B$ have the same base manifold $M$ then~\Ref{Fa92} implies
   $$(\tT (A\ot_M B), \widetilde{X}_{A\ot B}) = ((\tT A,
\widetilde{X}_A)\ot_{\sT M} (\tT B, \widetilde{X}_B) . \tag\label{Fa93}$$

   The tangent lift $\tT\zD_A$ of the pairing $\zD_A\colon \bA^\#\ot
\overline{\bA}\rightarrow I$  is a morphism
      $$ \tT\zD_A \colon \tT (\bA^\#\ot_M \overline{\bA})\rightarrow I,
                                          \tag \label{Fa94}$$
   in local coordinates (see \Ref{Fa56} and \Ref{Fa84} )
   $$ ((x^i,y^a,\dot x^j,\dot y ^b,\dot s),(x^i,f_a,\dot x^j,\dot f_b,\dot t) )
\mapsto  \dot t -\dot s +\dot f_ay^a +f_a\dot y^a . $$
   It follows from this formula that $\tT\zD_A$ is a pairing
      $$ \tT\zD_A \colon (\tT A,\widetilde{X}_A)\ot _{\sT M} \overline{({\tT
A}^\#, \widetilde{X}_{A^\#})} \rightarrow I
                                          \tag \label{Fa95}$$
   and from~\Ref{Fa12} we get an isomorphism

      $$ (\tT A, \widetilde{X}_A)^\# \simeq(\tT A^\#, \widetilde{X}_{A^\#})  .
                                          \tag \label{Fa96}$$

   The same procedure can be applied to the tangent bundle $\bT A$ and we
get two special affine structures on $\bT A$: with respect to the
canonical projection $\bT A\rightarrow \underline{A}$ and the tangent
projection $\bT A\rightarrow \sT M$.

   \hl2{The isomorphism $\eza_\bZ$}

   Let $\bZ$ be like in the previous sections. The tangent bundles $\tT Z$ and
$\bT Z$ have been defined as reductions of $\sT Z$. This implies that iterated
tangent manifolds $\tT \tT Z$, $\tT\bT Z$ and $\bT\tT Z$ are  reductions  of the
iterated tangent manifold $\sT\sT Z$.
      \claim \c{t}{Theorem}{}                           \label{Ca16}
   The canonical flip $\ezk_Z$
$$ \CD  \sT\sT Z @() \L{\ezk_Z} @(1,0) @()
\L{\ezt_{\sT Z}}@(0,-1)& \sT\sT Z @()\L{\sT\ezt_Z} @(0,-1) \\
   \sT Z @() \L{\text{id}} @(1,0) & \sT Z
         \endCD ,                 \tag \label{Fa96a}$$
 projects to $\widetilde{\ezk}_Z$: 
$$ \CD  \tT\tT Z @() \L{\widetilde{\ezk}_Z} @(1,0) @()
\L{\widetilde{\ezt}_{\tT Z}}@(0,-1)& \tT\tT Z @()\L{\tT\widetilde{\ezt}_Z}
@(0,-1) \\
   \sT M @() \L{\text{id}} @(1,0) & \sT M
         \endCD ,                 \tag \label{Fa96b}$$
 and to $\overline{\ezk}_Z$:  
$$ \CD  \bT\tT Z @() \L{\overline{\ezk}_Z} @(1,0) @()
\L{\overline{\ezt}_{\tT Z}}@(0,-1)& \tT\bT Z @()\L{\tT\overline{\ezt}_Z}
@(0,-1) \\
   \sT M @() \L{\text{id}} @(1,0) & \sT M
         \endCD .                 \tag \label{Fa96c}$$
      \endclaim
      \proof
   The canonical $\R$-action $\zvf$ on $Z$ induces the tangent action of $\sT \R
= \R^2$ on $\sT Z$  and the iterated tangent action of $\sT\sT \R$ on $\sT\sT Z$.
Since $\zk$ is a  natural transformation of functors, we have
         $$   \sT\sT\zvf\circ  (\zk_Z\times \zk_\R) = \zk_Z\circ \sT\sT \zvf.
                                       \tag\label{Fa97}$$
   The manifold $\tT Z$ is the manifold of orbits of the kernel  subgroup of a
homomorphism  $\sT \R\rightarrow \R\colon (r,\dot r)\mapsto \dot r$.
Similarly, $\bT Z$ is the manifold of orbits of the kernel of a homomorphism
$\sT\R\rightarrow \R\colon (r,\dot r) \mapsto r+\dot r$.  The canonical
$\R$-actions are induced by these homomorphisms. It follows that $\sT \tT Z$ is
$\sT\sT Z$ reduced by the tangent group homomorphism
   $$\sT\sT \R\rightarrow \sT \R\colon (r,\dot r, r', \dot r')\mapsto (\dot
r,\dot r')$$
   and $\tT\tT Z$ is $\sT\sT Z$ reduced by the group homomorphism
   $$\zq_1\colon \sT\sT\R\rightarrow \R\colon (r,\dot r, r', \dot r')\mapsto \dot r'.$$
   In the same way we interpret $\bT\tT Z$ as $\sT\sT Z$ reduced by a group homomorphism
   $$\zq_2\colon \sT\sT\R\rightarrow \R\colon (r,\dot r, r', \dot r')\mapsto
\dot r'+ \dot r$$
   and $\tT\bT Z$ as $\sT\sT Z$ reduced by the group homomorphism
   $$\zq_3\colon \sT\sT\R\rightarrow \R\colon (r,\dot r, r', \dot r')\mapsto
\dot r' + r'.$$
   We see that $\zq_1\circ \zk_Z =\zq_1$, $\zq_2\circ \zk_Z =\zq_3$ and,
consequently $\zk_Z$ projects to diffeomorphisms
   $$\widetilde{\zk}_Z\colon \tT\tT Z\rightarrow \tT\tT Z, \ \text{and} \
\overline{\zk}_Z\colon \bT\tT Z \rightarrow \tT\bT Z.$$
   The affine structures of $\tT\tT Z$, $\tT\bT Z$ and $\bT\tT Z$ are obtained
by reductions from the vector bundle structures of $\sT\sT Z$. It follows that
$\widetilde{\zk}_Z$ and $\overline{\zk}_Z$ are isomorphism of corresponding
special affine  bundles.
      \endproof

      \claim \c{t}{Theorem}{}                           \label{Ca17}
   There are canonical isomorphisms
    $$ \widetilde{\eza}_\bZ \colon \tT \tT^\# \bZ\rightarrow \tT^\#\tT \bZ \
\ \text{dual to} \ \ \widetilde{\zk}_\bZ
                                                \tag \label{Fa98}$$
\noindent such that the diagram $$ \CD  \tT\tT^\# Z @()
\L{\widetilde{\eza}_\bZ} @(1,0) @() \L{\tT\widetilde{\ezp}_Z }@(0,-1)&
\tT^\#\tT Z @()\L{\widetilde{\ezp}_{\tT Z}} @(0,-1) \\
   \sT M @() \L{\text{id}} @(1,0) & \sT M
         \endCD                  \tag \label{Fa98a}$$
is commutative and
      $$ \overline{\eza}_\bZ \colon \tT\sC \bZ \rightarrow \sC\tT \bZ
                                                \tag \label{Fa99}$$
    with the commutative diagram
$$ \CD  \tT\sC Z @() \L{\overline{\eza}_Z} @(1,0) @()
\L{\tT \sC\zz}@(0,-1)& \sC\tT Z @()\L{\sC\tT\zz} @(0,-1) \\
   \sT M @() \L{\text{id}} @(1,0) & \sT M
         \endCD .                 \tag \label{Fa99a}$$
      \endclaim
      \proof
   Follows directly from  the previous theorem and from~\Ref{Fa96}
      \endproof
   Both isomorphisms, $\widetilde{\eza}_\bZ$ and $\overline{\eza}_\bZ$ project to the
same diffeomorphism
      $$ \underline{\eza}_\bZ \colon \sT \sP Z \rightarrow \sP\tT Z.
                                                \tag \label{Fa100}$$
   Let  $\zf\colon M\rightarrow \sP Z$ be an affine 1-form, then
$\underline{\eza}_\bZ\circ \sT \zf \colon \sT M\rightarrow \sP\tT Z$ is an
affine 1-form on $\sT M$. We call it the {\it complete} or {\it tangent
lift} of $\zf$ and we will denote it $\xd_{\tT}\zf$.  As in the classical
case we have for a section $\zs$ of $\zz$ and for an affine 1-form $\zf$
      $$ \aligned
   \xd_{\tT \bZ}\xd\zs&=\xd \xd_{\tT \bZ}\zs \\
   \xd_\sT \xd \zf&=\xd \xd_{\tT \bZ}\zf.
                              \endaligned         \tag \label{Fa101}$$

      \claim \c{p}{Proposition}{}                           \label{C1}

      $$ \aligned \overline{\eza}_Z^{\textstyle *} \zvy_{\tT \bZ} &=
\xd_{\tT \bZ} \zvy_\bZ\\
         \underline{\eza}_Z ^{\textstyle *} \ezw_{\tT\bZ}&= \xd_\sT \ezw_\bZ.
                        \endaligned               \tag \label{F1}$$
      \endclaim
      \proof
   In a  trivialization given by a section of $\zz$, $\tT \sC Z \simeq
\sT\sT^{\textstyle *} M \times \R$, $\sC \tT Z \simeq  \sT^{\textstyle *} \sT M \times
\R$ and $\widetilde{\eza}_\bZ$ is the trivial lift of $\eza_M$. The first equality
follows from the well known equality $\eza_M^{\textstyle *} \zvy_{\sT M} = \dt
\zvy_M$. The second equality is a consequence of the first one and of \Ref{Fa101}.
      \endproof

        \hl1{The dynamics of a charged particle [\Ref{TU}]}
    In gauge theories potentials are interpreted as connections on
principal bundles. In the electrodynamics  the gauge group is $(\R, +)$
and a potential is a connection on a principal bundle $\zz\colon
Z\rightarrow M$ over the space-time $M$. The bundle $Z$ can be considered
a special affine bundle modeled on $M\times \R$ with the distnguished
section $\bold 1_M$. An electromagnetic potential is a section $A\colon
M\rightarrow \sP Z$.

    According to [\Ref{We}] the phase manifold for a particle with the
charge $e\in \R$ is obtained by the symplectic reduction of
$\sT^{\textstyle *}  Z$ with respect to the coisotropic submanifold
        $$  K_e = \{p\in \sT^{\textstyle *} Z\colon \langle p, X_1\rangle
=e\}.
                                                    \tag \label{Fa106}$$
    Let us denote by $\sP_e Z$ the reduced phase space. It is easy to see
that it is an affine bundle modeled on $\sT ^{\textstyle *} M$. We show that
$\sP_e Z$ is the phase bundle for certain special affine bundle $\bZ_e$.

    First, let $Y=Z\times \R$ and $\bY = (Y,\bold 1_Z)$ as in
Example~\Ref{Ca10}. We define an $\R$-action on $Y$ by the formula
        $$ (Z\times \R)\times \R\ni ((z,r),t) \mapsto (z+t, r +te)\in
Z\times \R =Y.
                                                            \tag \label{Fa107}$$
     The space of orbits is an affine bundle modeled on $M\times \R$ and
denoted by $Z_e$. We denote by $\zz_e$ the canonical projection
$Z_e\rightarrow M$. The distinguished section of $\sV(Y)$ (the function
$\bold 1_Z$) projects to the constant function $\bold 1_M$ and the
canonical projection $\zl_e\colon Y\rightarrow Z_e$ is a morphism of
special affine bundles $\bY\rightarrow \bZ_e = (Z_e, \bold 1_M)$. The
induced $R$-action on $Z_e$ has the form
        $$ \zl_e(z,r) +s = \zl_e(z,r+s) = \zl_e(z+t, r+ s +te)
                                                        \tag \label{Fa107a}$$
    For $e=0$ the bundle $Z_e$ is trivial: $Z_0 =M\times \R$ and for
$e\neq 0$ we have a diffeomorphism
        $$\align
    \zF_e&\colon Z\rightarrow Z_e  \\
    &\colon z \mapsto \zl_e(z,0).
                                                    \tag \label{Fa108}\endalign$$
The diffeomorphism $\zF_e$ is not a morphism of special affine bundles:
        $$ \zF_e(z+r) = \zl_e(z+r,0) = \zl_e(z,- er)= \zl_e(z,0)-
er = \zF_e(z) -er,
                                                    \tag \label{Fa109}$$
    except for $e =-1$ (Example~\Ref{Ca10}).

    Let $\zs$ be a section of $\zz_e$. The function $\zl_e^{\textstyle *}
\zs$ on $Z$ has the property
        $$  X_1(\zl_e^{\textstyle *} \zs) = e
                                                            \tag \label{Fa110}$$

     We conclude that the induced by $\zl_e$ relation $\sP Y\rightarrow
\sP Z_e$ is the symplectic reduction with respect to a coisotropic submanifold
        $$ K_e = \{p\in \sT^{\textstyle *} Z\colon \langle p, X_1\rangle
=e\}.
                                                            \tag \label{Fa111}$$
    Thus we have proved.
        \claim \c{p}{Proposition}{}         \label{Ca30}
    The phase manifold $P_e Z$ for a particle with the charge $e$ is the
phase bundle for the special affine bundle $\bZ_e$.
        \endclaim

    The diffeomorphism $\zF_e$ gives a one-to-one correspondence between
sections of $\zz$ and sections  of $\zz_e$, for $e\neq 0$. It follows that
a choosen section of $\zz$ provides a trivialization of $Z$ and also of
$Z_e$. In such trivializations, a section $\zs$ of $\zz$ and the
corresponding section $\zF_e\circ \zs$ of $\zz_e$ are functions on $M$ related by
the formula
        $$ \zF_e\circ \zs (m) = -e\zs(m).
                                                            \tag \label{Fa112}$$
     The correspondence $\zs\rightarrow \zF_e\circ \zs$ of sections
projects to a correspondence of affine covectors and consequently gives a
correspondence of affine 1-forms. Let $A$ be a section of $\sP \zz$ and
$A_e$ the corresponding  section $\sP\zz_e$. In a given trivializations, the
sections $A$ and $A_e$ are 1-forms related by the formula
        $$ A_e = -eA.
                                                            \tag \label{Fa113}$$

    The Lagrangian of a charged particle is a section $L_e$ of the bundle
$\tT\zz_a\colon \tT Z_e\rightarrow \sT M$ over the open set $C=\{v\in \sT
M\colon g(v,v)>0\}$ and  given by the formula
        $$ L(v)= \langle A_e,v\rangle + \slanted  m\sqrt{g(v,v)}
                                                    \tag \label{Fa113a}$$
     where $g$ is the metric tensor, $\slanted m$ is the mass of the
particle, and $\langle \,,\,\rangle $ has been defined in~\Ref{Fa42}. The
Lagrangian section $L_e$ generates a Lagrangian submanifold $D_{l,e}$ of
$\sP Z_e$.

The dynamics of the system is the Lagrangian submanifold $D_e$ of $\sT\sP
Z_e$
        $$ D_e = \underline{\eza}^{-1}_{\bZ_e} (D_{l,e}).
                                                            \tag \label{Fa114}$$
    We have
        $$ D_e = \{w\in \sT \sP Z_e \colon v=\sT\sP\zz_{Z_e}(w)\in C \
\text{and} \ \slanted m v = \sqrt{g(v,v)}g(p - A_e(m))\},
                                                            \tag \label{Fa115}$$
     where $p = \ezt_{\sP Z_e}(w)$ and $m = \sP\zz_e(p)$.

In the Hamiltonian formulation of the dynamics, the generating object is a
Morse family $F\colon N\rightarrow \R$ defined on $N =\sT M\times _M\sP Z$
with the canonical projection
        $$ \zvs \colon N= \sT M\times _M\sP Z \rightarrow \sP Z
                                                            $$
    by
        $$ F(v, p) = L(v) - \langle p,v\rangle .
                                                            $$
\newcounter\bibno
\def\bb{\bib\key \bibno}
        \makebib
        \bb \label{Be} \by  S. Benenti \paper Fibr\'es affines canoniques
et m\'ecanique newtonienne \jour S\'eminaire Sud-Rhodanien de G\'eometrie,
Journ\'ees S.M.F \paperinfo Lyon, 26-30 mai, 1986
        \endbib
    \bb \label{GU1} \by J.~Grabowski and P.~Urba\'nski \pp 6743--6777
    \paper Tangent lifts of Poisson and related structures
    \yr  1995 \vol 28
    \jour J. Phys. A
        \endbib

        \bb \label{Ko} \by K. Konieczna \paper Affine formulation of the
Newtonian analytical mechanics (in Polish)\paperinfo Thesis, University of
Warsaw 1995
        \endbib
    \bb \label{KU} \by K.~Konieczna and P.~Urba\'nski \paper Double vector
bundles and duality \jour Archivum Mathematicum \vol 35 \yr 1999 \pp 59--95
        \endbib
       \bb \label{Sa2} \by Mart\'\i nez, E.,  Mestdag, T.,  and Sarlet,W.,
\paper
Lie algebroid structures and Lagrangian systems on affine bundles \paperinfo
math.DG/0203178
                \endbib
        \bb \label{MP}  \by E.~Massa, E.~Pagani,  and P.~Lorenzoni  \paper
On the gauge structure of classical mechanics  \pp 69--91  \jour Transport
Theory Stat. Phys.  \yr 2000 \vol 29
        \endbib
        \bb \label{MV}  \by E.~Massa, S.~Vignolo, and D.~Bruno  \pp
6713--6742    \paper Non-holonomic Lagrangian and Hamiltonian mechanics:
an intrinsic approach \jour J. Phys. A \yr 2002
        \endbib
        \bb \label{MT} \by M.~R.~Menzio and W.~Tulczyjew \pp 349--367
\paper Infinitesimal symplectic relations and generalized Hamiltonian
dynamics \yr 1978 \vol 28 \jour Ann. Inst. H. Poincar\'e.
        \endbib
    \bb \label{PT} \by G.~Pidello and W.~Tulczyjew \pp 249--265 \paper
Derivations of differential forms on jet bundles \yr 1987 \vol 147 \jour
Ann. Mat. Pura Appl.
        \endbib

        \bb \label{Sa} \by Sarlet,W.,  Mestdag, T.,  and Mart\'\i nez, E.
\paper
Lie algebroid structures on a class of affine bundles \paperinfo
math.DG/0201264
                \endbib
         \bb \label{Sa1} \by Sarlet, W., Mestdag, T., and Mart\'\i nez, E.
\paper Lagrangian equations on affine Lie algebroids \paperinfo In: {\it
Differential Geometry and its Applications}, Proc. 8th Int. Conf. (Opava
2001), D.~Krupka et al. (Eds.) \toappear
                \endbib
    \bb \label{Tu} \by W.~Tulczyjew \pp 101--114 \paper Hamiltonian systems,
Lagrangian systems, and the Legendre transformation \yr 1974 \vol 14 \jour
Symposia Math.
        \endbib

        \bb \label{Tu1}
\by Tulczyjew, W.M.
\paper Frame independence of analytical mechanics
\jour Atti Accad. Sci. Torino,
\vol 119
\yr 1985
        \endbib
        \bb \label{TU} \by Tulczyjew, W.~M. and Urba\'nski, P.  \pp 257--265
\paper
An affine framework for the dynamics of char\-ged particles \yr 1992 \vol
126 \jour Atti Accad. Sci. Torino Suppl. n. 2
        \endbib
        \bb \label{TUZ} \by Tulczyjew, W.~M., Urba\'nski, P., and
S.~Zakrzewski \pp 66-72 \paper A pseudocategory of principal bundles \yr 1988 \vol
122 \jour Atti Accad. Sci. Torino
        \endbib
        \bb \label{U1} \by Urba\'nski, P. \pp 123--129 \paper Affine Poisson
structure in analytical mechanics \inbook Quantization and
Infinite-Dimensional Systems \eds J-P. Antoine and others
\yr 1994 \publ Plenum Press \publaddr New York and London
        \endbib
        \bb \label{Vi}  \by Vignolo, S.  \pp 5117--5135  \paper A new
presymplectic geometrical framework for time-dependent Lagrangian systems
\jour J.~Phys.~A.  \yr 2000 \vol 33
        \endbib
        \bb \label{We} \by Weinstein, A. \paper A universal phase space for
particles in Yang-Mills fields \jour Lett. Math. Phys. \yr 1978\pp 417--420
\vol 2
        \endbib

\enddocument